\documentclass[11pt]{article}
\usepackage[applemac]{inputenc}% Pour saisie sous Mac 
\usepackage{geometry}
\usepackage[square,authoryear]{natbib}
\usepackage[colorlinks]{hyperref}
\usepackage[T1]{fontenc}% Pour utiliser les fontes codees sur 8 bits
\usepackage{latexsym,amsmath,amssymb,amsfonts, amscd,theorem, graphics}% Packages de base ecrits par
\usepackage{array}% Supplement pour les tableaux
\usepackage{graphicx}% Fondamental pour integrer du graphique
\usepackage{float}% Pour placer les figure plus facilement
\usepackage{multicol}% Indispensable pour l'index
\usepackage{fancyhdr}% Permet d'adapter les tÉtes de page
\usepackage[arrow, matrix, curve]{xy}
\usepackage{url}
\usepackage{colordvi}
\usepackage{framed}
\usepackage{color}
\usepackage{eucal}

\usepackage{mathrsfs}
\usepackage{multirow}

\newcommand{\mathsym}[1]{{}}

\newtheorem{theorem}{Theorem}[section]

\newtheorem{lemma}[theorem]{Lemma}
\newtheorem{remark}[theorem]{Remark}
\newtheorem{proposition}[theorem]{Proposition}

\begin{document}
\title{Multisymplectic Lie group variational integrator for a geometrically exact beam in $\mathbb{R}^3$}
\author{Fran\c{c}ois Demoures$^{1}$, Fran\c{c}ois Gay-Balmaz$^{2}$, 
Marin Kobilarov$^{3}$, and Tudor S. Ratiu$^{4}$}
\addtocounter{footnote}{1} 
\footnotetext{Section de Math\'ematiques, 
\'Ecole Polytechnique F\'ed\'erale de
Lausanne, CH-1015 Lausanne. Switzerland. Supported by 
Swiss NSF grant 200020-137704.
\texttt{francois.demoures@epfl.ch}
\addtocounter{footnote}{1} }
\footnotetext{CNRS / LMD, \'Ecole Normale Sup\'erieure, 
Paris, France. Partially supported by a ÊProjet Incitatif de RechercheË contract from ENS-Paris.
\texttt{francois.gay-balmaz@lmd.ens.fr}
\addtocounter{footnote}{1} }
\footnotetext{Laboratory for Computational Sensing and Robotics, 
Johns Hopkins University, 112 Hackerman Hall,
3400 N. Charles Street,
Baltimore, MD 21218, USA. Partially supported by XXXXX. 
\texttt{marin@jhu.edu}
\addtocounter{footnote}{1} }
\footnotetext{Section de Math\'ematiques and Bernoulli
Center, \'Ecole Polytechnique F\'ed\'erale de
Lausanne, CH-1015 Lausanne. Switzerland. Partially supported 
by Swiss NSF grant 200021-140238 and by the government grant 
of the Russian Federation for support of research projects 
implemented by leading scientists, Lomonosov Moscow State 
University under  agreement No. 11.G34.31.0054. 
\texttt{tudor.ratiu@epfl.ch}
\addtocounter{footnote}{1} }

\date{}

\maketitle

\begin{abstract}
In this paper we develop, study, and test a Lie group multisymplectic integrator for geometrically exact beams based on the covariant Lagrangian formulation. We exploit the multisymplectic character of the integrator to analyze the energy and momentum map conservations associated to the temporal and spatial discrete evolutions.
\end{abstract}

\tableofcontents

\section{Introduction}

In this paper we develop, study, and test a Lie group multisymplectic integrator for geometrically exact beams.
Multisymplectic integrators, as developed in \cite{MaPaSh1998}, are 
based on a discrete version of spacetime covariant variational 
principles in field theory (e.g., \cite{GoMa2012}) 
and are extensions of the well-known symplectic integrators for 
Hamiltonian ODEs (\cite{HaLuWa2010}). The geometrically exact 
model for elastic beams, in the spirit of classical mechanics,
has been developed in \cite{Si1985} and \cite{SiMaKr1988}. In 
this paper, we shall employ the field theoretic covariant description 
of geometrically exact beams, as developed in \cite{ElGa-BaHoPuRa2010} 
through covariant Lagrangian reduction. A noteworthy feature of our proposed multisymplectic point of view is that it allows us to describe
not only the behavior of the beam during an interval of time, which is a classical dynamical point of view, but also to discover the evolution 
in space of the deformations of the beam when the evolution ``in time'' of the strain located at a boundary node is known.

At the core of our approach lies a discrete variational principle 
in convective representation defined directly on the configuration 
manifold without resorting to local coordinates. This is accomplished 
by exploiting its Lie group structure and regarding velocities as 
elements of its Lie algebra, a well established approach in geometric 
integration (e.g., \cite{IsMu-KaNoZa2000}) and discrete optimal 
control on Lie groups (e.g., \cite{KoMa2011}). A central point in 
our development is to perform a unified Lie group discretization, both spatially and temporally, in a geometrically consistent manner. The discrete equations of motion then take a surprisingly simple-to-implement form using retraction maps, such as the Cayley map and its derivatives.
The advantages of Lie group formulations have been explored in a number of recent works in multibody dynamics, such as 
\cite{MuTe2009}, \cite{PaBoPl1995}, \cite{PaCh2005}, \cite{MuMa2003}. 
A quaternion-based formulation has been employed for geometrically 
exact models in \cite{CeSa2010}. Nonlocal geometrically exact models (charged molecular strands) have been developed in 
\cite{ElGa-BaHoPuRa2010}. 

A family of Lie group time integrators for the simulation of flexible multibody systems has been proposed in \cite{bruehls10}, while a 
Lie group extension of the generalized-$\alpha$-time integration 
method for the simulation of flexible multibody systems was developed 
in \cite{bruehls12}. Regarding the control theory perspective, in \cite{SoBr2013} it is shown that the Lie group setting allows an efficient development and implementation of semi-analytical methods for sensitivity analysis. Directly associated to the geometrically exact 
beam model, two different structure preserving integrators are derived
in \cite{DeGBLeOBRaWe2013}, using a Lie group time variational integrator with favorable comparisons to energy-momentum schemes.
 
In contrast to these previous approaches, our key contribution is to exploit the multisymplectic point of view to develop a new family of structure-preserving algorithms for the geometrically exact beam model.
 We investigate the quality of the resulting methods both analytically and numerically through the evolution of the associated discrete energy and discrete momentum maps. In addition, we consider the role of symplecticity associated to the displacement in both time and space. 
In particular, we highlight, through two examples, the relationship between the discrete covariant Noether theorem associated to time
evolution and the discrete Noether theorem for space evolution.

The present paper builds upon the discrete multisymplectic variational theory developed in \cite{DeGBRa2013} which is based on a discretization 
of the configuration bundle, the jet bundle, the density Lagrangian, and the variational principle, following \cite{MaPaSh1998}. In particular, the paper \cite{DeGBRa2013} studies the link between discrete multisymplecticity and usual 
symplecticity, the relationship between the discrete covariant Noether theorem and the discrete 
standard Noether theorem (in the Lagrangian formulation), and the role of the boundary conditions.

\medskip

The paper is organized as follows. In Section \ref{Sec2}, we briefly 
review the covariant continuum formulation of the geometrically exact beam model. The discrete problem is formulated in Section \ref{MVI}. Spacetime discretization and the discrete Lagrangian are introduced, the discrete covariant principle is stated, and the integrator is obtained. The discrete covariant formulation of the Lagrange-d'Alembert principle 
with forcing is recalled. In Section \ref{sec4}, Noether's theorem 
and multisymplectic discrete covariant Euler-Lagrange equations are 
developed. We recall the relationship between the symplectic nature of 
the variational integrator for time and space evolution from the 
point of view of multisymplectic geometry. The main accomplishments of this paper are illustrated by the results of our tests in Section \ref{sec5} which illustrate our new point of view by presenting 
time and space simulations using multisymplectic integrators.

\section{Covariant (spacetime) formulation of geometrically exact beams}\label{Sec2} 

\paragraph{Geometrically exact beams.} Our developments are based on the geometrically exact beam model as developed in \cite{Re1972}, \cite{Si1985}, and \cite{SiMaKr1988}. This model, regarded as an extension of the classical Kirchhoff-Love rod model (see \cite{An1974}), provides a convenient parametrization from an analytical and computational point of view. In geometrically exact models, the instantaneous configuration of a beam is described by its line of centroids, as a map $\mathbf{r} :[0,L]\rightarrow\mathbb{R}^3$, and the orientation of all its cross-sections at points $\mathbf{r}(s)$, $s \in [0,L]$, by a moving orthonormal basis $\{\textbf{d}_1(s), \textbf{d}_2(s), \textbf{d}_3(s)\}$. The attitude of this moving basis is described by a map $\Lambda:[0,L]\rightarrow SO(3)$ satisfying
\begin{equation}
\label{d_definition}
\mathbf{d}_I(s)=\Lambda(s)\mathbf{E}_I,\qquad  I=1,2,3,
\end{equation}
where $\{\textbf{E}_1, \textbf{E}_2, \textbf{E}_3\}$ is a fixed orthonormal basis, 
the \textit{material frame}.

The motion of the beam is thus described by the the configuration 
variables $ \Lambda (t,s)$ and $\mathbf{r} (t,s)$, solutions of a 
critical action principle associated to the Lagrangian of the beam.
Two mathematical interpretations can be made in the variational principle. 
First, one can view these configuration variables as curves 
$t \mapsto ( \Lambda (t), \mathbf{r}(t))$ in the infinite dimensional 
space $ \mathcal{F} (\mathcal{B} ,\mathcal{M} )$ of maps from 
$\mathcal{B} :=[0,L]$ to $\mathcal{M} :=SO(3) \times \mathbb{R}^3$. 
This approach is referred to as the dynamic formulation. Secondly, 
one can view the configuration variables as spacetime maps 
$(s,t) \mapsto  (\Lambda(t,s), \mathbf{r}(t,s))$ from  
$X:=[0,T] \times [0,L]$ to $\mathcal{M} =SO(3) \times \mathbb{R}^3$. 
This approach is referred to as the covariant formulation.
We quickly comment on these two approaches below.
Note that we have introduced above the notations $\mathcal{B} $, 
$X$, $\mathcal{M}$ for the spatial domain, the spacetime, and 
the space of all possible deformations, respectively.

\paragraph{Dynamic formulation.} In the traditional Lagrangian formulation of continuum mechanics, the motion of the mechanical 
system is described by a time-dependent curve $q(t)$ in the (infinite dimensional)
configuration space $Q= \mathcal{F} (\mathcal{B} ,\mathcal{M} )$ of the system. The Lagrangian function 
is a given map $L:TQ \rightarrow \mathbb{R}$ defined on the 
tangent bundle $TQ$ of $Q$, to which one associates the action 
functional
\[
\mathcal{A} (q(\cdot))=\int_0^TL(q(t), \dot{q}(t)) \operatorname{d}\!t. 
\]
Hamilton's Principle $ \delta \mathcal{A} =0$ for variations $\delta q$ vanishing at the endpoints $t=0,T$ yields the classical Euler-Lagrange equations
\[
\frac{d}{dt} \frac{\partial L}{\partial \dot{q}}- 
\frac{\partial L}{\partial q}=0.
\]

When the Lagrangian admits Lie group symmetries then, by Noether's theorem, the associated momentum map is conserved. More precisely, consider the action of a Lie group $G$ with Lie algebra $ \mathfrak{g}  $ and dual $ \mathfrak{g}  ^\ast $ on the configuration manifold $Q$. If the Lagrangian is invariant under the action of a Lie group $G$ on $Q$, then the momentum map
\begin{equation}\label{Lagr_momap} 
  \mathbf{J}_L:TQ \rightarrow \mathfrak{g}^\ast, \qquad \left\langle \mathbf{J}_L( q, \dot q), \xi \right\rangle = \left\langle \mathbb{F}  L(q, \dot q), \xi _Q (q) \right\rangle 
\end{equation} 
is a conserved quantity along the solutions of the Euler-Lagrange equations. The inner product $\left\langle \mu , \xi \right\rangle$ denotes the pairing between $\mu \in \mathfrak{g}^\ast$ and $\xi \in \mathfrak{g}$. The term $\xi_Q $ denotes the infinitesimal generator of the action associated to the Lie algebra element $\xi \in \mathfrak{g}$,
\[
\xi _Q (q) := \left.\frac{d}{d\varepsilon}\right|_{\varepsilon=0} \operatorname{exp}( \varepsilon \xi )q,
\]
where $ \operatorname{exp}: \mathfrak{g}  \rightarrow G$ is the Lie group exponential map and $gq$ the group action of $ g \in G$ on $ q \in Q$. The map $ \mathbb{F}  L:TQ \rightarrow T ^\ast Q$ denotes the Legendre transform associated to $L$ which, in standard tangent
bundle coordinates, has the expression 
$\mathbb{F}L(q, \dot{q})= \left(q,\frac{\partial L}{\partial \dot{q}}\right) \in T^*_qQ$.

\medskip

In the case of the geometrically exact beam we have $Q=\mathcal{F} \left([0,L],SO(3) \times \mathbb{R}^3\right)$ and the Lagrangian $L:TQ\rightarrow \mathbb{R}$ reads
\begin{align}\label{continuous_Lagrangian_L0}
L(\Lambda, \mathbf{r}, \dot{\Lambda}, \dot{\mathbf{r}}) :=&   \frac{1}{2} \int_0^L \left[ M \left\| \boldsymbol{\gamma}   \right\|^2 + \boldsymbol{\omega}  ^T J \boldsymbol{\omega}  \right]  \operatorname{d}\!s \nonumber\\
&  -
\frac{1}{2} \int_0^L \left[ (\boldsymbol{\Gamma}  - \mathbf{E}_3)^T \mathbf{ C}_1 (\boldsymbol{\Gamma}  - \mathbf{ E}_3) +  \boldsymbol{\Omega}  ^T \mathbf{ C}_2\boldsymbol{\Omega}    \right] \operatorname{d}\!s  -  \int_0^L \Pi( \Lambda , \mathbf{r} )\! \operatorname{d}\!s ,
\end{align}
where we defined the convective velocities and strains by
\begin{equation}
\label{omega_gamma_def}
( \omega , \boldsymbol{\gamma}  ):= ( \Lambda ^{-1} \dot \Lambda , \Lambda ^{-1} \dot{\mathbf{r}} ), \qquad ( \Omega , \boldsymbol{\Gamma} ):= ( \Lambda ^{-1} \Lambda ', \Lambda ^{-1} \mathbf{r}'),
\end{equation}
considering that the thickness of the rod is small compared to its length, and that the material is homogeneous and isotropic. Here $M$ is the mass and $J$ is the inertia tensor, both assumed to be constant. The matrices $ \mathbf{C} _1 , \mathbf{C} _2 $ are given by (see \cite{SiVu-Qu1986})
\begin{equation}\label{C1_C2}
\mathbf{C}_1:=\mathrm{Diag}\left( GA \ GA \ \ EA  \right) \quad\text{and} \ \ \mathbf{C}_2:=\mathrm{Diag}\left(  EI_1 \ \ EI_2 \ \ GI \right),
\end{equation}
where $A$ is the cross-sectional area of the rod, $I_1$ and $I_2$ are the principal moments of inertia of the cross-section, $I=I_1+I_2$ is its polar moment of inertia, $E$ is Young's modulus, $G= E/[2(1+ \nu)]$ is the shear modulus, and $\nu$ is Poisson's ratio. The basic kinematic assumption of this model precludes changes in the cross-sectional area. Thus, for a given homogenous material, $\mathbf{C}_1$ and $\mathbf{C}_2$ are constant. In \eqref{continuous_Lagrangian_L0}, the three integrals correspond, respectively, to the kinetic energy, the bending energy, and the potential energy density due to the gravitational forces.

\paragraph{Covariant formulation.} In the covariant approach of continuum mechanics, one interprets the configuration variables as space-time dependent maps (or fields) $(t,s) \in [0,T] \times \mathcal{B} \mapsto \varphi (s,t) \in  \mathcal{M} $.

In this framework, the time and space variables are treated in the same way and we shall take advantage of this fact later when formulating the geometric discretization of the beam.
In order to obtain the intrinsic geometric description, the fields have to be interpreted as sections of the (here trivial) fiber bundle $ \pi : X \times \mathcal{M} \rightarrow X$, with $X=[0,T] \times \mathcal{B} $. The action functional is obtained by spacetime integration of the Lagrangian density $ \mathcal{L} $, i.e.,
\begin{equation} \label{action_functional}
\mathcal{A} ( \varphi (\cdot))=\int_X \mathcal{L} ( t,s, \partial _t \varphi , \partial _s \varphi ).
\end{equation}
The Covariant Hamilton Principle $ \delta \mathcal{A} =0$, for 
variations $\delta \varphi$ vanishing at the boundary, yields the covariant Euler-Lagrange equations
\[
\partial _t  \frac{\partial \mathcal{L} }{\partial (\partial_t  \varphi )}+\partial _s  \frac{\partial \mathcal{L} }{\partial (\partial _s\varphi)}- \frac{\partial      \mathcal{L} }{\partial\varphi }=0.
\] 
If the Lie group $G$ is a symmetry of the Lagrangian,
the corresponding covariant momentum map is a conserved quantity. In continuum solid mechanics problems there are two main symmetries: translation and rotation. For example, if $\mathcal{M} = \mathbb{R}^3$, 
the covariant linear and angular momentum maps are given as follows.

\textit{Linear momentum}: The action on $X\times \mathcal{M} $ by translation of 
$\mathrm{x}\in G=\mathbb{R}^3$ is given by $\mathrm{x}\cdot (s,t,\varphi) =(s,t,\varphi +\mathrm{x})$. For a given direction $\xi \in \mathfrak{g} = \mathbb{R}^3$, the covariant linear momentum map is 
\[ 
\mathfrak{T}_{\mathcal{L}}^\mu(\xi) = \frac{\partial \mathcal{L}}{\partial (\partial_\mu\varphi)^i} \xi^i, \quad \mu \in \{t,s\}, \ i \in \{1,2,3\}.
\]
\textit{Angular momentum}: The proper rotation group $G=SO(3)$ acts
on $X\times \mathcal{M} $ by $\Lambda\cdot (s,t,\varphi) =(s,t, \Lambda \varphi)$, 
$\Lambda \in SO(3)$, For a given direction $\widehat{\boldsymbol{\omega}} \in \mathfrak{g} = \mathfrak{so}(3)$ the covariant angular momentum map is
\[ 
\mathfrak{R}_{\mathcal{L}}^\mu(\widehat{\omega}) = \frac{\partial \mathcal{L}}{\partial (\partial_\mu\varphi)^i}\,\widehat{\boldsymbol{\omega}}_{\ j}^i \varphi^j, \quad \mu \in \{t,s\}, \ i,j \in \{1,2,3\}.
\]

The convective covariant formulation of geometrically exact beams has been developed in \cite{ElGa-BaHoPuRa2010}; see especially \S6 and
\S7 of this paper. In this approach, the maps $ \Lambda , \mathbf{r} $ 
are interpreted as space-time dependent fields
\[
(t,s) \in \mathbb{R}  \times [0,L] \longmapsto( \Lambda(t,s) , \mathbf{r}(t,s) )\in SO(3) \times \mathbb{R}  ^3
\]
 rather than time-dependent curves in the infinite dimensional configuration space $Q= \mathcal{F} \left([0,L],SO(3) 
 \times \mathbb{R}^3\right)$.
The fiber bundle of the problem is therefore given by
\[
X \times G \rightarrow X, \quad\text{with}\quad \quad X= \mathbb{R}  \times [0,L] \ni (t,s), \quad G=SE(3)\ni ( \Lambda , \mathbf{r} ),
\] 
and the Lagrangian density depends on $(\Lambda , \mathbf{r}, \dot{\Lambda }, \dot{\mathbf{r}}, \Lambda ', \mathbf{r}')$, where  
$\dot{} = \partial_t$ and ${}' = \partial_s$. 
Here, $SE(3)$ denotes the special Euclidean group of orientation
preserving rotations and translations and $\mathfrak{se}(3)$ is its
Lie algebra. In terms 
of the convective variables, the Lagrangian density (i.e., the integrand in \eqref{continuous_Lagrangian_L0}) can be written as
\begin{align} 
\label{cal_L_continuous}
 \mathcal{L}( \Lambda , \mathbf{r} ,\boldsymbol{\omega}  , \boldsymbol{\gamma}  , \boldsymbol{\Omega} , \boldsymbol{\Gamma} )&= \frac{1}{2}\langle \mathbb{J}\xi, \xi\rangle-\frac{1}{2} \langle \mathbb{C}\ (\eta - \mathbf{E}_6), 
 (\eta-\mathbf{E}_6) \rangle- \Pi(g) \nonumber\\
&= :K( \xi )- \Phi ( \eta )- \Pi(g)= \mathcal{L} (g, \xi , \eta ),
\end{align} 
where $g:= (\Lambda, \mathbf{r}) \in SE(3)$ are the configuration variables, $\xi :=(\omega, \gamma):= g^{-1}\dot{g}\in \mathfrak{se}(3)$,
are the convective velocities, $\eta :=(\Omega, \Gamma)= 
g^{-1}g'\in \mathfrak{se}(3)$ are the convective strains,  
$\mathbf{E}_6=(0,0,0,0,0,1) \in \mathbb{R}^6$, $\mathbb{J}$ is given
in \eqref{inertia_matrix}, and $\mathbb{C}$ by \eqref{strain_matrix}. Recall that $K$, $ \Phi $, and $ \Pi $ correspond, respectively, to the kinetic energy density, the bending energy density, and the potential energy density. The bold face
letters are the images of light faced letters under the standard
isomorphism $\mathfrak{se}(3) \cong \mathbb{R}^6$ given by
\[
\mathfrak{se}(3)= \mathfrak{so}(3) \times \mathbb{R}^3\ni (\omega, \gamma) \longmapsto (\boldsymbol{\omega},
\boldsymbol{\gamma})\in \mathbb{R}^6, \qquad \widehat{\boldsymbol{\omega}} =\omega, \;\; \boldsymbol{\gamma} = 
\gamma, 
\]
where $\widehat{\boldsymbol{\omega}} \mathbf{v}: = \boldsymbol{\omega}\times \mathbf{v}$ for any $\mathbf{v}\in\mathbb{R}^3$. In the 
expression above, we have used this isomorphism when we wrote
$\eta - \mathbf{E}_6$ to mean $\eta - ({\bf 0}, (0,0,1)) \in 
\mathfrak{so}(3) \times \mathbb{R}^3$. The dual $\mathfrak{se}(3)^*$  is identified
with $\mathbb{R}^6$ by using
$ \left\langle(\boldsymbol{\mu}, \boldsymbol{\nu}), (\omega, \gamma) 
\right\rangle =\boldsymbol{\mu}\cdot \boldsymbol{\omega} + 
\boldsymbol{\nu} \cdot \boldsymbol{\gamma}$ for any $(\omega, \gamma)
\in \mathfrak{se}(3)$.

Using the above isomorphisms of 
$\mathfrak{se}(3)^*$ and $\mathfrak{se}(3)$ with $\mathbb{R}^6$, the
map $\mathbb{J}:\mathfrak{se}(3) \rightarrow\mathfrak{se}(3) ^*$ is the linear operator on $\mathbb{R}^6$ with matrix in the standard basis
equal to 
\begin{equation}\label{inertia_matrix}
\mathbb{J} = \left[\begin{array}{cc}J & 0 \\ 0 & M
    \mathbf I_3 \end{array}\right],
\end{equation}
where $M \in \mathbb{R}$ is the mass by unit of length of the beam, $\mathbf{I}_3$ is
the identity $3 \times 3$ matrix,  and $J$ is inertia tensor.
Similarly, the linear operator $\mathbb C:\mathfrak{se}(3) \rightarrow \mathfrak{se}(3) ^*$ encodes the potential interaction which,
under the isomorphisms of $\mathfrak{se}(3)^*$ and  
$\mathfrak{se}(3)$ with $\mathbb{R}^6$, has matrix in the standard
basis equal to
\begin{equation}\label{strain_matrix}
\mathbb{C}\ = \left[\begin{array}{cc}\mathbf C_2 & 0 \\ 0 & \mathbf C_1 \end{array}\right],
\end{equation}
where $\mathbf{C}_1, \mathbf{C}_2$ are defined in \eqref{C1_C2}.

The Covariant Hamilton Principle becomes in this case
\begin{equation}\label{trivialized_HP}
 \delta \int_0^T \!\!\!\int_0^L \left( K(\xi) - \Phi(\eta)  - \Pi (g) \right) \operatorname{d}\!s \operatorname{d}\!t =0,
\end{equation}
for all variations $ \delta g$ of $g$, vanishing at the boundary.
It yields the trivialized covariant Euler-Lagrange equations,
\begin{equation}\label{Triv_CEL_BEam} 
\frac{d}{dt}\frac{\partial K}{\partial \xi }  - \mathrm{ad}_\xi^* \frac{\partial K}{\partial \xi }= \frac{d}{ds}\frac{\partial \Phi }{\partial \eta}- \mathrm{ad}_\eta^*\frac{\partial \Phi }{\partial \eta}   - g^{-1} \frac{\partial \Pi}{\partial g},
\end{equation}
where $\operatorname{ad}^*_ \xi : \mathfrak{g}^\ast \rightarrow 
\mathfrak{g}^\ast$ is the dual  map to $\operatorname{ad}_\xi : 
\mathfrak{g} \rightarrow \mathfrak{g}$, $\operatorname{ad}_\xi \eta 
:=[\xi, \eta]$. In the case of $G=SE(3)$, we have 
$\operatorname{ad}^*_{(\boldsymbol{\omega}, \boldsymbol{\gamma})}
(\boldsymbol{\mu}, \boldsymbol{\nu})=-(\boldsymbol{\omega} \times \boldsymbol{\mu} + \boldsymbol{\gamma} \times \boldsymbol{\nu} , \boldsymbol{\omega} \times \boldsymbol{\nu})$.  We refer to 
\cite{ElGa-BaHoPuRa2010} for a detailed derivation of these equations for the beam.

\paragraph{External forces.} External Lagrangian forces are added by using a covariant analogue of the Lagrange-d'Alembert principle. Namely, denoting the force density by $\mathfrak{F}( g, \xi, \eta )\in T_g^*G$ the principle is
\begin{align}\label{eq:LDAvar}
\delta \int_0^T \!\!\!\int_0^L \left( K(\xi) - \Phi(\eta) - \Pi(g)\right)  \operatorname{d}\!s \operatorname{d}\!t+ \int_0^T\!\!\!\int_0^L \mathfrak{F}( g, \xi, \eta ) \delta g \, \operatorname{d}\!s \operatorname{d}\!t = 0,
\end{align}
for all variations $ \delta g$ vanishing at the boundary. The resulting equations correspond to \eqref{Triv_CEL_BEam} with the term  $g ^{-1} \mathfrak{F}( g, \xi, \eta )$ added to the right hand side.

\section{Covariant variational integrator} \label{MVI}

We next develop the discrete variational counterpart of the continuous covariant beam formulation. The first step is to perform a space-time discretization which is accomplished in a unified way by regarding displacements in both space and time as Lie group transformations. A discrete covariant variational principle is then formulation based on this discretization. Finally, a structure-preserving integrator is obtained and the details of its implementation are provided.

\subsection{Space time discretization}

\begin{figure}[ht]
\centering
\includegraphics[width=2.3 in]{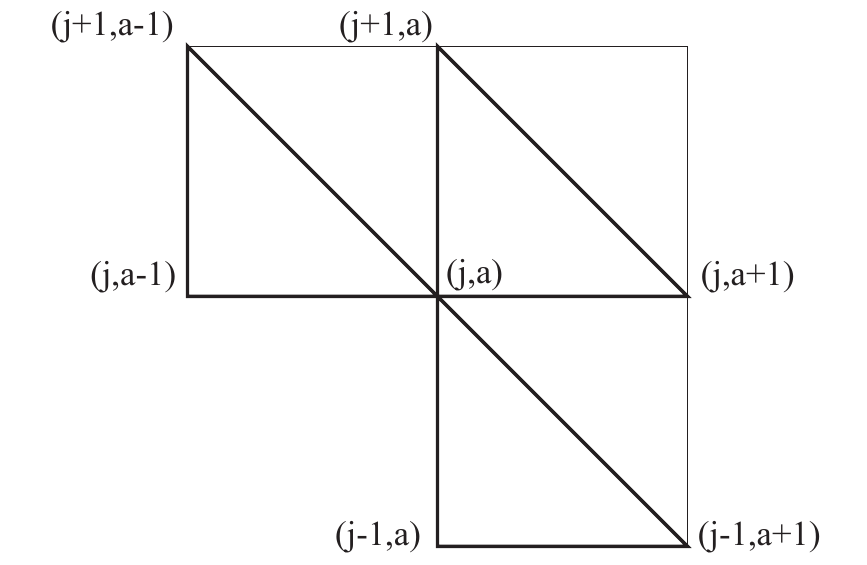}
\vspace{-3pt}
\caption{\footnotesize The triangles $\triangle_a^j, \triangle_a^{j-1}, \triangle_{a-1}^j$. }  
\label{triangles_touching_j_a}  
\end{figure}

Spacetime discretization is realized by fixing a time step $\Delta t$ and a space step $\Delta s$, and decomposing the intervals $[0,T]$ and $[0,L]$  into subintervals of length $\Delta t$ and $\Delta s$, respectively. We denote by  $j\in \{0,...,N\}$ and $a \in \{0,...,A\}$ the time and space indices, respectively.
The discretization of the spacetime domain is based on a triangular decomposition (Figure~\ref{triangles_touching_j_a}), where a \emph{triangle} $\triangle_a^j$ is defined by the three pairs of indices
\[
\triangle_a^j= ((j,a),(j+1,a),(j,a+1)), \quad j=0,...,N-1,\; a=0,...,A-1
\]

Small displacements in both space and time will be represented using Lie algebra elements with the help of a map $\tau : \mathfrak{g}  \rightarrow G$ which is a local diffeomorphism around the origin that satisfies $\tau (0)=e$. The discrete convective velocities $\xi_a^j$ and the discrete convective strains $\eta_a^j$ are defined then by $\tau (\xi^j_a\Delta t)=  (g_a^j)^{-1}g^{j+1}_{a}$ and $\tau (\eta_a^j\Delta s)= (g_a^j)^{-1}g_{a+1}^j$. In our case, $\xi_a^j = (\omega_a^j , \boldsymbol{\gamma}_a^j)$ and 
$\eta_a^j = (\Omega_a^j , \boldsymbol{\Gamma}_a^j)$.

\medskip

The action functional associated to the Lagrangian density \eqref{cal_L_continuous} is approximated on the square $((j,a),(j+1,a),(j,a+1), (j+1, a+1))$ by the discrete Lagrangian $\mathscr{L}_d:X^\triangle _d \times G \times \mathfrak{g} \times \mathfrak{g} \rightarrow \mathbb{R}$
given by 
\begin{equation}\label{DCL_beam}
\mathscr{L}_d(\triangle_a^j, g_a^j, \xi_a^j, \eta_a^{j}) =  \Delta t \Delta s  K(\xi^{j}_a) - \Delta t \Delta s \left[  \Phi(\eta_a^j) +  \Pi ( g _a^j ) \right],
\end{equation}
where $X_d ^\triangle  $ denotes the set of all triangles $\triangle$ in spacetime $[0, T] \times [0, L]$.

\begin{remark}{\rm The domain $X _d ^\triangle \times G \times  \mathfrak{g}  \times \mathfrak{g}  $ of the discrete Lagrangian is understood as the trivialization of the discrete first jet bundle. It is the discrete analogue of the first jet bundle $J ^1 (X \times G)\rightarrow X \times G$, which is the domain of the continuous Lagrangian density. We refer to \cite{MaPeShWe2001} for the detailed description of discrete jet bundles.}
\end{remark}

\subsection{Analogue of the Cayley map for $SE(3)$}
\label{Cayey_section} 

In what follows we shall work
with the Cayley map as an approximation of the exponential map. First, we
briefly recall the Cayley map for the rotation group $SO(3)$. Recall
that the Lie algebras $(\mathfrak{so}(3), [\,,])$ and 
$(\mathbb{R}^3, \times)$ are isomorphic via the map
$\mathfrak{so}(3) \ni \omega \mapsto\boldsymbol{\omega} \in \mathbb{R}^3$ given by $\omega \boldsymbol{\rho} :=    \boldsymbol{\omega}\times \boldsymbol{\rho}$ for all $\boldsymbol{\rho}\in\mathbb{R}^3$.

The classical Cayley map $\operatorname{cay}:\mathfrak{so}(3)\rightarrow SO(3)$ 
is defined by
\begin{equation}
\label{cayley_so_3}
\operatorname{cay}(\omega ):=\left(\mathbf{I}_3 - \frac{\omega }{2}\right)^{-1} 
\left(\mathbf{I}_3 + \frac{\omega }{2} \right)
=\mathbf I_{3}+\frac{4}{4+\|\boldsymbol{\omega} \|^{2}}\left(\omega +\frac{\omega ^{2}}{2}\right).
\end{equation}
The Cayley map is invertible on the set of all matrices
$\Lambda \in SO(3)$ which are not rotations by the angle $\pm \pi$; the
formula for the inverse is (\cite[formulas (10), (11)]{Selig2007}
\begin{equation}
\label{inverse_cayley_so_3}
\operatorname{cay}^{-1}(\Lambda) =\frac{2}{1+\operatorname{Tr}(\Lambda)}
\left(\Lambda - \Lambda^T\right)
=2 (\Lambda - \mathbf{I}_3)(\Lambda + \mathbf{I}_3)^{-1};
\end{equation}
note that $1+\operatorname{Tr}(\Lambda) = 0$ if and only if $\Lambda$
is a rotation by the angle $\pm \pi$.

The Rodrigues formula for the exponential of a matrix $\omega \in 
\mathfrak{so}(3)$ (see, e.g., \cite[formula (9.2.8)]{MaRa1999})
\[
e^\omega = \mathbf{I}_3 + 
\frac{\sin\|\boldsymbol{\omega}\|}{\|\boldsymbol{\omega}\|} \omega + 
\frac{1}{2}\left[\frac{\sin\frac{\|\boldsymbol{\omega}\|}{2}}{
\frac{\|\boldsymbol{\omega}\|}{2}}\right]^2\omega^2
\]
and \eqref{inverse_cayley_so_3} give a precise relation between the exponential and the Cayley maps,
namely
\[
\operatorname{cay}^{-1}\left(e^\omega\right) = 
\frac{2 \sin\|\boldsymbol{\omega}\|}{\|\boldsymbol{\omega}\|(1+ \cos\|\boldsymbol{\omega}\|)} \omega =\frac{\tan 
\frac{\|\boldsymbol{\omega}\|}{2}}{\frac{\|\boldsymbol{\omega}\|}{2}}
\omega
\]
which shows that these two maps are indeed close in a neighborhood of 
the origin of $\mathfrak{so}(3)$.
\medskip

For a general Lie group, given a local diffeomorphism 
$\psi:\mathfrak{g}\rightarrow G$, we denote 
by $\operatorname{d\psi}_\xi : \mathfrak{g}  \rightarrow \mathfrak{g}$ 
the right trivialized (or logarithmic) derivative of $\psi$ at 
$\xi\in \mathfrak{g}$, defined by 
$\operatorname{d\psi}_\xi (\eta)= \left(T_\xi\psi(\eta)\right) 
\psi(\xi)^{-1}$, where $T_\xi\psi: \mathfrak{g}\rightarrow T_{\psi(\xi)}G$ 
is the usual differential (tangent map) of $\psi$ at $\xi\in \mathfrak{g}$.
We compute now the right logarithmic derivative for $\operatorname{cay}:
\mathfrak{so}(3) \rightarrow SO(3)$. For $\omega, \delta\omega \in \mathfrak{so}(3)$ we have
\[
T_\omega\operatorname{cay}(\delta\omega) = 
\left(\mathbf{I}_3 -\frac{\omega}{2} \right)^{-1} 
\frac{\delta \omega}{2}\left(\operatorname{cay}(\omega) + 
\mathbf{I}_3\right)
= \left(\mathbf{I}_3 -\frac{\omega}{2} \right)^{-1} 
\delta\omega\left(\mathbf{I}_3 -\frac{\omega}{2} \right)^{-1} 
\]
and hence $\left(\operatorname{d\,cay}_\omega\right)^{\pm 1}: \mathfrak{so}(3) 
\rightarrow \mathfrak{so}(3)$ have the expressions
\begin{equation}
\label{log_right_der_cay}
\begin{aligned}
\operatorname{d\,cay}_\omega (\delta\omega) &= 
\left(\mathbf{I}_3 -\frac{\omega}{2} \right)^{-1} 
\delta\omega\left(\mathbf{I}_3 +\frac{\omega}{2} \right)^{-1} \\
\left(\operatorname{d\,cay}_\omega\right)^{-1} (\delta\omega) &= 
\left(\mathbf{I}_3 -\frac{\omega}{2} \right) 
\delta\omega\left(\mathbf{I}_3 +\frac{\omega}{2} \right). 
\end{aligned}
\end{equation}
It is useful to regard $\left(\operatorname{d\,cay}_\omega\right)^{\pm 1}: 
\mathbb{R}^3 \rightarrow \mathbb{R}^3$. A lengthy direct computation
using the first formula in \eqref{log_right_der_cay}
proves the first equality below (which recovers the one in 
\cite[Section VI(A)]{KoMa2011}) and the second is an easy verification
from the first:
\begin{equation}
\label{log_right_der_cay_vector}
\begin{aligned}
\operatorname{d\,cay}_\omega &=
\frac{2}{4+\|\boldsymbol{\omega}\|^2} (2\mathbf{I}_3 +\omega )\\
\left(\operatorname{d\,cay}_\omega\right)^{-1}
&=\mathbf{I}_3 - \frac{1}{2}\omega +
  \frac{1}{4}\boldsymbol{\omega} \boldsymbol{\omega} ^T
\end{aligned}
\end{equation}
where $\omega\in \mathfrak{so}(3)$.

\medskip

We need similar formulas for the Lie group $SE(3)$. The computations are simpler, if we embed the special Euclidean group $SE(3) \subset 
SL(4, \mathbb{R})$ and its Lie algebra $\mathfrak{se}(3) \subset 
\mathfrak{sl}(4, \mathbb{R})$ by
\begin{equation}
\label{identifications_se}
SE(3) \ni (\Lambda, \mathbf{r}) \mapsto
\begin{bmatrix}
\Lambda& \mathbf{r}\\
\mathbf{0}^T& 1
\end{bmatrix} \in SL(4, \mathbb{R}), \qquad 
\mathfrak{se}(3)\ni (\omega, \boldsymbol{\gamma}) \mapsto
\begin{bmatrix}
\omega&\boldsymbol{\gamma}\\
\mathbf{0}^T& 0
\end{bmatrix} \in \mathfrak{sl}(4, \mathbb{R}).
\end{equation}
This allows us to work with $SE(3)$ as a group of matrices. The usual
way to define a Cayley map for $SE(3)$ is to imitate the classical
formula, that is, to define the map $\tau:\mathfrak{se}(3) \rightarrow 
SE(3)$ by (see \cite[Section III]{Selig2007})
\begin{align}
\label{cay}
\tau(\omega, \boldsymbol{\gamma}): &= \left(\mathbf{I}_4 - \frac{1}{2}
\begin{bmatrix}
\omega& \boldsymbol{\gamma}\\ {\bf 0}^T&0
\end{bmatrix} \right) ^{-1}
\left(\mathbf{I}_4 + \frac{1}{2}
\begin{bmatrix}
\omega& \boldsymbol{\gamma}\\ {\bf 0}^T&0
\end{bmatrix} \right) = 
\begin{bmatrix}
\operatorname{cay}\omega& \left(\mathbf{I}_3 - 
\frac{\omega}{2}\right)^{-1} \boldsymbol{\gamma}\\
{\bf 0}^T& 1
\end{bmatrix}  \nonumber \\[6pt]
& = \mathbf{I}_4 + 
\begin{bmatrix}
\omega& \boldsymbol{\gamma}\\ {\bf 0}^T&0
\end{bmatrix} 
+\frac{2}{4+\|\boldsymbol{\omega}\|^2}
\begin{bmatrix}
\omega& \boldsymbol{\gamma}\\ {\bf 0}^T&0
\end{bmatrix}^2 +
\frac{1}{4+\|\boldsymbol{\omega}\|^2}
\begin{bmatrix}
\omega& \boldsymbol{\gamma}\\ {\bf 0}^T&0
\end{bmatrix}^3 \nonumber \\[6pt]
& = \begin{bmatrix}
\operatorname{cay}\omega& \frac{4}{4+ \|\boldsymbol{\omega}\|^2} \left(\mathbf{I}_3 +\frac{1}{2} \omega +
\frac{1}{4}\boldsymbol{\omega}\boldsymbol{\omega}^T\right) \boldsymbol{\gamma}\\
{\bf 0}^T& 1
\end{bmatrix};
\end{align}
the third equality is obtained from the formula
\[
\begin{bmatrix}
\omega& \boldsymbol{\gamma}\\ {\bf 0}^T&0
\end{bmatrix}^3 = 
\begin{bmatrix}
- \|\boldsymbol{\omega}\|^2\omega& \omega^2\boldsymbol{\gamma}\\ 
{\bf 0}^T&0
\end{bmatrix}.
\]
The map $\tau$ is invertible on the set of all elements $(\Lambda, \mathbf{r}) \in SE(3)$ for which $\Lambda$ is not a rotation by the
angle $\pm \pi$, namely
\begin{align*}
\tau^{-1}(\Lambda, \mathbf{r}) &= 
\begin{bmatrix}
\operatorname{cay}^{-1}(\Lambda)& 
2(\Lambda + \mathbf{I}_3)^{-1}\mathbf{r}\\
{\bf 0}^T& 0
\end{bmatrix} = 
2\begin{bmatrix}
(\Lambda + \mathbf{I}_3)^{-1}& 0\\
{\bf 0}^T& 1
\end{bmatrix}
\begin{bmatrix}
\Lambda - \mathbf{I}_3& \mathbf{r}\\
{\bf 0}^T& 0
\end{bmatrix} \\
& = -2\left(\mathbf{I}_4 + \left[
 \begin{array}{cc} \Lambda & \mathbf{r} \\0 & 1\end{array}\right]\right)^{-1} \left( \mathbf{I}_4 - \left[
 \begin{array}{cc} \Lambda & \mathbf{r} \\0 & 1\end{array}\right] \right)
\end{align*}
as a direct verification shows (see \cite[formula (21)]{Selig2007}).
\medskip

Finally, we compute the right logarithmic derivative of $\tau$. 
Proceeding as in the case of the Cayley map but using \eqref{cay}, we get
\[
T_{(\omega, \boldsymbol{\gamma})} \tau(\delta\omega, \delta\boldsymbol{\gamma}) = \left(\mathbf{I}_4 - \frac{1}{2}
\begin{bmatrix}
\omega& \boldsymbol{\gamma}\\
{\bf 0}^T & 0
\end{bmatrix}\right)^{-1}
\begin{bmatrix}
\delta\omega& \delta\boldsymbol{\gamma}\\
{\bf 0}^T&0
\end{bmatrix}
\left(\mathbf{I}_4 - \frac{1}{2}
\begin{bmatrix}
\omega&\boldsymbol{\gamma}\\
{\bf 0}^T & 0
\end{bmatrix}\right)^{-1}
\]
and hence
\begin{align*}
\operatorname{d\,\tau}_{(\omega, \boldsymbol{\gamma})}
(\delta\omega, \delta\boldsymbol{\gamma}) &=
\left(\mathbf{I}_4 - \frac{1}{2}
\begin{bmatrix}
\omega& \boldsymbol{\gamma}\\
{\bf 0}^T & 0
\end{bmatrix}\right)^{-1}
\begin{bmatrix}
\delta\omega& \delta \boldsymbol{\gamma}\\
{\bf 0}^T&0
\end{bmatrix}
\left(\mathbf{I}_4 + \frac{1}{2}
\begin{bmatrix}
\omega& \boldsymbol{\gamma} \\ 
{\bf 0}^T & 0
\end{bmatrix}\right)^{-1} \nonumber \\[6pt]
\left(\operatorname{d\,\tau}_{(\omega, \boldsymbol{\gamma})}\right)^{-1}
(\delta\omega, \delta\boldsymbol{\gamma}) &=
\left(\mathbf{I}_4 - \frac{1}{2}
\begin{bmatrix}
\omega& \boldsymbol{\gamma}\\
{\bf 0}^T & 0
\end{bmatrix}\right)
\begin{bmatrix}
\delta\omega& \delta \boldsymbol{\gamma}\\
{\bf 0}^T&0
\end{bmatrix}
\left(\mathbf{I}_4 + \frac{1}{2}
\begin{bmatrix}
\omega& \boldsymbol{\gamma}\\
{\bf 0}^T&0
\end{bmatrix} \right) \\[6pt]
& = \begin{bmatrix}
\left(\operatorname{d\, cay}_\omega \right)^{-1}(\delta\omega) & 
\left(\mathbf{I}_3 - \frac{1}{2}\omega\right)
\left(\frac{1}{2}\delta \omega \boldsymbol{\gamma} + 
\delta \boldsymbol{\gamma} \right)\\[6pt]
{\bf 0}^T& 0
\end{bmatrix}.
\end{align*}
Viewed as an operator $\left(\operatorname{d\,\tau}_{(\omega, \boldsymbol{\gamma})}\right)^{-1}: \mathbb{R}^3 \times \mathbb{R}^3
\rightarrow \mathbb{R}^3 \times \mathbb{R}^3$, the $6\times 6$ matrix
of this linear map has the expression (in agreement with 
\cite[Section VI(C)]{KoMa2011})
\begin{align}
\label{matrix_dcay_se3}
\left(\operatorname{d\,\tau}_{(\omega, \boldsymbol{\gamma})}\right)^{-1}
= \begin{bmatrix}
\mathbf{I}_3 - \frac{1}{2}\omega +
  \frac{1}{4}\boldsymbol{\omega} \boldsymbol{\omega}^T & \mathbf{0} 
  \\[6pt]
      -\frac{1}{2}\left(\mathbf{I}_3 -
        \frac{1}{2}\omega \right) \boldsymbol{\gamma} & \mathbf{I}_3 -
      \frac{1}{2} \omega 
\end{bmatrix}.
\end{align} 
Since we are using the pairing $\left<(\boldsymbol{\mu} , 
\boldsymbol{\eta} ),(\boldsymbol{\omega} , \boldsymbol{\gamma})\right> = 
\boldsymbol{\mu} \cdot \boldsymbol{\omega} + \boldsymbol{\eta} \cdot 
\boldsymbol{\gamma}$ between $\mathfrak{se}(3)\simeq \mathbb{R}^3 \times 
\mathbb{R}^3$ and its dual $\mathfrak{se}(3)^*\simeq \mathbb{R}^3 \times 
\mathbb{R}^3$, the matrix of $\left(\left(\operatorname{d\,\tau}_{(\omega, \boldsymbol{\gamma})}\right)^{-1}\right)^*$ is 
the transpose of \eqref{matrix_dcay_se3}.

These formulas are used in the implementation of the numerical algorithms
that we shall develop in the rest of the paper.

\color{black}

\subsection{Discrete covariant Hamilton principle}

The discrete covariant Euler-Lagrange equations (DCEL) are obtained from the Discrete Covariant Hamilton Principle
\begin{equation}\label{discrete_cov_HP} 
\delta  \sum_{j=0}^{N-1}\sum_{a=0}^{A-1} \mathscr{L}_d(\triangle_a^j, g_a^j, \xi_a^j, \eta_a^{j}) =0,
\end{equation} 
for arbitrary variations of $ g _a ^j $ vanishing at the boundary. This is the discrete version of the variational principle \eqref{trivialized_HP}.

The variations of $\xi_a^j$ and $\eta_a^j$ induced by variations of $ g _a ^j $ are computed as
\begin{equation}\label{variations}
\begin{aligned}
\delta \xi^j_a & =  {\rm d}\tau^{-1}_{\Delta t \xi^j_a }  \left(-\zeta ^j _a +
\operatorname{Ad}_{\tau(\Delta t\xi^j_a )}\zeta ^{j+1}_a \right)/\Delta t,  \\
\delta \eta ^j_a & =  {\rm d}\tau^{-1}_{\Delta s \eta ^j_a }  \left(-\zeta ^j _a +
\operatorname{Ad}_{\tau(\Delta s \eta ^j_a )}\zeta ^{j}_{a+1} \right)/ \Delta s,
\end{aligned}
\end{equation}
where $\zeta^j_a= (g^j_a)^{-1} \delta g_a^j$. Here 
$ \operatorname{Ad}_g \xi $ and $\operatorname{Ad}^*_{g^{-1}} \mu$ 
denote the (left) adjoint and coadjoint representations of $G$ on 
$\mathfrak{g}$ and $\mathfrak{g}^\ast$, respectively, where  
$g \in G$, $\xi \in \mathfrak{g}$, and $\mu \in  \mathfrak{g}^\ast$. 
For example, if $G=SE(3)$, which is the Lie group used later on in
the numerical algorithm for the beam, the concrete expressions for
the adjoint and coadjoint actions are
\begin{equation}
\label{ad_coad_se3}
\operatorname{Ad}_{(\Lambda , \mathbf{r})}(\boldsymbol{\omega} , \boldsymbol{\gamma}) =(\Lambda \boldsymbol{\omega}, 
\Lambda \boldsymbol{\gamma} + \mathbf{r} \times 
\Lambda \boldsymbol{\omega}),\qquad 
\operatorname{Ad}_{(\Lambda , \mathbf{r})^{-1}}^\ast (\boldsymbol{\mu} , \boldsymbol{\nu})=
(\Lambda \boldsymbol{\mu} + \mathbf{r} \times \Lambda \boldsymbol{\nu} , 
\Lambda \boldsymbol{\nu} ).
\end{equation}
Returning to the general case, using the notations
\begin{align}
\label{mu_lambda_formula}
 \mu_a^j:=\left( {\rm d}\tau^{-1}_{\Delta t \xi^j_a }\right)^*\partial_\xi K\left(\xi_a^j\right) \quad \text{and}\quad \lambda_a^j:=\left({\rm d}\tau^{-1}_{\Delta s \eta ^j_a } \right)^*\partial_\eta \Phi\left(\eta_a^j\right),
\end{align}
a direct computation shows that when arbitrary variations are allowed, \eqref{discrete_cov_HP} yields
\begin{equation} \label{CDEL_beam}
\begin{aligned}
 &  \frac{1}{\Delta t} \left(- \mu_a^j + \mathrm{Ad}^*_{\tau(\Delta t \xi_a^{j-1})} \mu_a^{j-1} \right)     \\
&  \hspace{2 cm}  + \frac{1}{\Delta s} \left( \lambda_a^j  - \mathrm{Ad}^*_{\tau(\Delta s \eta_{a-1}^j)} \lambda_{a-1}^j   \right) - (g_a^j)^{-1} D_{g_a^j}\Pi(g_a^j)  = 0,  \\
& \qquad  \text{for all $j=1,..., N-1$, and $a=1,..., A-1$,}
\end{aligned}
\end{equation}
\begin{equation}\label{boun_cond1}
 \begin{aligned}
  \frac{1}{\Delta t} \left(-\mu_0^j + \mathrm{Ad}^*_{\tau(\Delta t \xi_0^{j-1})} \mu_0^{j-1} \right) + \frac{1}{\Delta s} \lambda_0^j & = (g_0^j)^{-1} D_{g_0^j}\Pi(g_0^j),  \qquad \quad \quad  \\
 \qquad \qquad     - \frac{1}{\Delta s} \mathrm{Ad}^*_{\tau(\Delta s \eta_{A-1}^j)} \lambda_{A-1}^j & = 0,  \\
 \qquad  \text{for all $j=1,..., N-1$,} 
 \end{aligned}
 \end{equation}
 \begin{equation}\label{boun_cond2}
 \begin{aligned}
  - \frac{1}{\Delta t} \mu_a^0 +   \frac{1}{\Delta s} \left( \lambda_a^0 -  \mathrm{Ad}^*_{\tau(\Delta s \eta_{a-1}^0)} \lambda_{a-1}^0  \right) & = (g_a^0)^{-1} D_{g_a^0}\Pi(g_a^0),\qquad\qquad   \\
 \qquad \qquad  \frac{1}{\Delta t} \mathrm{Ad}^*_{\tau(\Delta t \xi_a^{N-1})} \mu_a^{N-1} & =0, \\
 \qquad  \text{for all $a =1,..., A-1$,}
 \end{aligned}
 \end{equation}
 \begin{equation}\label{boun_cond3}
 \begin{aligned}
   - \frac{1}{ \Delta t} \mu_0^0 + \frac{1}{\Delta s} \lambda_0^0 & = (g_0^0)^{-1} D_{g_0^0}\Pi(g_0^0), \qquad \qquad \qquad  \qquad \qquad \quad  \\
 \frac{1}{\Delta t}  \mathrm{Ad}^*_{\tau(\Delta t \xi_0^{N-1})} \mu_0^{N-1} &=0, \\
  -\frac{1}{\Delta s} \mathrm{Ad}^*_{\tau(\Delta S \eta_{A-1}^0)} \lambda_{A-1}^0 &= 0.
\end{aligned}
\end{equation}
Equations \eqref{CDEL_beam}, are associated to the interior variations 
$ \delta g _a ^j $, $j=1,...,N-1$, $a=1,...,A-1$. Equations \eqref{boun_cond1}, \eqref{boun_cond2}, \eqref{boun_cond3} are associated to the boundary variations 
$\delta g _a ^j $ for $j=0$, $j=N$, $a=0,...,A$ or for $j=0,...,N$, $a=0$, $a=A$. We refer to the Appendix for the explicit expressions of the equations \eqref{mu_lambda_formula}--\eqref{boun_cond3} for the case $G=SE(3)$. We thus obtain the following result.

\begin{proposition}
Given a discrete Lagrangian $\mathscr{L}_d: X _d ^\triangle \times G \times  \mathfrak{g}  \times \mathfrak{g}   \rightarrow \mathbb{R} $ 
and a discrete field $ g _d =\{g _a ^j \}$, the following conditions are equivalent.
\begin{itemize}
\item[{\rm (i)}] The Discrete Covariant Hamilton Principle
\[
\delta \sum_{j=0}^{N-1}\sum_{a=0}^{A-1} \mathscr{L}_d(\triangle_a^j, g_a^j, \xi_a^j, \eta_a^j)=0
\]
holds for arbitrary variations of $\delta g_a^j$ vanishing at the spacetime boundary : $j=0, N$, and $a=0, A$.
\item[{\rm (ii)}] The discrete field $g_d$ satisfies the DCEL equations
\begin{equation}\label{CDEL_beam_prop}
\begin{aligned}
 &  \frac{1}{\Delta t} \left(- \mu_a^j + \mathrm{Ad}^*_{\tau(\Delta t \xi_a^{j-1})} \mu_a^{j-1} \right)     \\
&  \hspace{1.5 cm}  + \frac{1}{\Delta s} \left( \lambda_a^j  - \mathrm{Ad}^*_{\tau(\Delta s \eta_{a-1}^j)} \lambda_{a-1}^j   \right) - (g_a^j)^{-1} D_{g_a^j}\Pi_{d}(g_a^j)  = 0,  \\
& \qquad  \text{for all $j=1,..., N-1$, and $a=1,... A-1$.}
\end{aligned}
\end{equation}
\end{itemize}
\end{proposition}

\begin{remark}[Discrete versus continuous equations] {\rm Note that the terms
\[
\frac{1}{\Delta t} \left(\mu_a^j - \mathrm{Ad}^*_{\tau(\Delta t \xi_a^{j-1})} \mu_a^{j-1} \right) \quad\text{and}\quad  \frac{1}{\Delta s} \left( \lambda_a^j  - \mathrm{Ad}^*_{\tau(\Delta s \eta_{a-1}^j)} \lambda_{a-1}^j   \right) 
\]
in the discrete equation \eqref{CDEL_beam} are, respectively, the discretization of the terms
\[
\frac{d}{dt}\frac{\partial K}{\partial \xi }  - \mathrm{ad}_\xi^* \frac{\partial K}{\partial \xi }\quad\text{and}\quad  \frac{d}{ds}\frac{\partial \Phi }{\partial \eta}- \mathrm{ad}_\eta^*\frac{\partial \Phi }{\partial \eta}
\]
of the continuous equation \eqref{Triv_CEL_BEam}. This reflects the covariant point of view we have used to derive the discrete equations of motion, namely, that the variables $t$ and $s$ are treated in the same way, so that discrete velocities and discrete gradients are approximated in the same geometry preserving way.}
\end{remark} 
 
\begin{remark}[Boundary conditions]\label{BC}{\rm
For simplicity we have considered above only the case when the configuration is fixed at the spacetime boundary, so that all variations vanish at the boundary. However, it is important to note that, as in the continuous case, if the configuration is not prescribed on certain 
subsets of the boundary, then natural discrete boundary conditions emerge from the variational principle. These are the discrete zero-traction boundary conditions (at the spatial boundary) and the discrete zero-momentum boundary conditions (at the temporal boundary), obtained from \eqref{boun_cond1}--\eqref{boun_cond3}. We refer to \cite{DeGBRa2013} 
for a detailed treatment.}
\end{remark} 
 
\begin{remark}\label{remark} 
 {\rm As it is usually done, and in a similar way with the continuous setting, the discrete equations are obtained by formulating the Discrete Hamilton Principle for a boundary value problem: the field is prescribed at the boundary of the spacetime domain. We will, however, solve these equations as initial value problems.
 
First, we will assume that the initial configuration $g(0,s)$ and initial velocity $\partial _t g(0,s)$ are given, and we compute the time evolution of the beam. In the discrete setting, this means that $ \boldsymbol{g} ^0 =(g _0 ^0 , ..., g _A ^0)$ and $ \boldsymbol{g} ^1 =(g _0 ^1 , ..., g _A ^1)$ are given and we solve for $\boldsymbol{g}^i$, $i=2,...,N$. 

Second, we will assume that the evolution and the deformation gradient of an extremity (say $s=0$) are known for all time, i.e., $g(t,0)$ and $\partial _s g(t,0)$ are given for all $t \in [0,T]$. We want to reconstruct the dynamics $g(t,s)$ in space  for all $s$. In the discrete setting, this means that $ \boldsymbol{g}_0 =(g _0 ^0 , ..., g _0 ^N)$ and $ \boldsymbol{g} _1 =(g _1 ^0 , ..., g _1 ^N)$ are given and we solve for $\boldsymbol{g}_a$, $a=2,...,N$.}
\end{remark}

\subsection{External forces}

External forces can be easily included in the discrete equations by considering the discrete analogue of the covariant Lagrange-d'Alembert principle \eqref{eq:LDAvar}. 
The discrete Lagrangian forces are maps $F_d^k : X^\triangle_d \times G \times G \times G  \rightarrow T^*G$, $k=1,2,3$, with
$F^1_d (\triangle _a ^j , g_a^j, g_a^{j+1}, g_{a+1}^j) \in T_{g_a^j}^*G$, $F^2_d (\triangle _a ^j, g_a^j, g_a^{j+1}, g_{a+1}^j ) \in 
T_{g_a^{j+1}}^*G$, $F^3_d (\triangle _a ^j, g_a^j, g_a^{j+1}, g_{a+1}^j ) \in T_{g_{a+1}^j}^*G$,
which are fiber preserving. 
Let $f^k : X^\triangle_d \times G \times \mathfrak{g} \times \mathfrak{g}  \rightarrow G \times \mathfrak{g}^*$
be the trivialized Lagrangian forces. 

The discrete covariant Lagrange-d'Alembert principle is
\begin{align*}  
\begin{split}
 & \, \delta\sum_{j=0}^{N-1} \sum_{a=0}^{A-1} 
 \mathscr{L}_d(\triangle_a^j, g_a^j, \xi_a^j, \eta_a^j)
  + \sum_{j=0}^{N-1} \sum_{a=0}^{A-1 } \Delta t \Delta s \left[ \langle f^1 (\triangle_a^j, 
  g_a^j, \xi_a^j, \eta_a^j), (g_a^j)^{-1} \delta g_a^j \rangle \right. \\
&  \hspace{.7 cm} \left. 
+ \langle  f^2 (\triangle_a^j, g_a^j, \xi_a^j, \eta_a^j), (g_{a}^{j+1})^{-1}\delta g_{a}^{j+1}\rangle  
+ \langle  f^3 (\triangle_a^j, g_a^j, \xi_a^j, \eta_a^j), (g_{a+1}^{j})^{-1}\delta g_{a+1}^{j}\rangle \right] = 0,
\end{split}
\end{align*}  
for arbitrary variations of $ g _a ^j $ vanishing at the boundary. It yields
\begin{equation*}
\begin{aligned}
&\frac{1}{\Delta t} \left(- \mu_a^j + \mathrm{Ad}^*_{\tau(\Delta t \xi_a^{j-1})} \mu_a^{j-1} \right)   + \frac{1}{\Delta s} \left( \lambda_a^j  - \mathrm{Ad}^*_{\tau(\Delta s \eta_{a-1}^j)} \lambda_{a-1}^j   \right) -  (g_a^j)^{-1} D_{g_a^j}\Pi_{d}(g_a^j)  \\
&\;\;\;+f^1(\triangle_a^j, g_a^j, \xi_a^j, \eta_a^j) + f^2(\triangle_a^{j-1}, 
g_a^{j-1}, \xi_a^{j-1}, \eta_a^{j-1})+f^3(\triangle_{a-1}^j, g_{a-1}^j, 
\xi_{a-1}^j, \eta_{a-1}^j)=0,
\end{aligned}
\end{equation*}
for all $j=1,..., N-1$, and $a=1,... A-1$.

\subsection{Time-stepping algorithm}
 
The complete algorithm obtained via the covariant
variational integrator can be implemented according to the following
steps. 

\begin{framed} 
\bf{Time Integrator}
\begin{align*} 
& \vspace{0.3cm}\text{Given: } g^j_a, \xi^{j-1}_a, \mu_a^{j-1}, f_a^{j}, \quad \text{for } a=0, ..., A, \phantom{\int }\\
& \text{Compute: } & \\
& \quad \eta_a^j=\frac{1}{s} \tau^{-1}\left((g_{a}^j)^{-1} g_{a+1}^j\right), \\ 
& \quad \lambda_a^j=\left( \operatorname{d} \tau_{s\eta_a^j}^{-1} \right)^*\partial_\eta
  \Phi\left(\eta_a^j\right), \\
  &\quad \mu_a^j = \operatorname{Ad}_{\tau(h\xi_a^{j-1})}^*\mu_a^{j-1}\\
  &\quad \quad \quad   + \left\{ \begin{array}{ll} 
     \vspace{0.2cm} h
      \left(\frac{1}{s} \lambda_{a}^j + f_a^j\right),  &\text{for } a=0,\\
      \vspace{0.1cm}h
      \left(\frac{1}{s} \left(\lambda_{a}^j-\operatorname{Ad}_{\tau(s
          \eta_{a-1}^{j})}^*\lambda_{a-1}^{j}\right) + f_a^j\right), &\text{for } a=1,...,A-1,
    \end{array}  \right.  \\ 
  &\text{Solve the discrete Legendre transform: }
  \mu_a^j=\left( \operatorname{d}\tau_{h\xi_a^j}^{-1} \right)^*\partial_\xi
  K\left(\xi_a^j\right), \text{ for } \xi^j_a
& \phantom{\int}\\
  &\text{Update:} \ g_a^{j+1} = g_a^j\tau\left(h\xi_a^j\right).
\end{align*}
\end{framed}
We obtain $g_A^{j+1}$ through the boundary condition 
$\operatorname{Ad}_{\tau(s\eta_{A-1}^{j})}^*\lambda_{A-1}^{j+1}=0$ as defined in \eqref{boun_cond1}. Note that, for clarity, we have denoted the total external force at point $(a,j)$ by $f_a^j$. The algorithm requires the implicit solution of the Legendre transform which is locally invertible for an appropriately chosen time-step. The rest of the update is performed explicitly.

Through the proposed unified space-time description it is now possible to study the symplectic-momentum preservation properties of the algorithm in both space and time directions. This is accomplished by developing a discrete analog of Noether's theorem as follows.

\section{Discrete momentum maps and Noether theorem}\label{sec4}

Recall that for the discretization of standard (i.e., based on trajectories evolving in time) Lagrangian mechanics, the tangent bundle $TQ$ is replaced by the Cartesian product $Q \times Q$. Given a discrete Lagrangian function $L_d:Q \times Q \rightarrow \mathbb{R}$, the discrete Legendre transforms are the two maps
\begin{equation}\label{discrete_momap_def}
\begin{aligned} 
\mathbb{F}  L_d ^\pm: Q \times Q \rightarrow T^*Q, \qquad 
&\mathbb{F}  L_d ^-( q ^j, q^{j+1})= - D _1 L_d( q ^j , q^{j+1})\\
&\mathbb{F}  L_d ^+( q ^j, q^{j+1})=D_2L_d( q ^j , q^{j+1}).
\end{aligned}
\end{equation} 
The discrete Euler-Lagrange equations can be written as 
\[
\mathbb{F}  L_d ^+( q ^{j-1}, q^{j}) = \mathbb{F}L_d^-( q ^j, q^{j+1}).
\]
Given a Lie group action, the discrete momentum maps are defined
analogously to the continuous expression \eqref{Lagr_momap}, namely
\begin{equation}\label{discrete_momap} 
\mathbf{J}_{L_d}^\pm: Q \times Q \rightarrow \mathfrak{g}^\ast , 
\qquad \left\langle\mathbf{J}_{L_d}^\pm( q ^j, q^{j+1}), 
\xi \right\rangle=
\left\langle \mathbb{F}  L_d ^\pm( q ^j, q^{j+1}), 
\xi_Q \right\rangle, \quad \forall \xi \in \mathfrak{g}.
\end{equation} 
If the discrete Lagrangian is $G$-invariant, then the two momentum maps coincide, $\mathbf{J}_{L_d}^+=\mathbf{J}_{L_d}^-=: \mathbf{J} _{L_d}$, and $\mathbf{J} _{L_d}$ is conserved along the solutions of the discrete Euler-Lagrange equations.

With this is mind, we will next construct more general Legendre transforms and associated momentum maps which extend to the space-time domain.

\subsection{Discrete covariant Legendre transforms}

In the present covariant discretization of Lagrangian mechanics, it is natural to define three Legendre transforms $\mathbb{F}^k\mathcal{L}_d: 
X^\triangle_d \times G \times  G \times  G \rightarrow T^*G$, 
$k=1,2,3$ associated to a discrete Lagrangian 
$\mathcal{L}_d = \mathcal{L}_d(\triangle_a^j, g_a ^j , 
g_a^{j+1}, g_{a+1}^j )$, defined in an analogous way to \eqref{discrete_momap_def}, namely
\begin{equation*} 
\mathbb{F}^k\mathcal{L}_d(\triangle_a^j, g_a ^j , 
g_a^{j+1}, g_{a+1}^j ) = \left(  g(k),\partial _{g(k)}  \mathcal{L}_d(\triangle_a^j)\right)\in T^*_{g(k)}G, \quad k=1,2,3,
\end{equation*}
where $g(1)=g_a^j $, $g(2)= g_a^{j+1}$, $g(3)=g_{a+1}^j$.

In terms of the discrete Lagrangian $\mathscr{L}_d= \mathscr{L}_d(\triangle_a^j, g _a ^j , \xi _a ^j , \eta _a ^j )$ in \eqref{DCL_beam}, the corresponding discrete Legendre transforms $\mathbb{F}^k\mathscr{L}_d: X_d^\triangle \times G \times  \mathfrak{g} \times \mathfrak{g} \rightarrow G  \times \mathfrak{g}^* $, $k=1,2,3$, are
\begin{align*}
\mathbb{F}^1 \mathscr{L}_d(\triangle_a^j, g_a ^j , 
\xi_a^{j}, \eta_{a}^j ) & = \left(g_a^j,  -  \Delta s  \mu_a^j +  \Delta t \lambda_a^j - \Delta t \Delta s (g_a^j)^{-1} D_{g_a^j}\Pi_{d}(g_a^j) \right) , \nonumber\\
\mathbb{F}^2\mathscr{L}_d(\triangle_a^{j}, g_a ^j , 
\xi_a^{j}, \eta_{a}^j ) & = \left( g_a^{j+1}, \Delta s \mathrm{Ad}^*_{\tau(\Delta t \xi_a^j)} \mu_a^j \right) , \\
\mathbb{F}^3 \mathscr{L}_d(\triangle_{a}^j, g_a ^j , 
\xi_a^{j}, \eta_{a}^j ) & = \left(g_{a+1}^{j}, - \Delta t \mathrm{Ad}^*_{\tau(\Delta s \eta_a^j)} \lambda_a^j  \right).\nonumber
\end{align*}
We note that the CDEL equations \eqref{CDEL_beam} can be written as
\[
\mathbb{F}^1 \mathscr{L}_d(\triangle_a^j, g_a ^j , 
\xi_a^{j}, \eta_{a}^j ) +  
\mathbb{F}^2 \mathscr{L}_d(\triangle_a^{j-1}, g_a ^{j-1} , 
\xi_a^{j-1}, \eta_{a}^{j-1} ) + \mathbb{F}^3 \mathscr{L}_d(
\triangle_{a-1}^j, g_{a-1}^j , 
\xi_{a-1}^{j}, \eta_{a-1}^j )  =0.
\]

\subsection{Discrete covariant momentum maps}\label{DCMP} 
In order to study the integrator preservation properties, we consider symmetries given by the action of a subgroup $H \subset G$ acting on $G$ by multiplication on the left. The infinitesimal generator of the left multiplication by $H$ on $G$, associated to the Lie algebra element $\zeta \in \mathfrak{h}$, is expressed as $\zeta_G(g) = \zeta g$.
The three discrete Lagrangian momentum maps 
$J^k_{\mathscr{L}_d} : X_d^\triangle \times G \times \mathfrak{g} \times \mathfrak{g}  \rightarrow \mathfrak{h}^*$, $k= 1,2,3$, are defined, analogously to \eqref{discrete_momap}, by
\[
\left< J^k_{\mathscr{L}_d}(\triangle_a^j, g_a ^j , 
\xi_a^{j}, \eta_{a}^j ), \zeta \right> = \left< 
\mathbb{F}^k\mathscr{L}_d (\triangle_a^j , g_a ^j , 
\xi_a^{j}, \eta_{a}^j ) ,  \ (g(k))^{-1} \zeta_{G}(g(k)) \right>,\quad \zeta \in \mathfrak{h},
\]
where, as before, we use the notations $g(1)=g_a^j$, $g(2)=g_a^{j+1}$, and $g(3)=g_{a+1}^j$. Moreover, $\mathbb{F}^k\mathscr{L}_d$ is seen as an element of $\mathfrak{g}^*\cong (G\times \mathfrak{g}^*)/G$. 
For the discrete Lagrangian $\mathscr{L}_d$ in \eqref{DCL_beam}, these momentum maps are
\begin{align}
\label{discrete_momentum_map}
J^1_{\mathscr{L}_d}(\triangle_a ^j, g_a ^j , 
\xi_a^{j}, \eta_{a}^j )&= i ^\ast \mathrm{Ad}_{(g_a^j)^{-1}}^* \left( -  \Delta s \mu_a^j +  \Delta t \lambda_a^j - \Delta t \Delta s (g_a^j)^{-1} \mathbf{D}_{g_a^j}\Pi_{d}(g_a^j) \right), \nonumber \\
J^2_{\mathscr{L}_d}(\triangle_a ^j, g_a ^j , 
\xi_a^{j}, \eta_{a}^j )&= i ^\ast \operatorname{Ad}_{(g_a^{j+1})^{-1}}^*\left( \Delta s \mathrm{Ad}^*_{\tau(\Delta t \xi_a^j)} \mu_a^j \right) , \\
J^3_{\mathscr{L}_d}(\triangle_a ^j, g_a ^j , 
\xi_a^{j}, \eta_{a}^j )&= i ^\ast\operatorname{Ad}_{(g_{a+1}^j)^{-1}}^* \left( - \Delta t \mathrm{Ad}^*_{\tau(\Delta s \eta_a^j)} \lambda_a^j  \right), \nonumber
\end{align}
where $ i ^\ast : \mathfrak{g}  ^\ast \rightarrow \mathfrak{h} ^\ast $ is the dual map to the Lie algebra inclusion $i: \mathfrak{h}  \rightarrow \mathfrak{g}  $.

The discrete Noether theorem follows from  the invariance of the discrete action under the left action of a Lie group $H$. In order to obtain the associated conservation law, we shall proceed exactly as in the continuous setting. The variations of $ g _d $ induced from this action are $ \delta g _a ^j =\zeta g _a ^j $, where $ \zeta \in \mathfrak{h}$; we thus get $(g _a ^j) ^{-1} \delta g _a ^j=  \operatorname{Ad}_{ ( g _a ^j ) ^{-1} } \zeta $. Assuming $H$-invariance of the discrete covariant Lagrangian $\mathscr{L}_d$ and assuming that the discrete Euler-Lagrange equations are satisfied, we get (from similar computations that lead to \eqref{CDEL_beam}--\eqref{boun_cond3}), the discrete version of the global Noether theorem: for all $0\leq B<C\leq A-1$, $0\leq K<L\leq N-1$, we have the conservation law
\begin{equation}\label{DCN}  
\mathscr{J}_{B,C}^{K,L}(g _d  ) = 0,
\end{equation}
where,
\begin{align}\label{Def_DN}  
\mathscr{J}_{B,C}^{K,L}(g _d  ) :=&\sum_{j=K+1}^L \left( J ^1_{\mathcal{L}_d}( \triangle _B ^j )+J ^2_{\mathcal{L}_d}(  \triangle _B ^{j-1} )+ J ^3_{\mathcal{L}_d}(\triangle _C ^j ) \right) \nonumber\\
&+ \sum_{a=B+1}^C \left( J ^1_{\mathcal{L}_d}( \triangle _a ^K )+J ^2_{\mathcal{L}_d}( \triangle _a ^L)+ J ^3_{\mathcal{L}_d}( \triangle _{a-1} ^K ) \right) \\
&+ J ^1_{\mathcal{L}_d}( \triangle _B ^K )+J ^2_{\mathcal{L}_d}( \triangle _B ^L)+ J ^3_{\mathcal{L}_d}(  \triangle _C  ^K ),\nonumber
\end{align}
where we have abbreviated $J^k_{\mathscr{L}_d}(\triangle_a ^j, g_a ^j , 
\xi_a^{j}, \eta_{a}^j )$ by $J^k_{\mathscr{L}_d}(\triangle_a ^j)$.
\medskip

\noindent \textbf{Notation:} From now on, in order to simplify
notation, we shall adopt this abbreviation, i.e., $J^k_{\mathscr{L}_d}(\triangle_a ^j): = J^k_{\mathscr{L}_d}(\triangle_a ^j, g_a ^j , 
\xi_a^{j}, \eta_{a}^j )$.

\subsection{Symplectic properties of the time and space discrete evolutions} \label{time_space_symplectic_property}
 
In this subsection, we shall verify the symplectic character of the integrator in both time and space evolution. This is achieved by defining, from the discrete  covariant Lagrangian density $\mathscr{L}_d$, two ``classical'' Lagrangians, namely one associated to time evolution  and one associated to space evolution, as done in \cite{DeGBRa2013}.

The time-evolution and space-evolution Lagrangians are defined, respectively, by
\begin{equation}\label{def_L_N_d}
\begin{aligned}  
 \mathsf{L}_d(\mathbf{g}^j, \boldsymbol{\xi}^j) & :=\sum_{a=0}^{A-1}  \mathscr{L}_d (\triangle_a^j, g_a^j, \xi_a^j, \eta_a^j)\quad\text{and}\quad   \mathsf{N}_d(\mathbf{g}_a, \boldsymbol{\eta}_a) & :=\sum_{j=0}^{N-1} \mathscr{L}_d (\triangle_a^j, g_a^j, \xi_a^j, \eta_a^j),
\end{aligned}
\end{equation} 
where $\mathbf{g}^j =(g_0^j, ..., g_A^j)$ and 
$\mathbf{g}_a=(g^0_a, ...,g^N_a)$.
 
It can be verified that the discrete Euler-Lagrange equations for both $ \mathsf{L}_d$ and $\mathsf{N}_d$ are equivalent to the covariant discrete Euler-Lagrange equations \eqref{CDEL_beam_prop} for $\mathscr{L}_d$.
For this verification, one has to take into account the boundary conditions involved in the problem. For example, concerning $\mathsf{L}_d$, boundary conditions in time can be treated as usual by allowing (or not) variations at the temporal boundary. However, if boundary conditions are imposed at the spatial boundary, then one has to incorporate them in the configuration space of the Lagrangian $\mathsf{L}_d$ and, as such, they have to be time independent.  The same comments hold for $\mathsf{N}_d$ with the role of time and space reversed. We refer to \cite{DeGBRa2013} for a detailed discussion.

From this discussion, and from the general result that discrete Euler-Lagrange equations yield symplectic integrators (see, e.g., \cite{MaWe2001}), it follows that both discrete flows $\mathbf{g}^j $, $j=0,...,N$ and $\mathbf{g}_a $, $a=0,...,A$, are symplectic. Here again, the space on which this symplecticity occurs depends on the boundary conditions assumed. These two flows correspond to the discrete time and the discrete space evolutions associated to the discrete field $g_a^j$, $j=0,...,N$, $a=0,...,A$, respectively.

Furthermore, if $\mathscr{L}_d$ is $H$-invariant, then both $ \mathsf{L}_d$ and $\mathsf{N}_d$ inherit this $H$-invariance. Adapting the general formula \eqref{discrete_momap} to the trivialized 
Lagrangians $\mathsf{L}_d$ and $\mathsf{N}_d$, we get 
the discrete Lagrangian momentum maps $\mathbf{J}_{\mathsf{L}_d}^\pm: 
SE(3)^A \times \mathfrak{se}(3)^A \rightarrow \mathfrak{h}^\ast$ 
and $\mathbf{J}_{\mathsf{N}_d}^\pm: SE(3)^N \times \mathfrak{se}(3)^N 
\rightarrow \mathfrak{h}^\ast$
\begin{equation} \label{time_momentum_map}
\mathbf{J} _{\mathsf{L}_d}^-(\mathbf{g}^j, \boldsymbol{\xi}^j)  = 
-\sum_{a=0}^{A-1}\left( J^1_{ \mathscr{L}_d}(\triangle_a^j)  + J^3_{\mathscr{L}_d}(\triangle_a^j) \right), \qquad
\mathbf{J} _{\mathsf{L}_d}^+(\mathbf{g}^j, \boldsymbol{\xi}^j)  = \sum_{a=0}^{A-1} J^2_{\mathscr{L}_d}(\triangle_a^j),
\end{equation}
\begin{equation} \label{space_momentum_map}
\mathbf{J} _{\mathsf{N}_d}^-(\mathbf{g}_a, \boldsymbol{\eta}_a) = 
-\sum_{j=0}^{N-1}\left( J^1_{\mathscr{L}_d}(\triangle_a^j)  + J^2_{\mathscr{L}_d}(\triangle_a^j) \right), \qquad
\mathbf{J} _{\mathsf{N}_d}^+(\mathbf{g}_a, \boldsymbol{\eta}_a)  = \sum_{j=0}^{N-1}J^3_{\mathscr{L}_d}(\triangle_a^j).
\end{equation}

However, from this fact, one cannot conclude that  $\mathbf{J} _{\mathsf{N}_d}$ and $\mathbf{J} _{\mathsf{L}_d}$ are necessarily conserved by the discrete dynamics. Indeed, as discussed earlier, the imposition of 
boundary conditions can break the $H$-invariance since the space on which 
$\mathsf{L}_d$, $\mathsf{N}_d$ have to be redefined is no longer 
$H$-invariant. We refer to \cite{DeGBRa2013} for a detailed account. Such 
a phenomenon is not surprising since it already occurs in the continuous 
setting. In the discrete setting, it is explained via the following 
lemma that relates the expression of the covariant and classical discrete momentum maps.

\begin{lemma}\label{important_lemma}  
When $B=0$ and $C=A-1$, or $K=0$ and $L=N-1$, we have, respectively
\begin{equation}\label{JL}
\begin{aligned}
\mathscr{J}_{0,A-1}^{K,L}( g _d )=&\sum_{j=K+1}^L \left( J ^1_{\mathcal{L}_d}( \triangle _0 ^j )+J ^2_{\mathcal{L}_d}(   \triangle _0 ^{j-1} )+ J ^3_{\mathcal{L}_d}( \triangle _{A-1} ^j ) \right)\\
&+ \mathbf{J}^+_{\mathsf{L}_d}( \mathbf{g} ^L , \boldsymbol{\xi} ^{L})- \mathbf{J}^-_{\mathsf{L}_d}( \mathbf{g} ^K , \boldsymbol{\xi}  ^{K})
\end{aligned} 
\end{equation}
\begin{equation}\label{JN}   
\begin{aligned}
\mathscr{J}_{B,C}^{0,N-1}( g _d )=&\sum_{a=B+1}^C \left( J ^1_{\mathcal{L}_d}( \triangle _a ^0 )+J ^2_{\mathcal{L}_d}(  \triangle _a ^{N-1} )+ J ^3_{\mathcal{L}_d}(  \triangle _{a-1} ^0 ) \right)\\
&+ \mathbf{J}^+_{\mathsf{N}_d}( \mathbf{g} _C , \boldsymbol{\eta}  _{C})- \mathbf{J}^-_{\mathsf{N}_d}(\mathbf{g} _B , \boldsymbol{\eta} _{B}).
\end{aligned}
\end{equation} 
\end{lemma} 

From this result, we deduce the following theorem that explains the dependence of the validity of the Noether theorem for $ \mathcal{L} _d $, $\mathsf{L}_d$, and $\mathsf{N}_d$ on the boundary conditions imposed.

\begin{theorem}\label{summary_Noether}  
Consider the discrete Lagrangian density $\mathcal{L}_d$ and the 
associated discrete Lagrangians $\mathsf{L}_d$ and $\mathsf{N}_d$ as 
defined in \eqref{def_L_N_d}. Consider the discrete covariant momentum 
maps $J^k_{\mathcal{L}_d}$ {\rm (}see \eqref{discrete_momentum_map}{\rm)} 
and the discrete momentum maps $\mathbf{J}^\pm_{\mathsf{L}_d}$, 
$\mathbf{J}^\pm_{\mathsf{N}_d}$ {\rm (}see \eqref{time_momentum_map}, 
\eqref{space_momentum_map}{\rm )} associated to the action of $H$.
Suppose that the discrete covariant Lagrangian density $\mathcal{L}_d$ 
is $H$-invariant. While the discrete covariant Noether theorem 
$\mathscr{J}_{B,C}^{K,L}(g_d) = 0$ {\rm (}see \eqref{Def_DN}{\rm )} is always verified, independently on the imposed boundary conditions, the validity of the discrete Noether theorems for $\mathbf{J}_{\mathsf{L}_d}^\pm$ and $\mathbf{J}_{\mathsf{N}_d}^\pm$ depends on the boundary conditions, in a similar way with the continuous setting.

More precisely, if the configuration is prescribed at the temporal extremities and zero-traction boundary conditions are used, then the discrete momentum map $\mathbf{J}_{\mathsf{L}_d}=\mathbf{J}_{\mathsf{L}_d} ^+= \mathbf{J}_{\mathsf{L}_d}^-$ is conserved. In general, conservation of $\mathbf{J}_{\mathsf{N}_d}^\pm$ does not hold in this case.

If the configuration is prescribed at the spatial extremities and zero-momentum boundary conditions are used, then the discrete momentum map 
$\mathbf{J}_{\mathsf{N}_d}=\mathbf{J}_{\mathsf{N}_d}^+=\mathbf{J}_{\mathsf{N}_d}^-$ is conserved. In general, conservation of 
$\mathbf{J}_{\mathsf{L}_d}^\pm$ does not hold in this case.
\end{theorem}

\section{Numerical examples}\label{sec5}
  
The covariant variational integrator derived in this paper stems from a 
unified geometric treatment of both time and space evolution. As a result, 
the integrator has analogous preservation properties in time \emph{and} space. It is multisymplectic and verifies the discrete covariant Noether theorem. In addition, under appropriate boundary conditions, it is symplectic in space or in time and preserves exactly the discrete momentum maps associated to space or time evolutions.

We shall illustrate these properties by considering first the usual initial value problem for beam dynamics, namely, the case when the position $g(s,0)$ and velocity  $\partial_t g(s,0)$ are given at time $t=0$, for all $s\in[0,L]$. Then, by switching the role of space and time variables, we will attempt to ``reconstruct'' the spatial motion of the beam starting with the spatial boundary condition (at node $0$) computed over the time interval $[0,T]$ and evolve it in space towards node $A$.

In addition, we consider the inverse problem, i.e., to spatially integrate the trajectory from the knowledge of $g(0,t), \partial_s g(0,t)$, at the extremity $s=0$, for all $t \in [0,T]$. We can then compare the computed global motions depending on whether the integration is performed first in time and then in space, or vice versa. The following terminology is used to clarify these points:
\begin{itemize}
\item \emph{time-integration} and \emph{space-reconstruction}: the case when standard time-stepping is performed first, after which the 
obtained time trajectory of node $a=0$ is used as an initial condition 
for ``reconstructing'' the full spatial motion towards node $a=A$.
\item \emph{space-integration} and \emph{time-reconstruction}: the 
case when the time-trajectory of node $a=0$ over the interval $[0,T]$ 
is given as a boundary condition and used to spatially evolve the 
motion. After that, the trajectory of the beam at $t=0$ is reconstructed 
in time towards time $t=T$.
\end{itemize}

We next analyze the behavior of these schemes through numerical simulation and computation of their motion invariants, i.e., the discrete energy and discrete momentum maps defined next.

\subsection{Discrete Lagrangian}
\label{subsec:discrete_Lagrangian}

Given is the discrete Lagrangian defined in \eqref{DCL_beam}, on a triangular mesh. In the case of the beam, the discrete Lagrangian 
$\mathsf{L}_d: M^{A+1} \times M^{A+1} \rightarrow \mathbb{R}$ is given by
\begin{equation} \label{model_L_d}
\mathsf{L}_d( \mathbf{g} ^j , \boldsymbol{\xi }^j ) =   \sum_{ a=0}^{A-1} \Big\{ \Delta t \Delta s  K(\xi^{j}_a) - \Delta t \Delta s \left[  \Phi(\eta_a^j) +  \Pi ( g _a^j ) \right] \Big\}
\end{equation}
and the discrete Lagrangian $\mathsf{N}_d: M^{N+1} 
\times M^{N+1} \rightarrow \mathbb{R}$ by 
\begin{equation} \label{model_N_d}
\mathsf{N}_d( \mathbf{g} _a , \boldsymbol{\eta }_a ) =  \sum_{ j=0}^{N-1} \Big\{ \Delta t \Delta s  K(\xi^{j}_a) - \Delta t \Delta s \left[  \Phi(\eta_a^j) +  \Pi ( g _a^j ) \right] \Big\}.
\end{equation}

Since $\xi_a^{N-1} = \tau ^{-1}\left((g_a^{N-1})^{-1}g_a^{N}\right)/\Delta t$ and 
the nodes $(N, a)$ lie outside of the domain, we impose
$K(\xi^{N-1}_a) =0$. Thus, \eqref{mu_lambda_formula} implies that the second set of zero-momentum boundary conditions $\mathrm{Ad}^*_{\tau(\Delta t \xi_a^{N-1})} \mu_a^{N-1}=0$  
in \eqref{boun_cond2} are verified. So, for the spatial evolution,
we impose the zero-momentum boundary conditions at $t=0$ an $t=T$
and the algorithm consists hence of \eqref{CDEL_beam} and the first
set in \eqref{boun_cond2}. As a consequence of 
$\mathrm{Ad}^*_{\tau(\Delta t \xi_a^{N-1})} \mu_a^{N-1}=0$, 
in the algorithm \eqref{CDEL_beam} we have $\mu_a^{N-1} = 0$, so for $j=N-1$, we get
\[
\frac{1}{\Delta t} \mathrm{Ad}^*_{\tau(\Delta t \xi_a^{N-2})} \mu_a^{N-2} +   \frac{1}{\Delta s} \left(   \lambda_a^{N-1} -   \mathrm{Ad}^*_{\tau(\Delta s \eta_{a-1}^{N-1})} \lambda_{a-1}^{N-1}  \right)   =  (g_{a}^{N-1})^{-1} D_{g_{a}^{N-1}}\Pi(g_{a}^{N-1}),
\]
for all $a =1,..., A-1$.

In the same way, since $\eta_{A-1}^j =  
\tau ^{-1}\big((g^j_{A-1})^{-1}g^j_A \big)/\Delta s$ and 
the nodes $(j,A)$ lie outside of the domain,  we impose
$\Phi(\xi^{j}_{A-1}) =0$, so, in view of \eqref{mu_lambda_formula},
the second set of zero traction boundary conditions 
$\mathrm{Ad}^*_{\tau(\Delta s \eta_{A-1}^j)} \lambda_{A-1}^j =0$ 
in \eqref{boun_cond1} holds. Thus, for the temporal evolution,
we impose the zero traction boundary conditions at $s=0$ and
$s=L$. The algorithm consists hence of \eqref{CDEL_beam} and the
first set in \eqref{boun_cond1}. As a consequence of
$\mathrm{Ad}^*_{\tau(\Delta s \eta_{A-1}^j)} \lambda_{A-1}^j =0$,
in the algorithm \eqref{CDEL_beam} we have $\lambda_{A-1}^j = 0$, so for $a=A-1$, we get
\[
\frac{1}{\Delta t} \left(-  \mu_{A-1}^j +   \mathrm{Ad}^*_{\tau(\Delta t \xi_{A-1}^{j-1})} \mu_{A-1}^{j-1} \right) - \frac{1}{\Delta s} \mathrm{Ad}^*_{\tau(\Delta s \eta_{A-2}^j)} \lambda_{A-2}^j   =  (g_{A-1}^j)^{-1} D_{g_{A-1}^j}\Pi(g_{A-1}^j),
\]
for all $j=1,..., N-1$.

\subsection{Evaluation criteria: discrete momentum and energy preservation}

Recall that since the discrete Lagrangian density $\mathscr{L}_d$ of the beam is $SE(3)$-invariant (see \eqref{DCL_beam}), the discrete covariant Noether theorem is verified, i.e.,
$\mathscr{J}_{B,C}^{K,L}(g _d) = 0$, see Theorem \ref{summary_Noether}. Recall also that the conservation of the discrete momentum maps $\mathbf{J}_{\mathsf{L}_d}$
and $\mathbf{J}_{\mathsf{N}_d}$ depends on the boundary conditions.

These conservation laws follow from the following particular cases of the covariant discrete Noether theorems, namely, $\mathscr{J}_{0,A-1}^{K,L}(g_d)=0$ and
$\mathscr{J}_{B,C}^{0,N-1}(g_d)=0$; see \eqref{JL} and \eqref{JN}.

We shall now consider the energies associated to the temporal and spatial evolutions.
The standard continuous \textit{time}-evolving energy is defined by
\begin{equation*}
\mathsf{E}_{ \mathsf{L}}(g, \xi )= \int_0^L \left( K(\xi) + \Phi(g ^{-1} \partial _s g) + \Pi(g)\right) \operatorname{d}\!s
\end{equation*}
and hence the corresponding discrete energy $\mathsf{E}_{\mathsf{L}_d}$ 
evaluated on the discrete time trajectory $\mathbf{g}^1,..., 
\mathbf{g}^{N-1}$ has the expression
\begin{equation} \label{formulas_reconstruction_space}
\mathsf{E}_{\mathsf{L}_d}(\mathbf{g}^j, \mathbf{g}^{j+1} ) =  
\sum_{a=0}^{A-1} K(\xi_a^j)  + \sum_{a=0}^{A-2} \left( \Phi(\eta_a^j) 
+   \Pi(g_a^j)\right) +  \Pi(g_{A-1}^j).
\end{equation}
On the other hand, the continuous ``energy'' evolving in \emph{space} 
is given by
\begin{equation*}
\mathsf{E}_{\mathsf{N}}(g, \eta)
 = \int_0^T \left( -K(g ^{-1} \partial _t g)  - \left\langle \mathbb{C}  (\eta - \mathbf{E} _6), \mathbf{E} _6 \right\rangle -\Phi(\eta ) 
 + \Pi(g)\right) \operatorname{d}\!t,
\end{equation*}
where $\mathbb{C}$ is the strain matrix \eqref{strain_matrix},
while the corresponding discrete energy $\mathsf{E}_{\mathsf{N}_d}$ is
\begin{align} \label{formulas_reconstruction_time_space}
  \mathsf{E}_{\mathsf{N}_d}(\mathbf{g}_a, \mathbf{g}_{a+1}) 
  &= - \sum_{j=0}^{N-2}  K(\xi_a^j)  + \sum_{j=1}^{N-2}  \left( - \left\langle \mathbb{C}  
  ( \eta_a^j - \mathbf{E} _6), \mathbf{E} _6 \right\rangle -\Phi(\eta_a^j ) + \Pi(g_a^j)\right) \nonumber \\
  & \quad - \frac{1}{2} \left\langle \mathbb{C}  
  ( \eta_a^0 - \mathbf{E} _6), \mathbf{E} _6 \right\rangle - \Phi(\eta_a^0 ) + \Pi(g_a^0) \nonumber \\
  & \quad - \frac{1}{2} \left\langle \mathbb{C}  
  ( \eta_a^{N-1} - \mathbf{E} _6), \mathbf{E} _6 \right\rangle -   \Phi(\eta_a^{N-1} ).
\end{align}

The symplecticity in time and in space of the discrete scheme obtained by discrete covariant Euler-Lagrange equations was studied in detail in \cite{DeGBRa2013}. The discussion depends on the boundary conditions considered and parallels the situation of the continuous setting. As a consequence, if the continuous energy $\mathsf{E}_\mathsf{L}$, resp., 
$\mathsf{E}_\mathsf{N}$, is preserved (which depends on the boundary conditions used), then the discrete energy $\mathsf{E}_{\mathsf{L}_d}$, resp., $\mathsf{E}_{\mathsf{N}_d}$, is approximatively preserved (i.e., it oscillates around its nominal value) due to the symplectic character of the scheme. The situation is summarized in Table~\ref{fig:methods} below. 
\begin{figure}[H]
  \begin{center}
    \begin{tabular}{|c|l|l|c|}
    \hline
Algorithm & Momentum & Energy behavior & Global Noether \\
\hline
\multirow{2}{*}{$\left\{\begin{tabular}{l} Time-integration  \\ 
Space-Reconstruction \end{tabular}\right.$} & 
$\mathbf{J}^\pm_{\mathsf{L}_d}$ exact & 
$\mathsf E_{\mathsf{L}_d}$ approx. pres.& 
exact \\
& $\mathbf{J}^\pm_{\mathsf{N}_d}$ not preserved & 
$\mathsf E_{\mathsf{N}_d}$ not preserved &
exact\\
\hline
\multirow{2}{*}{$\left\{\begin{tabular}{l} Space-integration  \\ 
Time-Reconstruction \end{tabular}\right.$} & 
$\mathbf{J}^\pm_{\mathsf{N}_d}$ exact & 
$\mathsf E_{\mathsf{N}_d}$ approx. pres. & 
exact \\
& $\mathbf{J}^\pm_{\mathsf{L}_d}$ not preserved & 
$\mathsf E_{\mathsf{L}_d}$ not preserved &
exact\\
\hline
\end{tabular}\caption{\footnotesize{Summary of algorithm preservation properties using the boundary conditions \eqref{boun_cond1}--\eqref{boun_cond3}.}}\label{fig:methods}
\end{center}
\end{figure}

\subsection{Time-integration and space-reconstruction}\label{TI_SR} 
 
The situation treated here corresponds to an usual initial value 
problem, namely, we assume that the initial configuration and 
velocity of the beam are known. We assume that the extremities 
evolve freely in space, which corresponds to zero-traction 
boundary conditions
\begin{equation}\label{cont_boun_cond}
\left. (\Gamma-\mathbf{E}_3)\right|_{s=0}=0, \quad 
\left. (\Gamma-\mathbf{E}_3)\right|_{s=L}=0, \quad \Omega  (0)=\Omega  (L)=0.
\end{equation}
In \S\ref{space_setting} we shall consider a different initial setting.

\begin{figure}[ht]
  \centering
  \includegraphics[width=1.5 in]{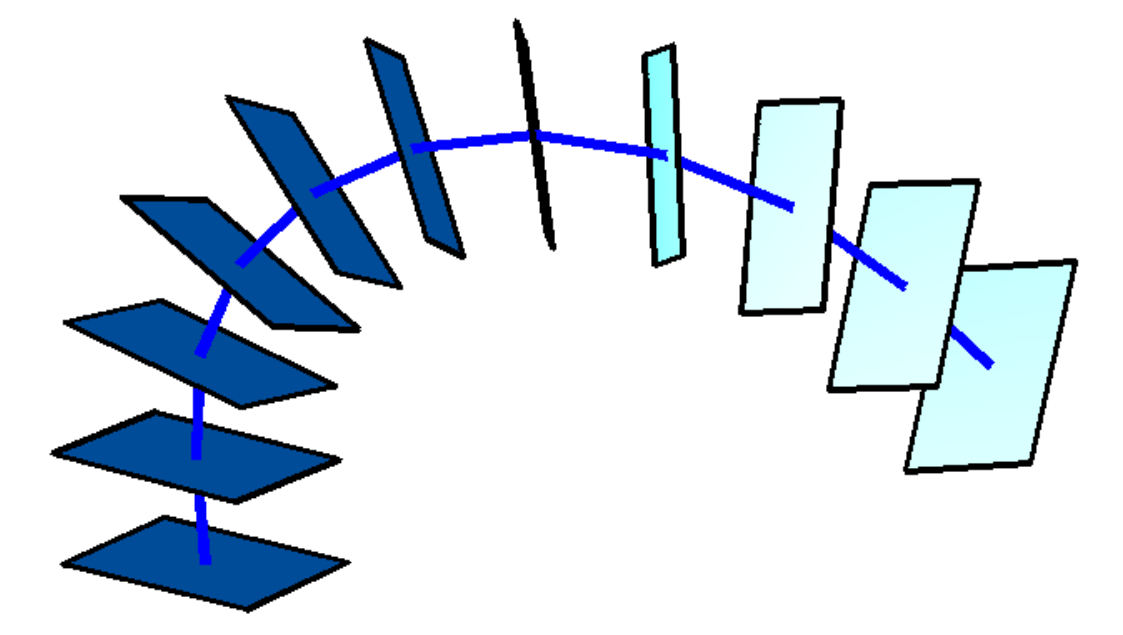}
  \vspace{-3pt}
  \caption{\footnotesize The geometrically exact beam model defined by the position $\mathbf{r} \in \mathbb{R}^3$ of the line of centroids and by the attitude matrix 
    $\Lambda \in SO(3)$ of each cross section. }  
  \label{fig:model_explanation}  
\end{figure}

Consider a beam (Figure~\ref{fig:model_explanation}) with length $L=1 \,m$ and square cross-section with side $a=0.01\,m$ that is free of tractions and body forces. A mesh size of $\Delta s=0.1$ and time step $\Delta t =0.0005$ are chosen with total simulation time $T = 3\,s$. The beam parameters are set to: $\rho = 10^3 \,kg/m^3$, $M=10^{-1}\,kg/m$, $E=5.10^3 \,N/m^2$, $\nu=0.35$. 

\medskip

We first consider the time-integration of the beam followed by space-reconstruction.

\paragraph{Time evolution.} 
The initial conditions are given by the configurations $g_a^0=(\Lambda_a^0, \mathbf{r}_a^0)$ and the initial speed $\xi_a^0=(\Omega_a^0, \Gamma_a^0)$ at time $t=t^0$ for all positions $s_0,...,s_A$. Given $\tau$ as defined 
in \eqref{cay}, we choose 
\[
g_0^0=(\text{Id},(0,0,0)), \quad g_{a+1}^0 = g_a^0\,
\tau(\Delta s \, \eta_a^0), \quad \text{for all} \ a=0, ..., A-1,
\]
where $\eta_a^0 =( 1,1.5,1,0,0,1)$, for all $a=0, ..., A-1$, and
\[
\xi_a^0 = \frac{1}{\Delta t} \tau^{-1}\left((g_a^0)^{-1}g_a^1\right), \quad \text{for all} \ a=0, ..., A-1,
\]
where $g_0^1=(\text{Id}, (0,0, \Delta t))$, and $g_{a+1}^1= g_a^1 \,
\tau(\Delta s \, \eta_a^1)$, for all $a=0, ..., A-1$, with $\eta_a^1=(1.004, 1.52, 1.005, -0.01, 0, 1)$.
For this problem, the discrete zero-traction boundary conditions (i.e., the discrete version of \eqref{cont_boun_cond}) are  imposed at the two extremities of the beam.
The implemented scheme is \eqref{CDEL_beam} with boundary conditions \eqref{boun_cond1} at the extremities.

The algorithm produces the configurations $\mathbf{g}^1,..., 
\mathbf{g}^{N-1}$,  where $\mathbf{g}^j= (g_0^j,...,g_{A-1}^j)$, plotted in Figure \ref{time_evolution} below.
\begin{figure}[H]
\centering
\begin{center}
\begin{tabular}{cc}
\includegraphics[width=1.1 in]{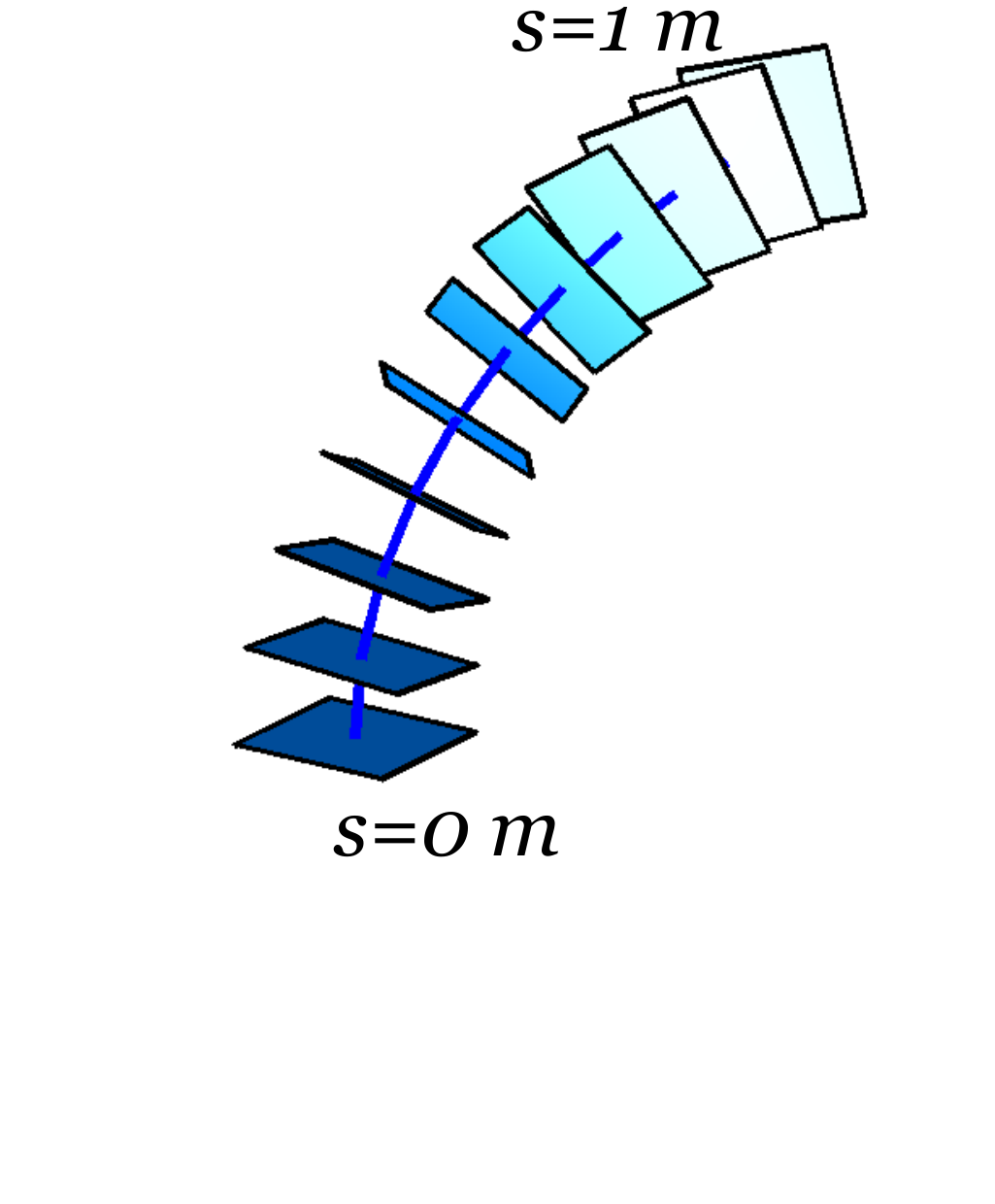}  
\includegraphics[width=1.1 in]{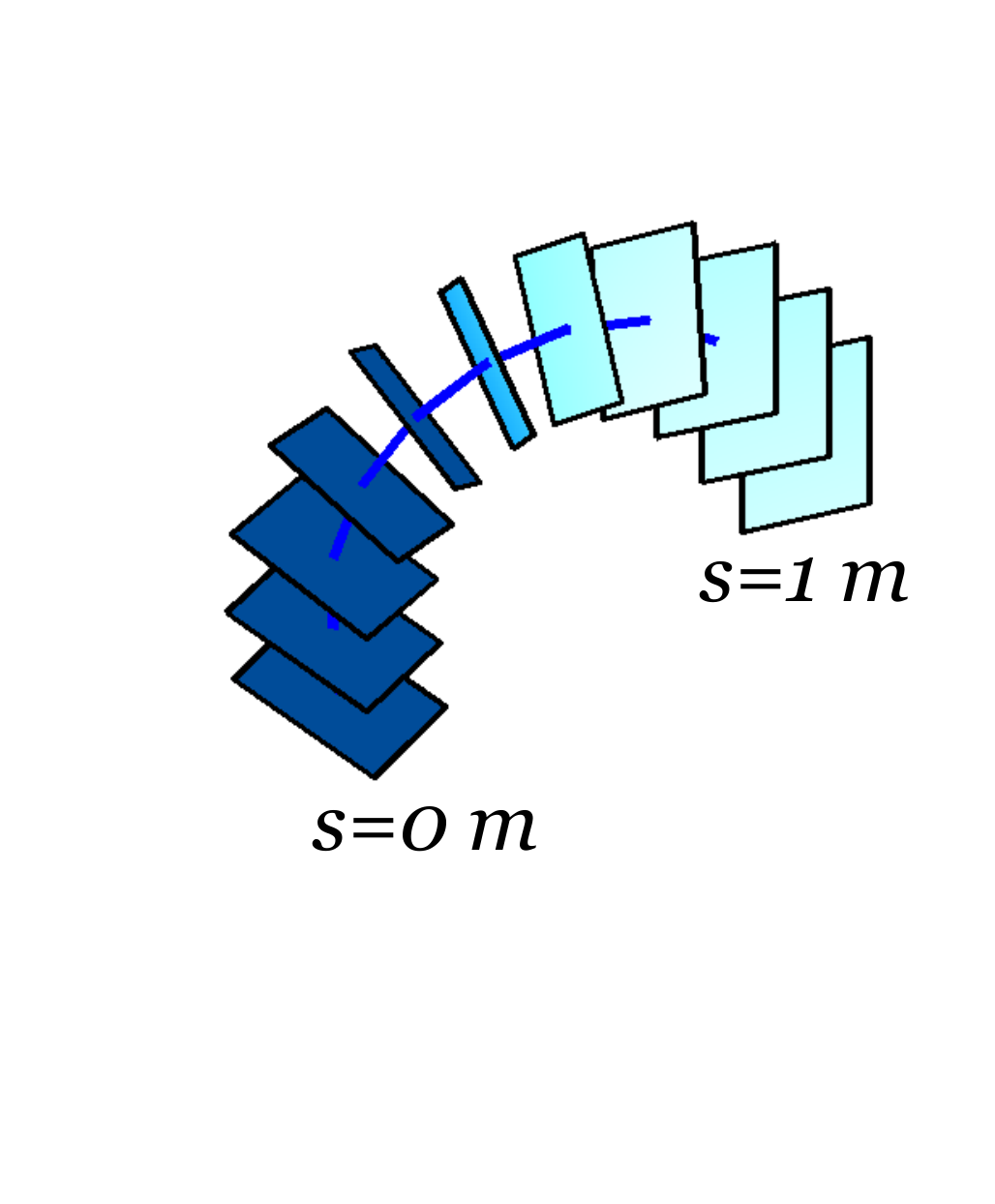}  
\includegraphics[width=1.1 in]{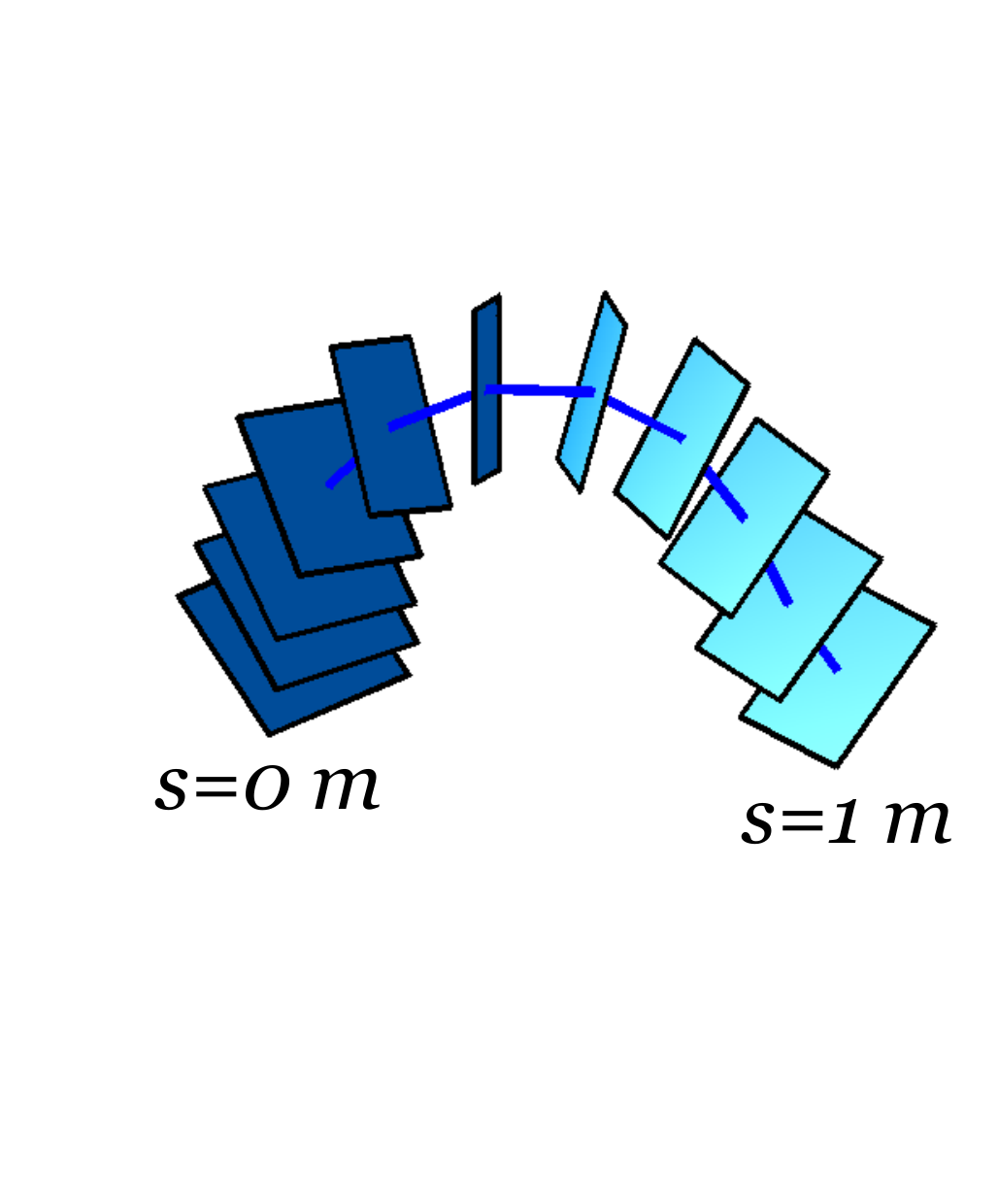}    
\includegraphics[width=1.1 in]{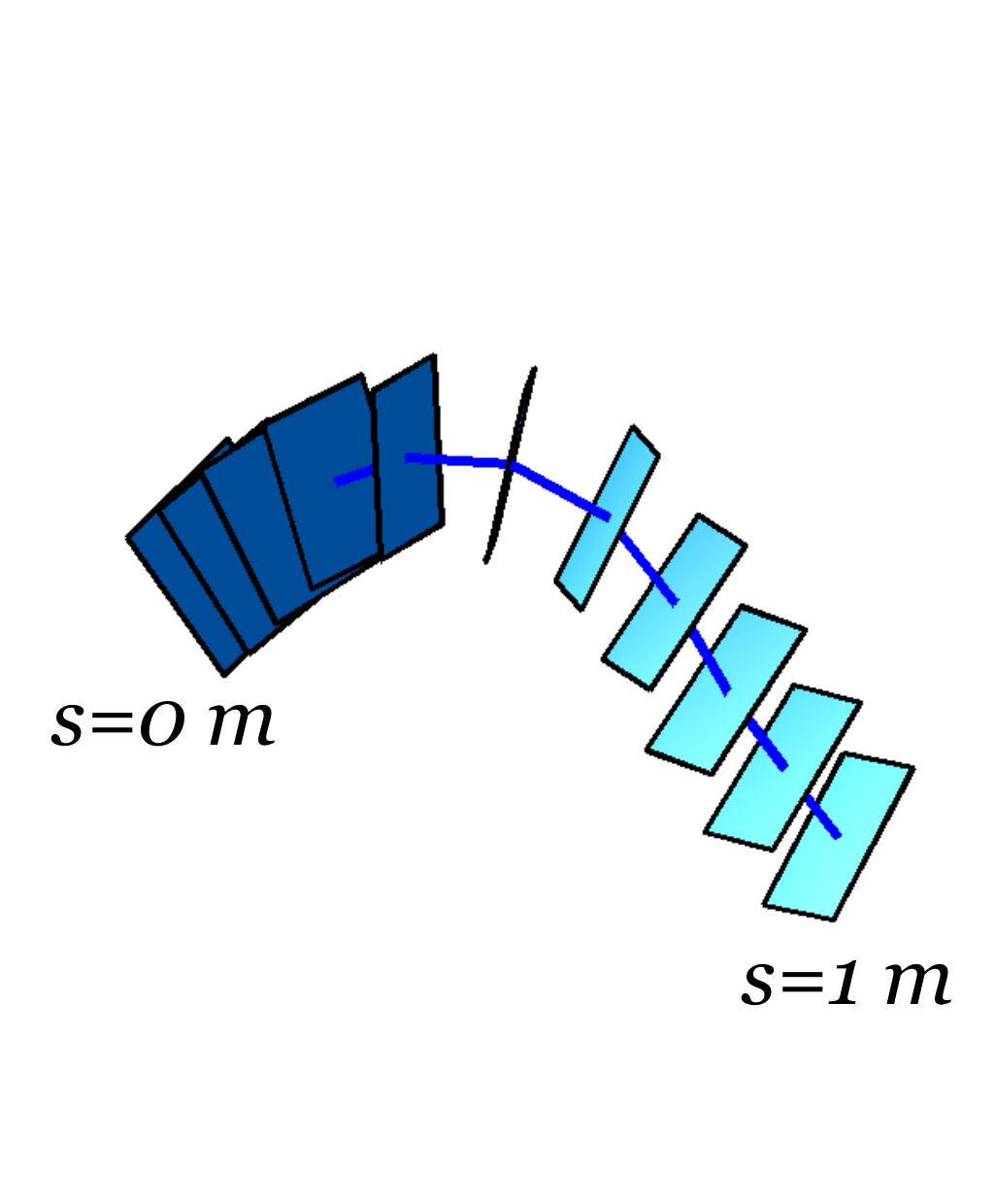}  
\includegraphics[width=1.1 in]{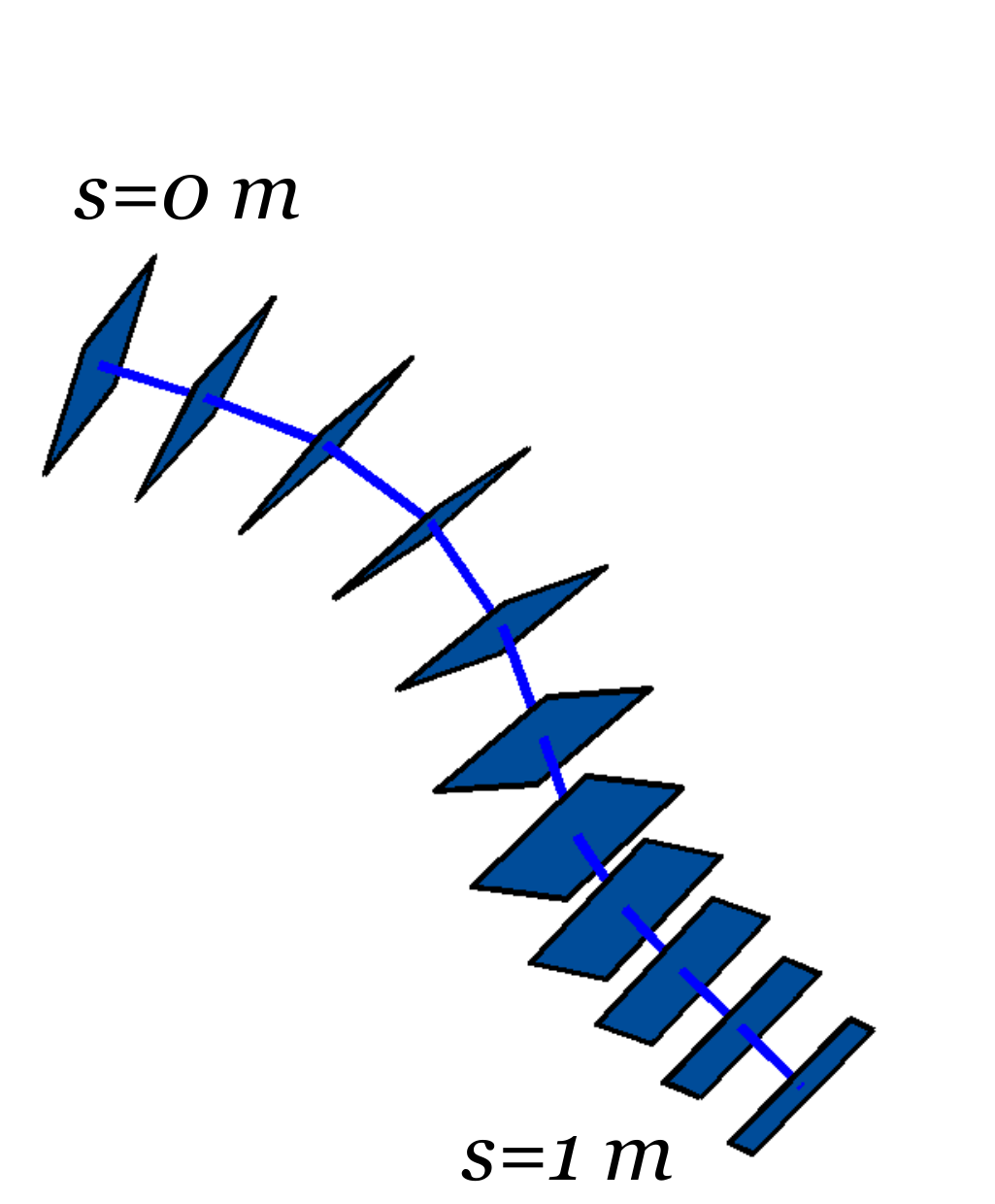}
\end{tabular}
\caption{\footnotesize  Each figure represents the space evolution $\mathbf{g}^j =\{g_a ^j, \ a=1,...,A\} $ with $s_A=1m$, at a given time evolution $t$ of the beam.
The chosen times correspond to $t=0\,s$, $0.1\,s$, $0.2\,s$, $0.3\,s$, $0.5\,s$.  }\label{time_evolution}  
  \end{center}
\end{figure}

\noindent\textit{Energy behavior.} The above DCEL equations, 
together with the boundary conditions, are equivalent to the 
DEL equations for $\mathsf{L}_d( \mathbf{g}^j , \boldsymbol{\xi}^j)$ 
in \eqref{model_L_d}; see the discussion in 
\S\ref{subsec:discrete_Lagrangian}. In particular, the solution 
of the discrete scheme defines a discrete symplectic flow in 
time $( \mathbf{g}^j   , \boldsymbol{\xi }^j ) \mapsto 
(\mathbf{g}^{j+1}   , \boldsymbol{\xi }^{j+1} )$. As a consequence, the energy $\mathsf{E}_{\mathsf{L}_d}$ (see \eqref{formulas_reconstruction_space}) of the Lagrangian $\mathsf{L}_d$ associated to the temporal evolution description is approximately conserved, as illustrated in Fig.~\ref{energy_momentum_time} left, below.

\medskip

\noindent\textit{Momentum map conservation.} Since the discrete Lagrangian density is $SE(3)$-invariant, the discrete covariant 
Noether theorem
$\mathscr{J}_{B,C}^{K,L}(g _d  ) = 0$ is verified; see \S\ref{DCMP}. 
Since the discrete Lagrangian $\mathsf{L}_d$ is $SE(3)$-invariant, the 
discrete momentum maps coincide: $\mathsf{J}^+_{\mathsf{L}_d}=
\mathsf{J}^-_{\mathsf{L}_d}= \mathsf{J}_{\mathsf{L}_d}$, and we have
\[
\mathsf{J} _{\mathsf{L}_d}(\mathbf{g}^j , \boldsymbol{\xi} ^j )
= 
\sum_{a=0}^{A-1}  \Delta s \mathrm{Ad}^*_{(g_a^j)^{-1}} \mu_a^j;
\]
see  \eqref{time_momentum_map} and \eqref{discrete_momentum_map}. 
In view of the boundary conditions used here, it follows that the 
discrete momentum map $\mathsf{J} _{\mathsf{L}_d}$ is exactly 
preserved as illustrated in Fig.~\ref{energy_momentum_time} right, consistently with Theorem \ref{summary_Noether}. This can be seen as a consequence of the covariant discrete Noether theorem $\mathscr{J}_{0,A-1}^{K,L}( g _d )=0$. We also checked numerically that the discrete covariant Noether theorem \eqref{DCN} is verified. For example, for $B=K=0$, $C=A-1$, $L=N-1$, we found
\[
\sum_{j=1}^{N-1}\Delta t \left( \mathrm{Ad}^*_{(g_0^{j})^{-1}} \lambda_0^{j} - \mathrm{Ad}^*_{(g_{A-1}^{j})^{-1}} \lambda_{A-1}^j \right) +\Delta s\left(-\mathrm{Ad}^*_{(g_0^{N-1})^{-1}} \mu_0^{N-1} +\mathrm{Ad}^*_{(g_0^0)^{-1}} \mu_0^0\right) = \mathbf{0},
\]
up to round-off error.
\begin{figure}[H]
\centering
\begin{center}
\begin{tabular}{cc}
 \includegraphics[width=1.7 in]{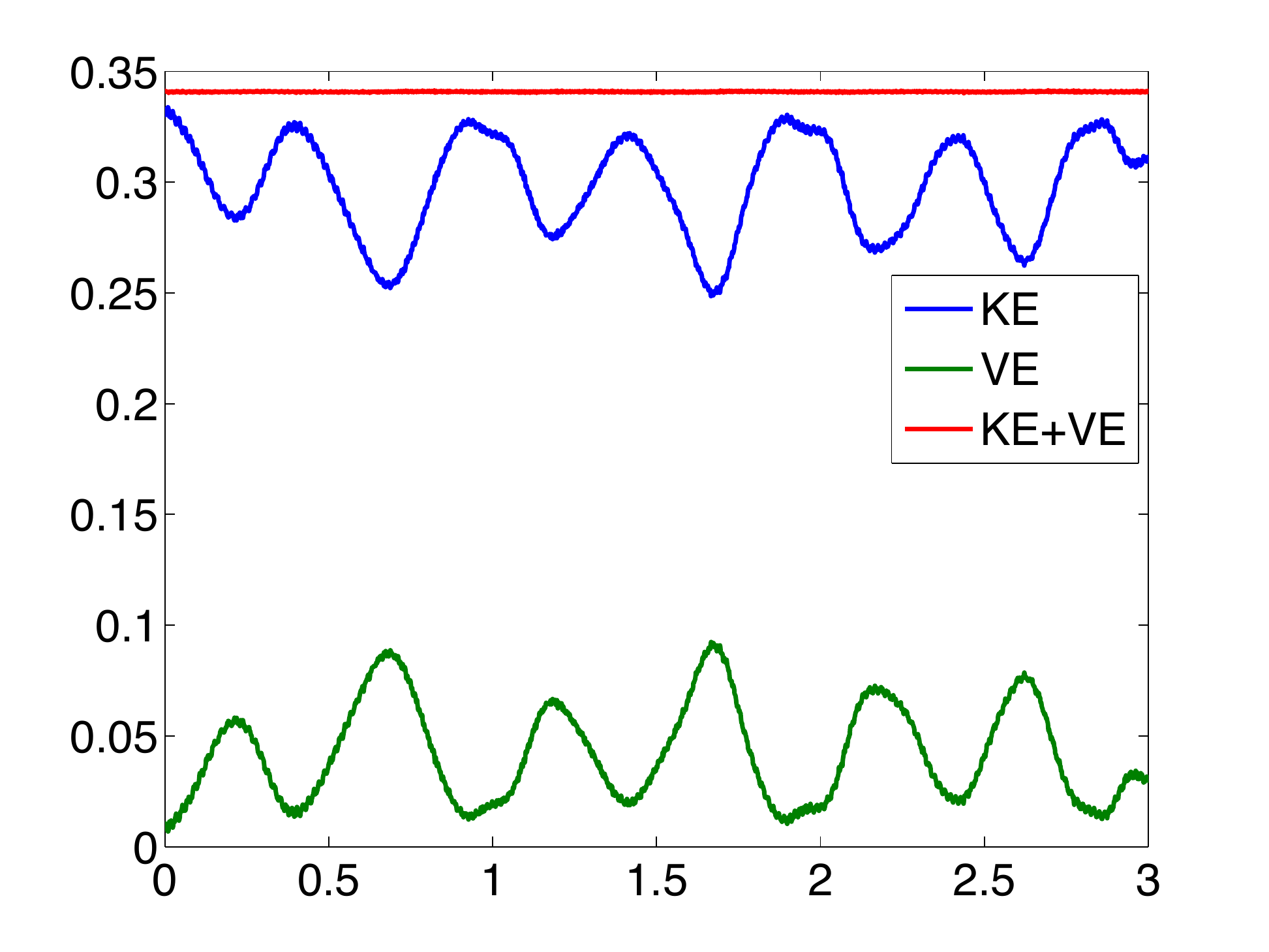}    \qquad   \includegraphics[width=1.7 in]{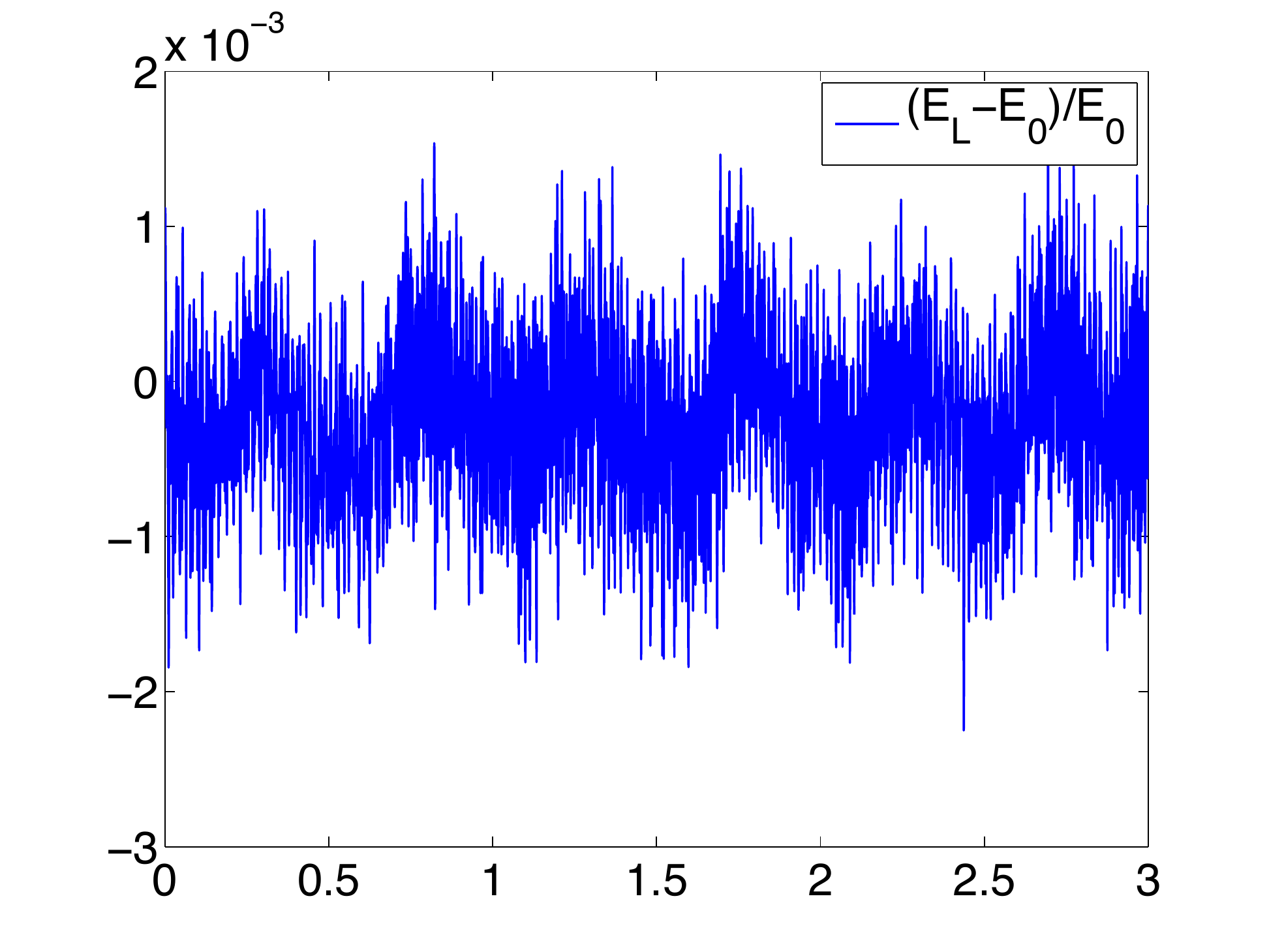}  
\end{tabular}
\vspace{-3pt}
\caption{\footnotesize Left: total energy behavior $\mathsf{E}_{\mathsf{L}_d}$. Right: relative error $(\mathsf{E}_{\mathsf{L}_d}(t^j) -\mathsf{E}_{\mathsf{L}_d}(t^0))/\mathsf{E}_{\mathsf{L}_d}(t^0)$. 
Both, during a time interval of $3$s.}\label{energy_momentum_time} 
\end{center}  
\end{figure}

\begin{figure}[H]
\centering
\begin{center}
\begin{tabular}{cc}
 \includegraphics[width=1.7 in]{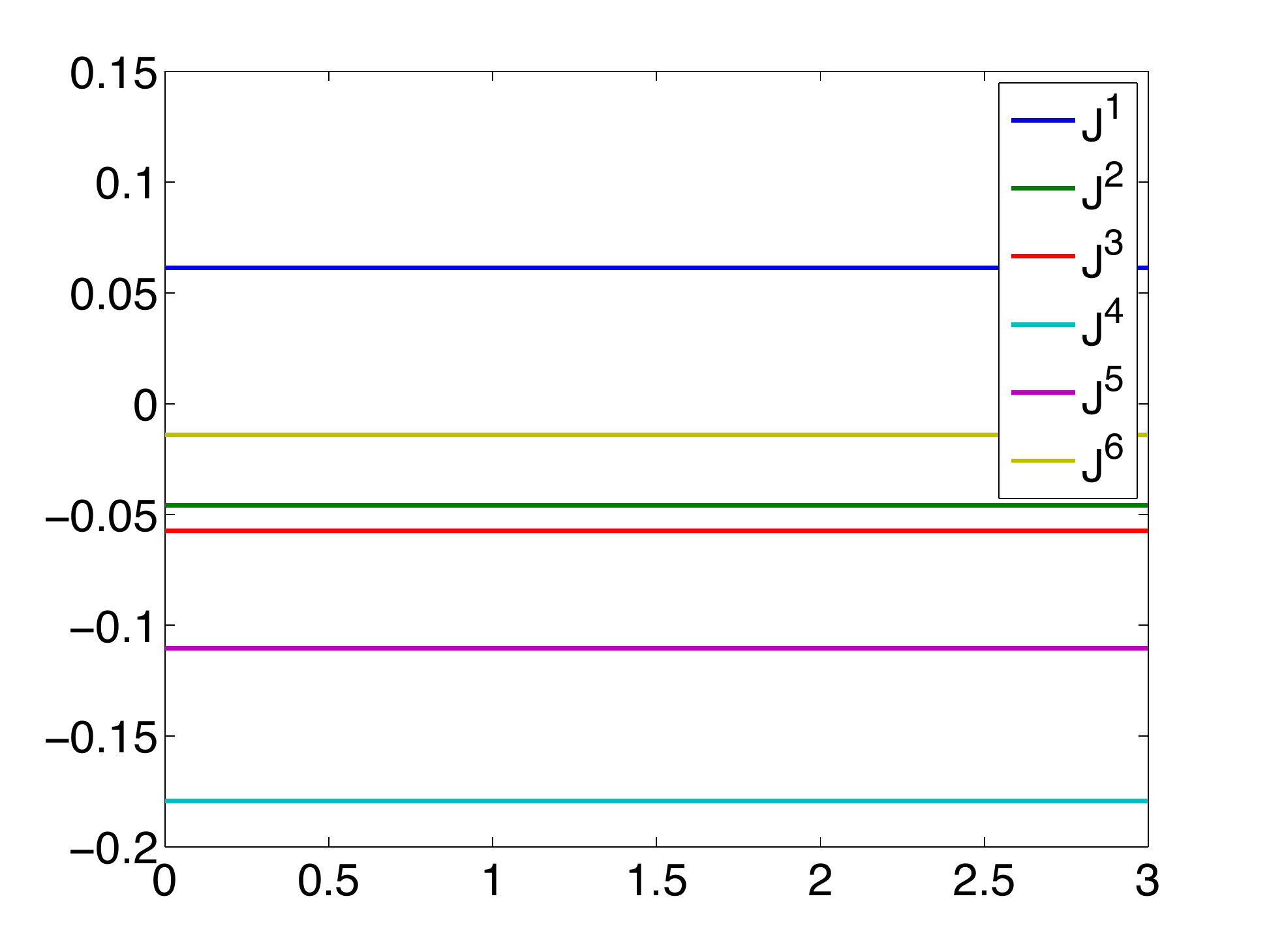}    \qquad   \includegraphics[width=1.7 in]{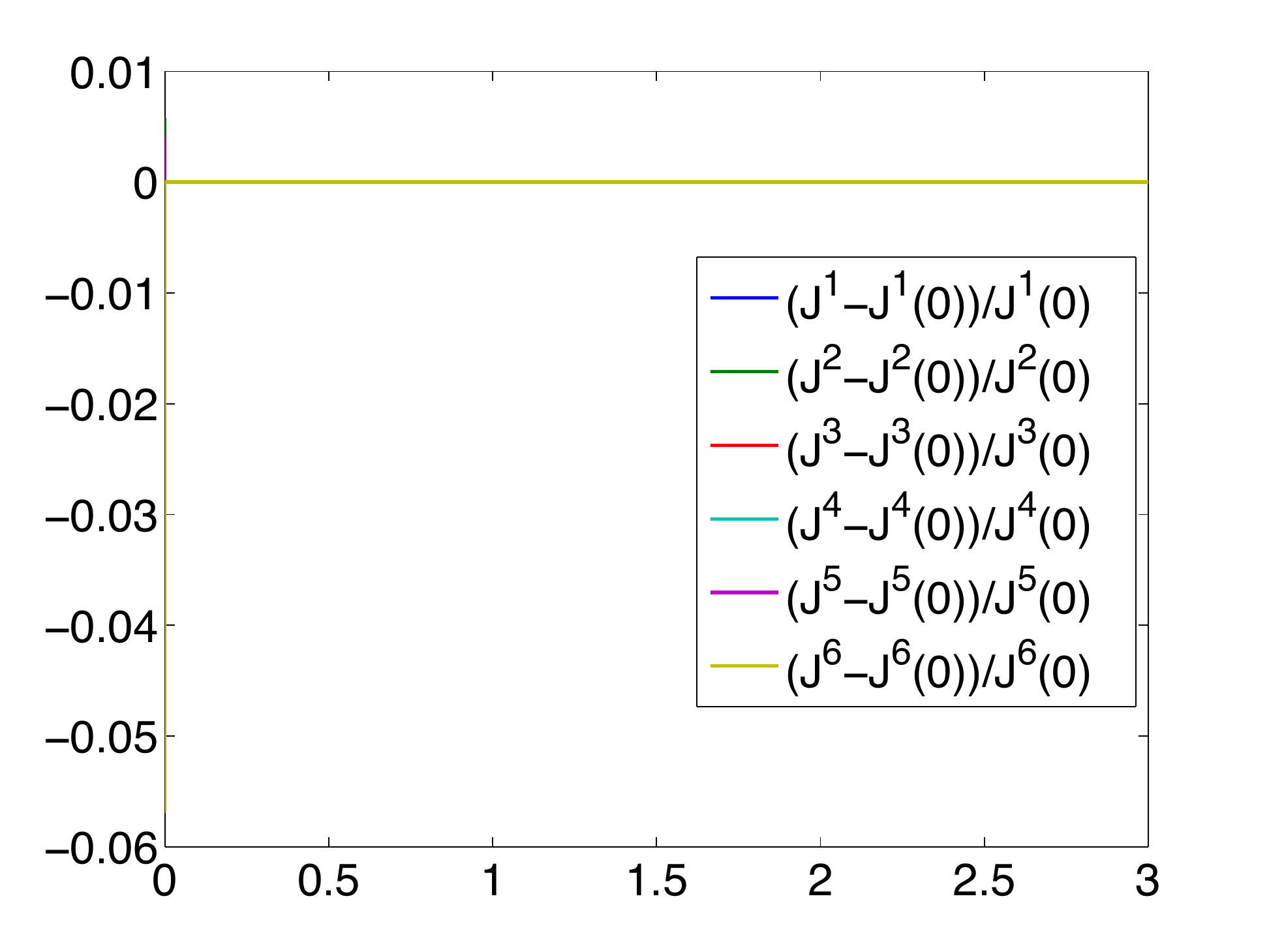} 
\end{tabular}
\vspace{-3pt}
\caption{\footnotesize Left: conservation of the discrete momentum map $\mathbf J_{\mathsf{L}_d}=(\mathsf{J}^1,...,\mathsf{J}^6) \in \mathbb{R}^6$. Right: relative error $(\mathbf J_{\mathsf{L}_d}(t^j)- \mathbf J_{\mathsf{L}_d}(t^0))/\mathbf J_{\mathsf{L}_d}(t^0)$.
Both, during a time interval of $3$s. }\label{energy_momentum_time} 
\end{center}  
\end{figure}

\begin{figure}[H]
\centering
\begin{center}
\begin{tabular}{cc}
 \includegraphics[width=3.8 in]{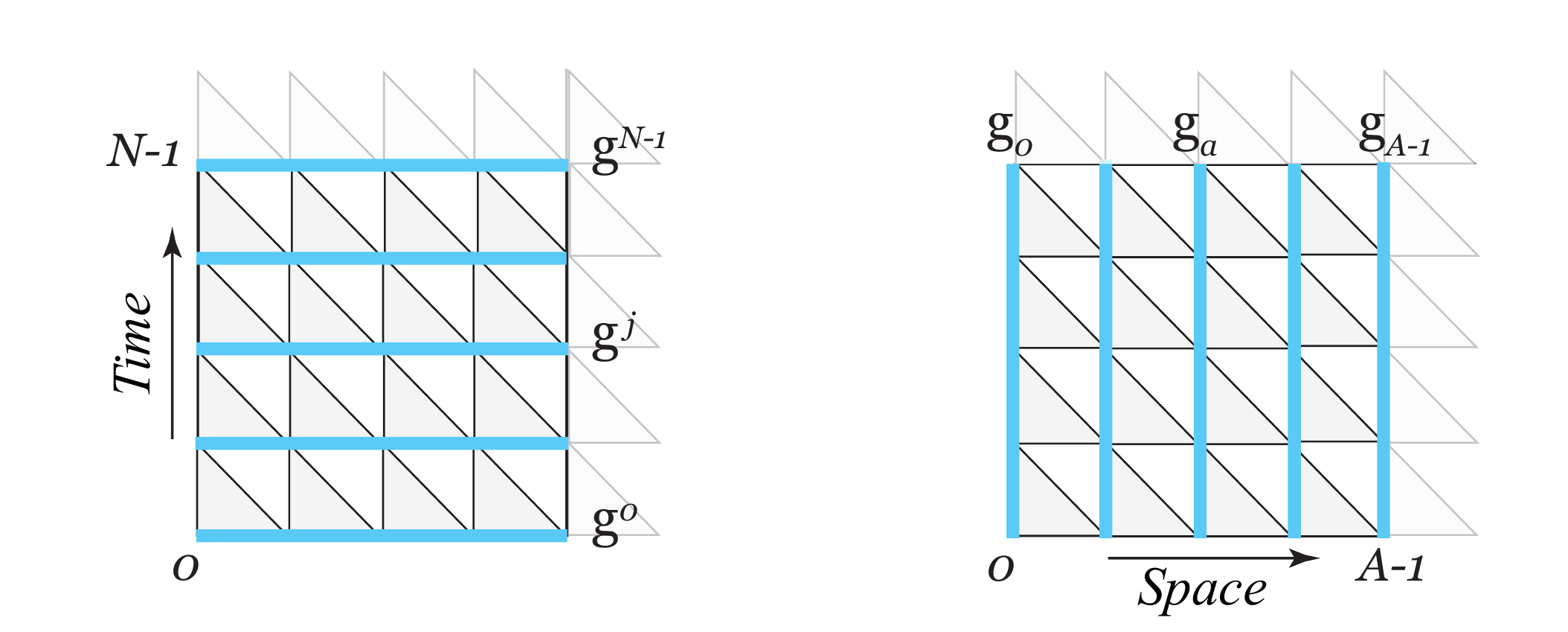}      
\end{tabular}
\vspace{-3pt}
\caption{\footnotesize Integration can be performed either in time-direction $\mathbf{g}^1 \rightarrow \cdots \rightarrow \mathbf{g}^{N-1}$ (left) or in space direction $\mathbf{g}_0 \rightarrow \cdots \rightarrow \mathbf{g}_{A-1}$ (right). } \label{mesh_recons}
\end{center}  
\end{figure}

\begin{figure}[H]
\centering
\begin{center}
\begin{tabular}{cc}
\includegraphics[width=1.1 in]{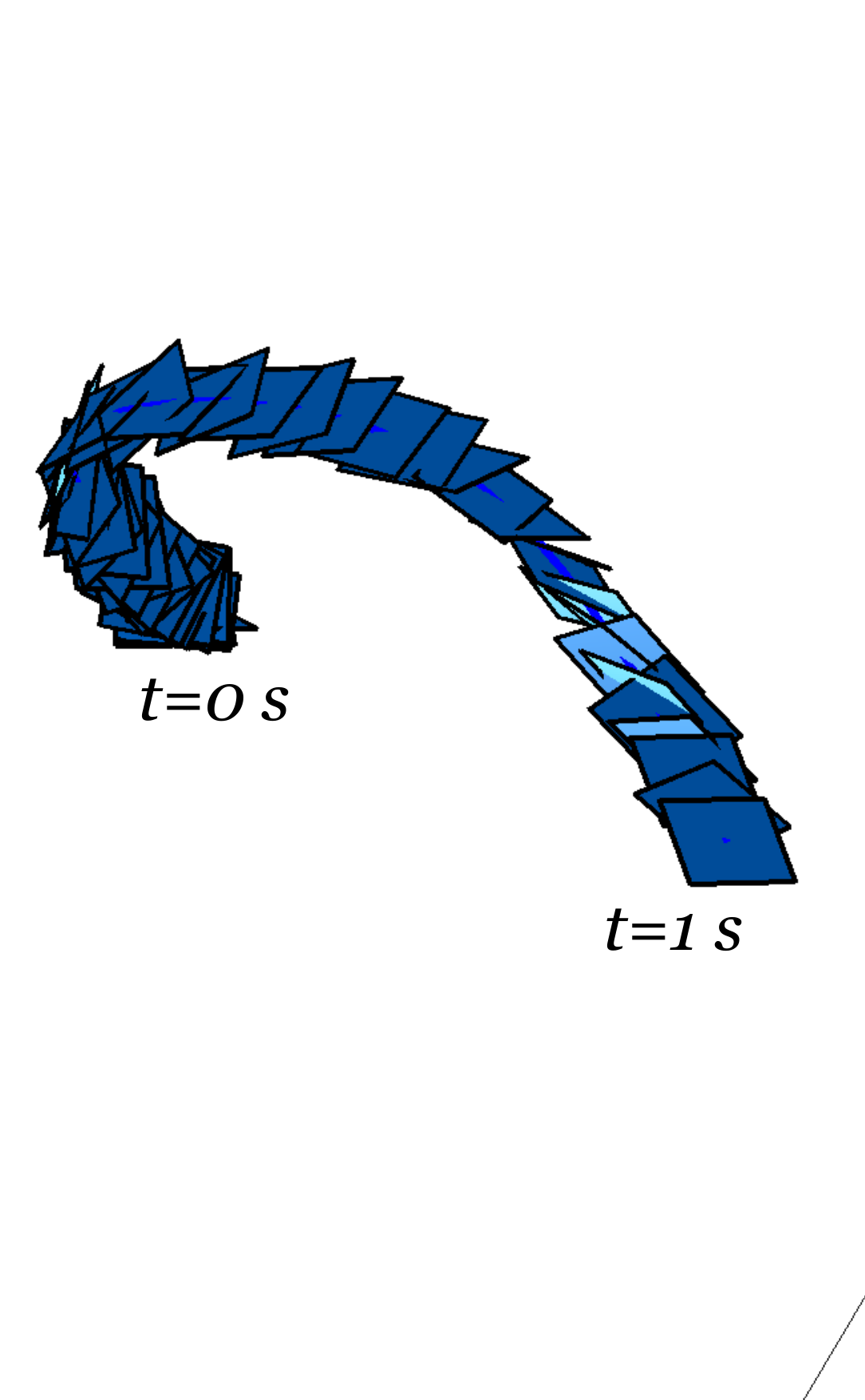} 
\qquad 
 \includegraphics[width=1.1 in]{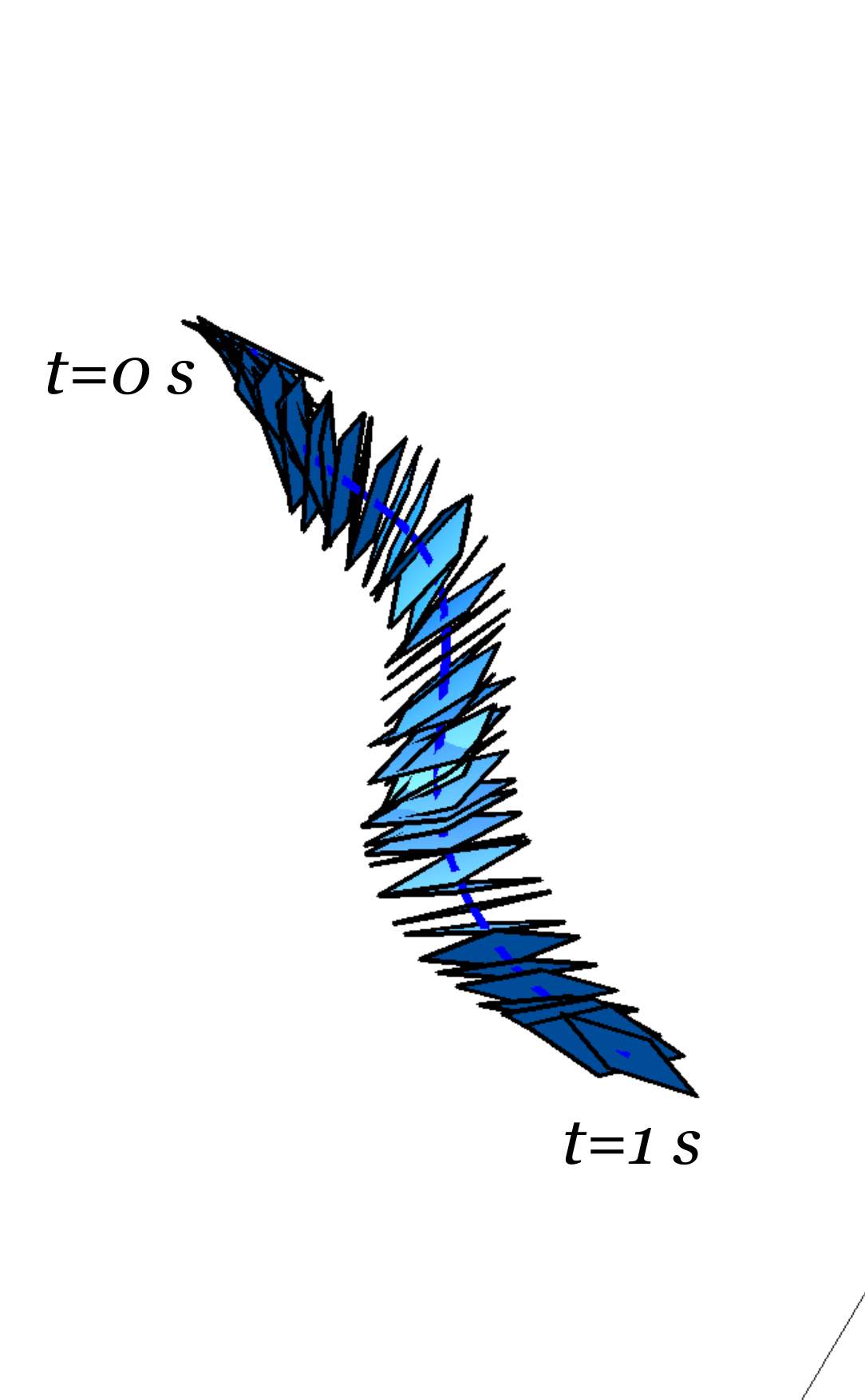}  
\qquad
\includegraphics[width=1.1 in]{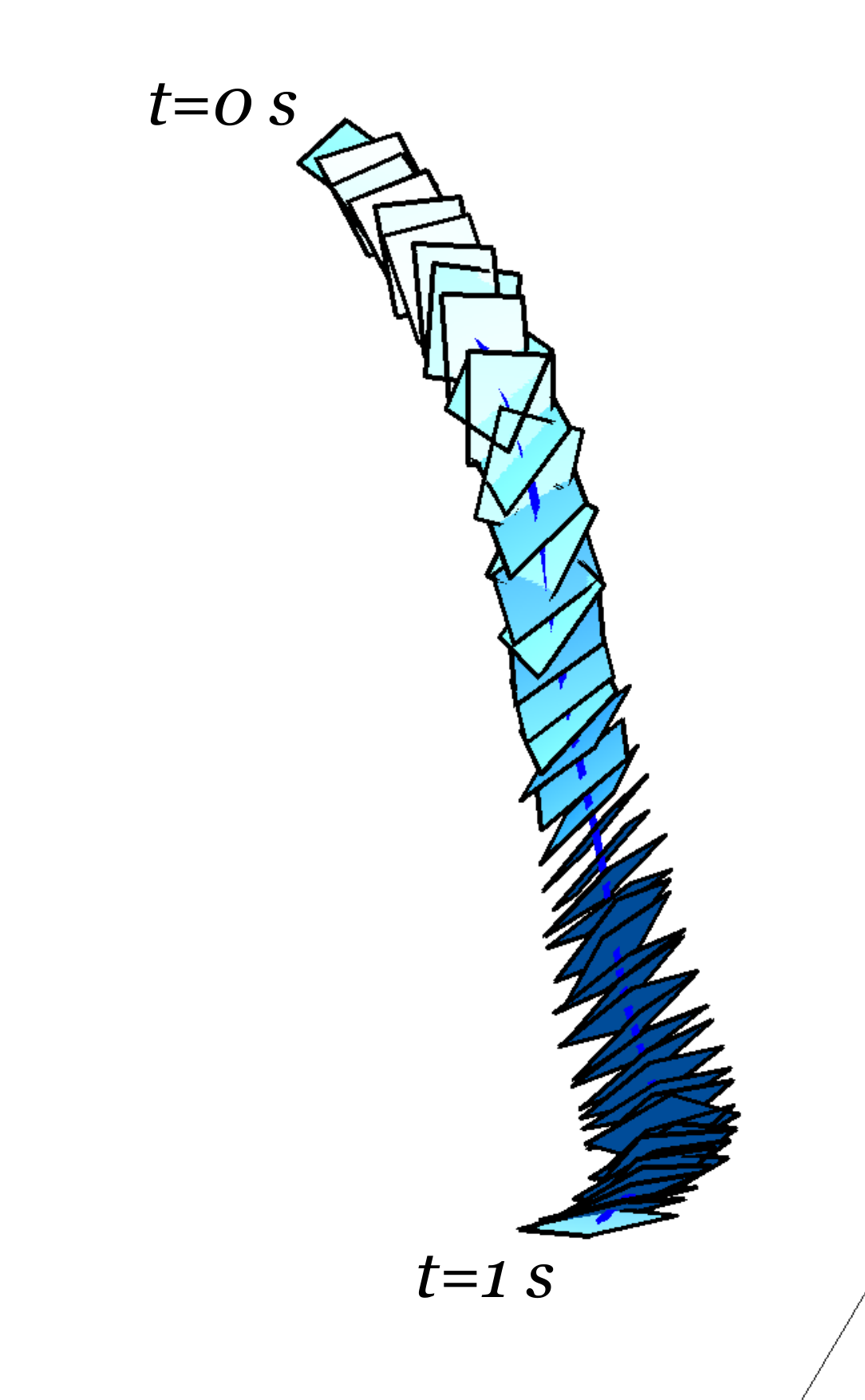} 
\qquad 
\includegraphics[width=1.1 in]{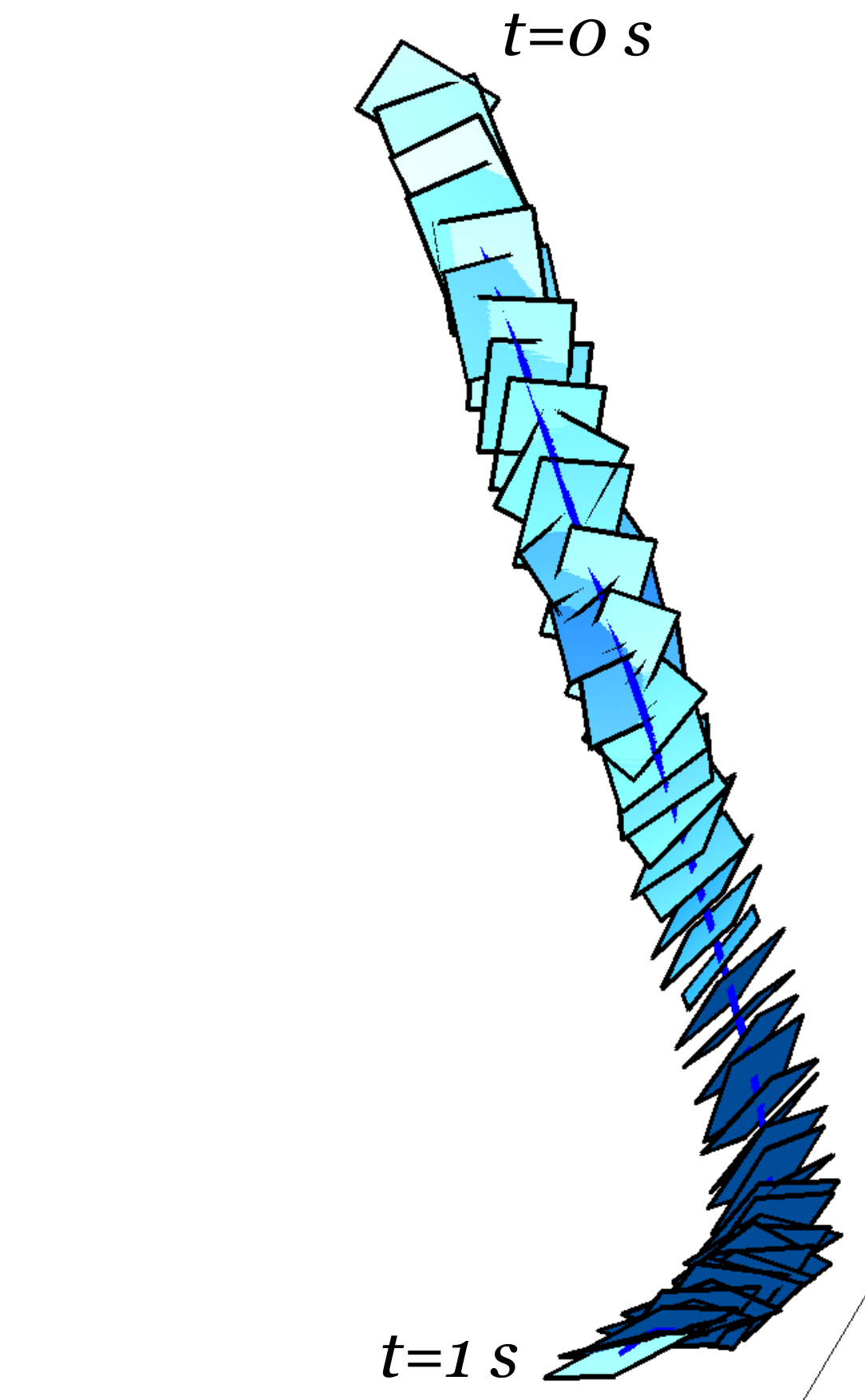}
\end{tabular}
\caption{\footnotesize  
From left to right: space-reconstruction of the trajectories in time of the beam sections at $s=0.1\,m, \,  0.5\,m, \,  0.8\,m,\, 1\,m$. }\label{reconstruction_time_space}  
\end{center}
\end{figure}

\paragraph{Reconstruction.} The above computed time evolution provides a set of configurations $\mathbf{g}^1,..., \mathbf{g}^{N-1}$  for the duration of $1\,s$. Repackaging these results to emphasize spatial
evolution yields $\mathbf{g}_0,..., \mathbf{g}_{A-1}$, where 
$\mathbf{g}_a= (g_a^0,...,g_a^{N-1})$ (see 
Figure~\ref{mesh_recons}). The actual reconstructed motions are 
shown in Figure~\ref{reconstruction_time_space}.

\subsection{Space-integration and time-reconstruction}\label{space_setting} 
 
The problem treated here corresponds to the following situation. 
We assume that we know the evolution (for all $ t\in [0,T]$) of 
one of the extremities, say $s=0$, as well as the evolution of 
its strain (for all $ t\in [0,T]$). We also assume that at the 
initial and final times $t=0,T$, the velocity of the beam is zero. 
This corresponds to zero-momentum boundary conditions. The spatial 
configuration of the beam at $t=0$ and $t=T$ is, however, unknown. 
The approach described in this paper, that makes use of both the 
temporal and spatial evolutionary descriptions in both the 
continuous and discrete formulations, is especially well designed 
to discretize this problem in a structure preserving way.

Note that we do not impose the zero-traction boundary conditions \eqref{cont_boun_cond}.
 
\subsubsection{Scenario A}

In this example, the mesh is defined by the space step $\Delta s=0.05$ and the time step $\Delta t =0.05$. The total length of the beam is $L=0.8\,m$ and the total simulation time is $T = 10\,s$. The characteristics of the material are: $\rho = 10^3 \,kg/m^3$, $M=10^{-1}\,kg/m$, $E=5.10^4 \,N/m^2$, $\nu=0.35$.
  
\begin{figure}[H]
\centering
\begin{center}
\begin{tabular}{cc}
\includegraphics[width=2 in]{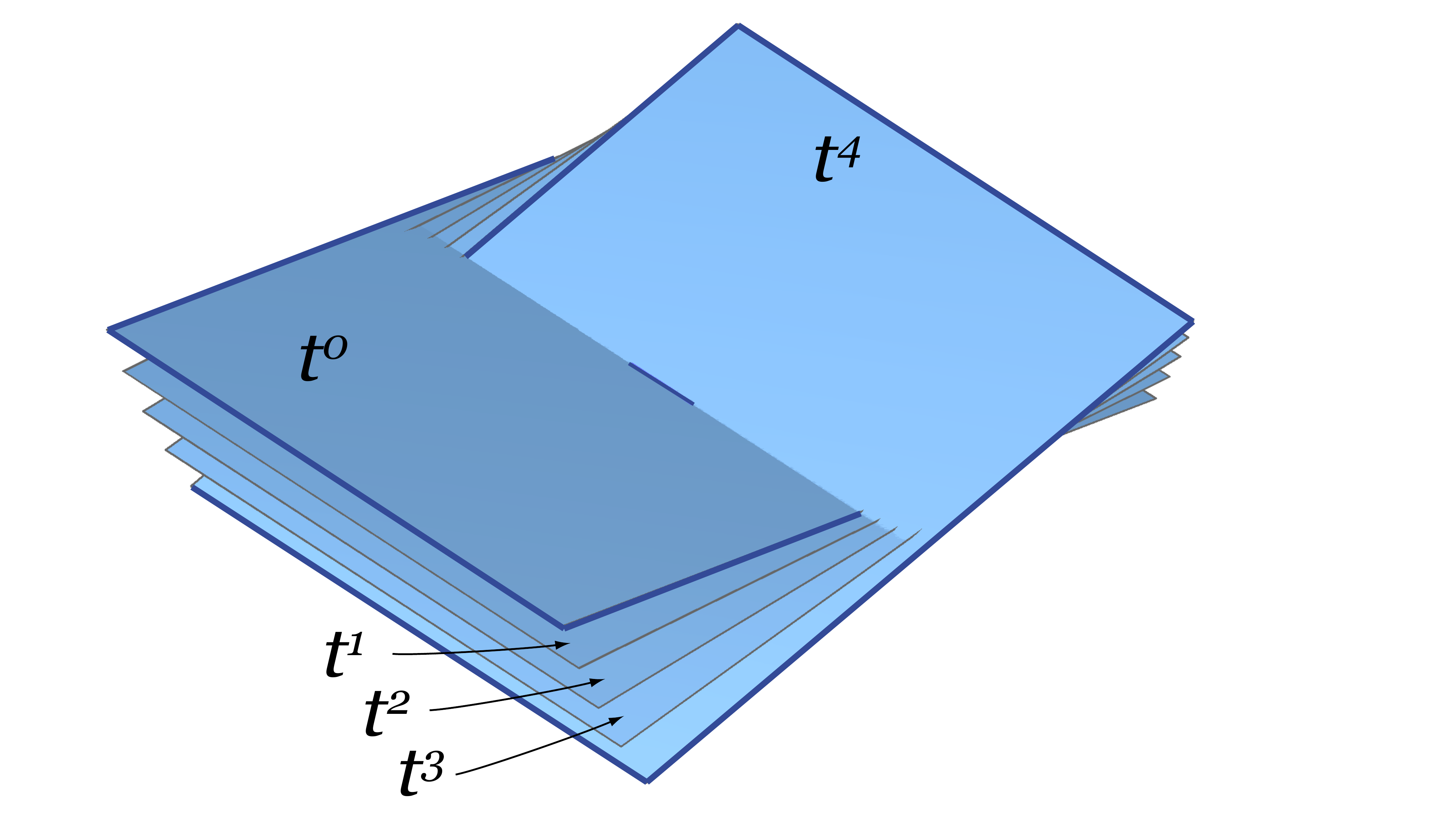} 
\end{tabular}
\caption{\footnotesize  
Initial conditions $\mathbf{g}_0^j$ (only five time-slices at $j\in\{0,1,2,3,4\}$ shown).}
\label{initial_condition}  
\end{center}
\end{figure}

\paragraph{Space-integration.} The initial conditions are given by the configuration $\mathbf{g}_0$ and the initial strain $\boldsymbol{\eta}_0$ at the extremity $s=0$.  We choose the following configuration and strain:
\[
g_0^0=(\mathrm{Id}, (0,0,0)),\quad g_0^{j+1} = g_0^j\, \tau
(\Delta t \xi _0 ^j ),\quad \text{for all $j=0,...,N-1$},
\]
where $ \xi _0 ^j = (0,-2,0,0,-0.1,0)$, for all $j=0,...,N-1$, and
\[
\eta_0^j= \frac{1}{\Delta s} \, \tau^{-1}\left((g_0^j)^{-1} g^j_1 \right),\quad \text{for all $j=0,...,N-1$},
\]
where $g^0_1  =(\mathrm{Id}, (0,0,\Delta s))$ and $g^{j+1}_1= g^j_1 \,
\tau(\Delta t \xi ^j _1 )$, for all $j=0,...,N-1$, with 
$\xi ^j _1  = (0.007,-1.998, -0.007,-0.08,-0.1,0)$, see Fig.~\ref{initial_condition}. For this problem, the discrete zero-momentum boundary conditions are imposed. The implemented scheme is
\eqref{CDEL_beam} with the boundary conditions
 \eqref{boun_cond2} at the temporal extremities.

\begin{figure}[H]
\centering
\begin{center}
\begin{tabular}{cc}
\includegraphics[width=1.3 in]{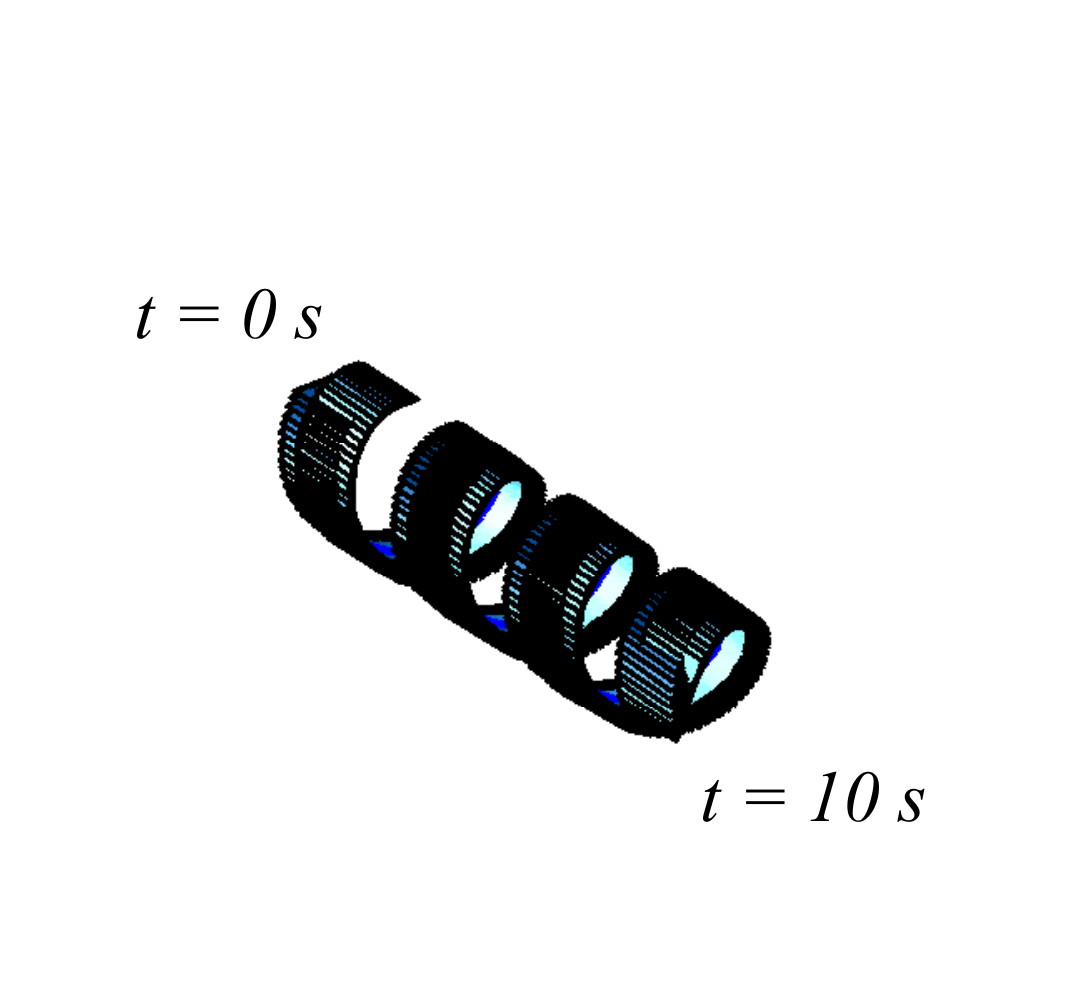}  
\includegraphics[width=1.3 in]{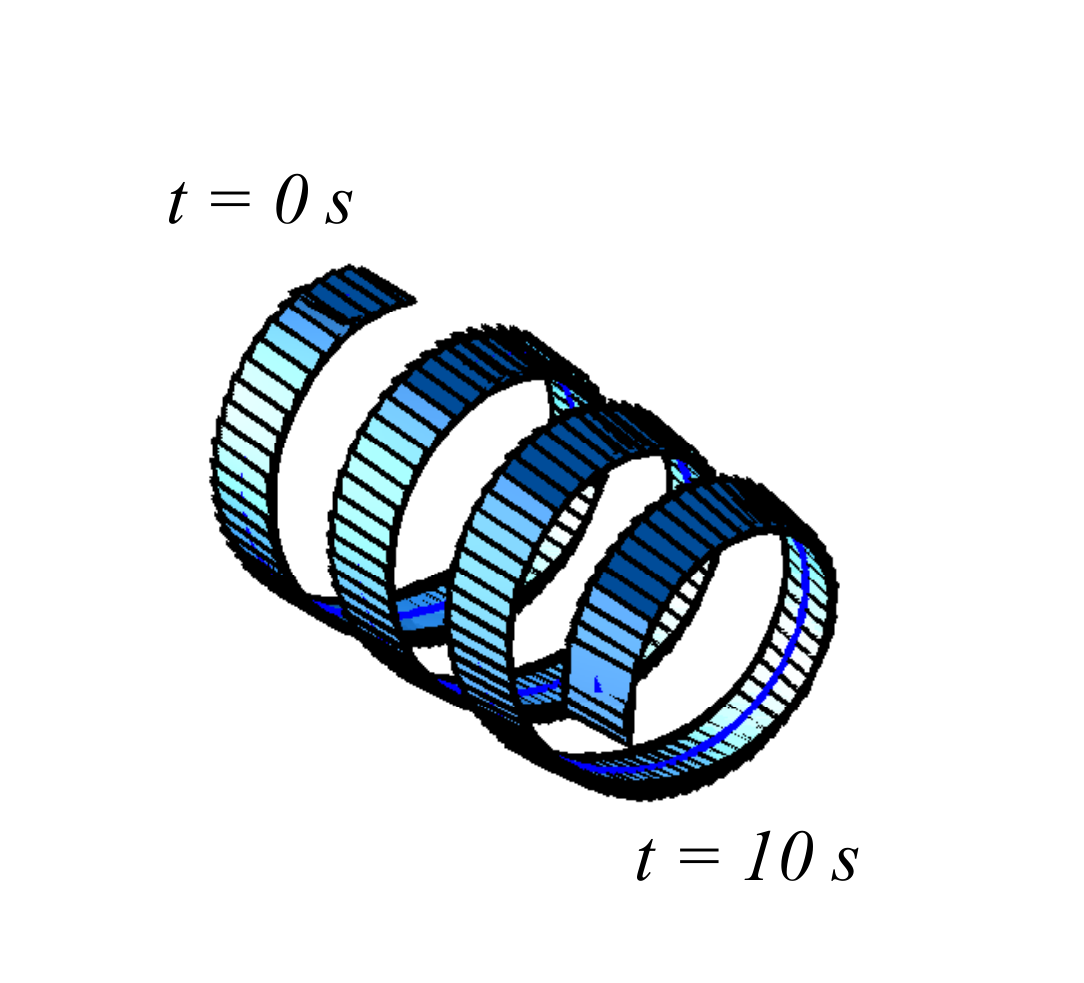}  
\includegraphics[width=1.3 in]{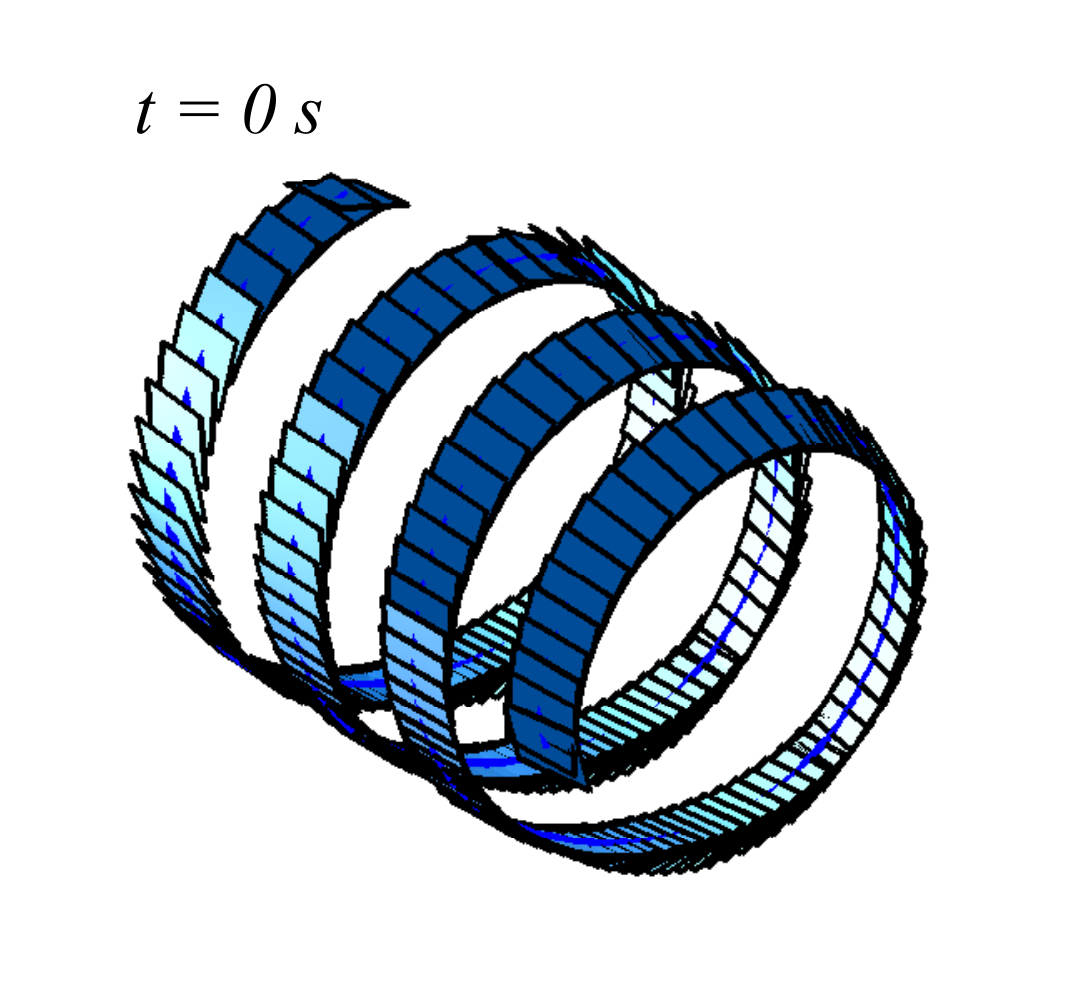}     \includegraphics[width=1.3 in]{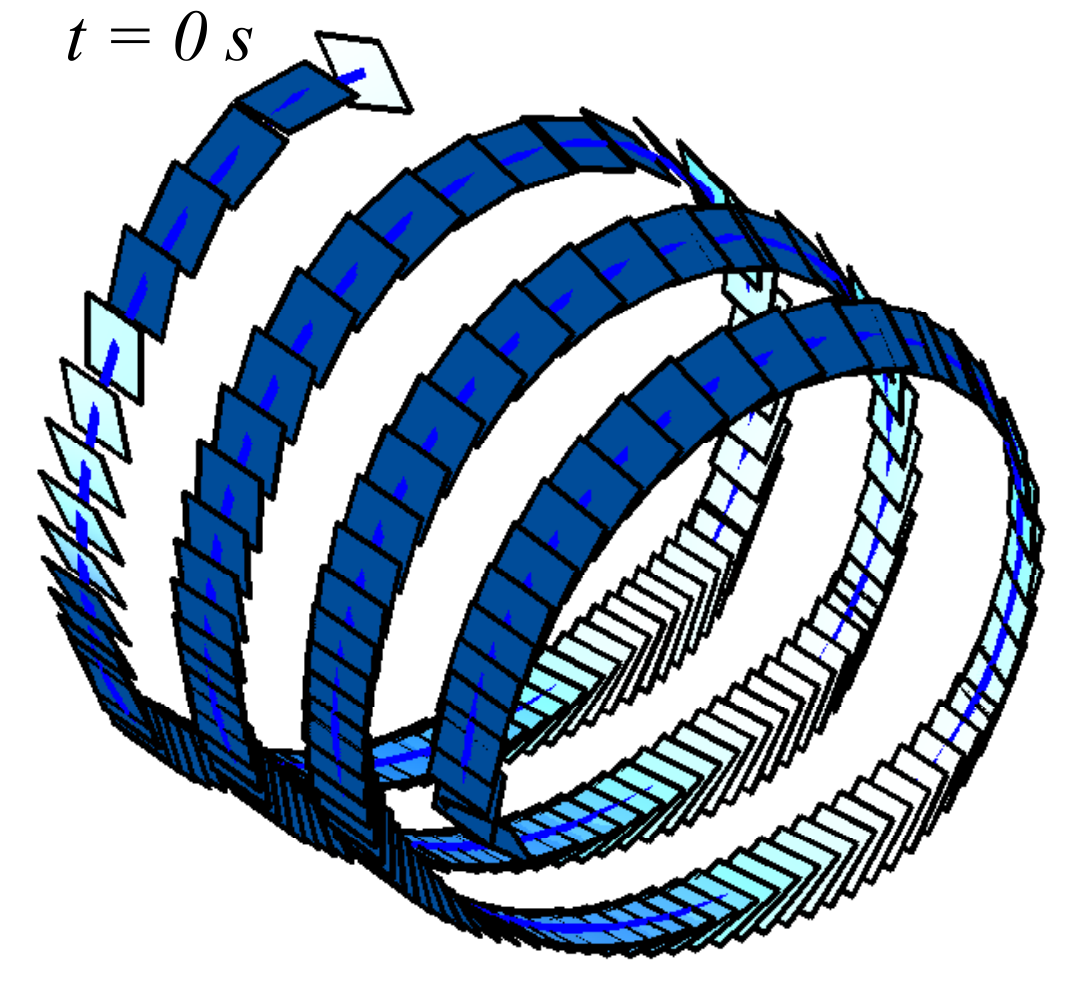}
\end{tabular}
\caption{\footnotesize  
Each figure represents the time evolution $\mathbf{g}  _a =\{g_a ^j, j=1,...,N-1\} $ of a given node $a$ of the beam, with $t^{N-1}=10s$.
The chosen nodes correspond to $s=0\,m$, $0.3\,m$, $0.55\,m$, $0.8\,m$. }\label{space_evolution1}  
\end{center}
\end{figure}

\begin{figure}[H]
\centering
\begin{center}
\begin{tabular}{cc}
\includegraphics[width=2.2 in]{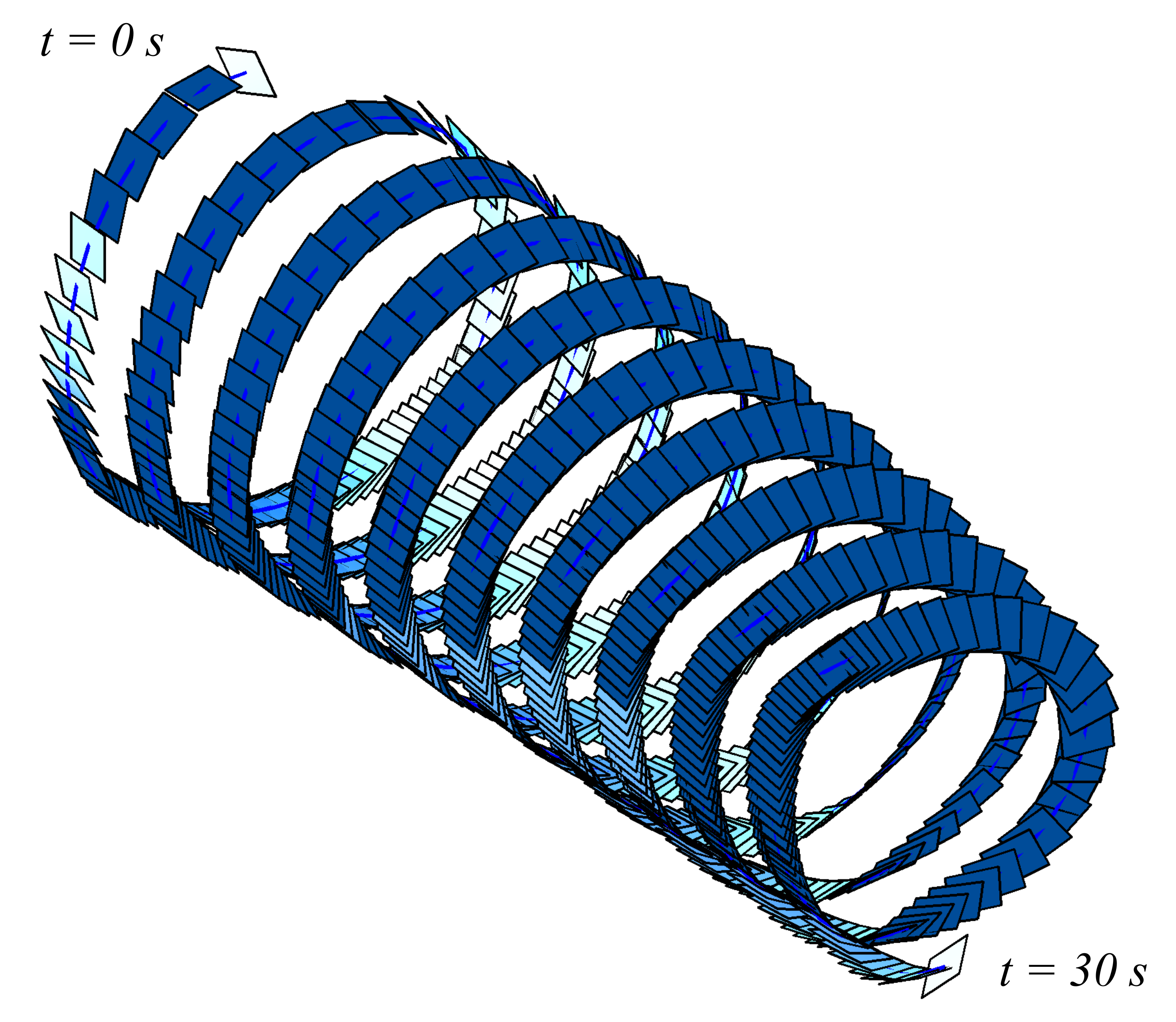} 
\end{tabular}
\caption{\footnotesize  
This figure represents the time evolution $\mathbf{g}  _a =\{g_a ^j, j=1,...,N-1\} $, of a given node $a$ of the beam, when $t^{N-1}=30s$.
The chosen node correspond to $s=0.8\,m$. }
\label{space_evolution10}  
\end{center}
\end{figure}

The algorithm produces the displacement in space $\mathbf{g}_1,..., \mathbf{g}_A$ which, in this example, corresponds to the rotation of a 
beam around an axis combined with a displacement like an air-screw, see Fig.~\ref{space_evolution1} and Fig.~\ref{space_evolution10}.

 
\begin{figure}[H]
\centering
\begin{center}
\begin{tabular}{cc}
\includegraphics[width=1.6 in]{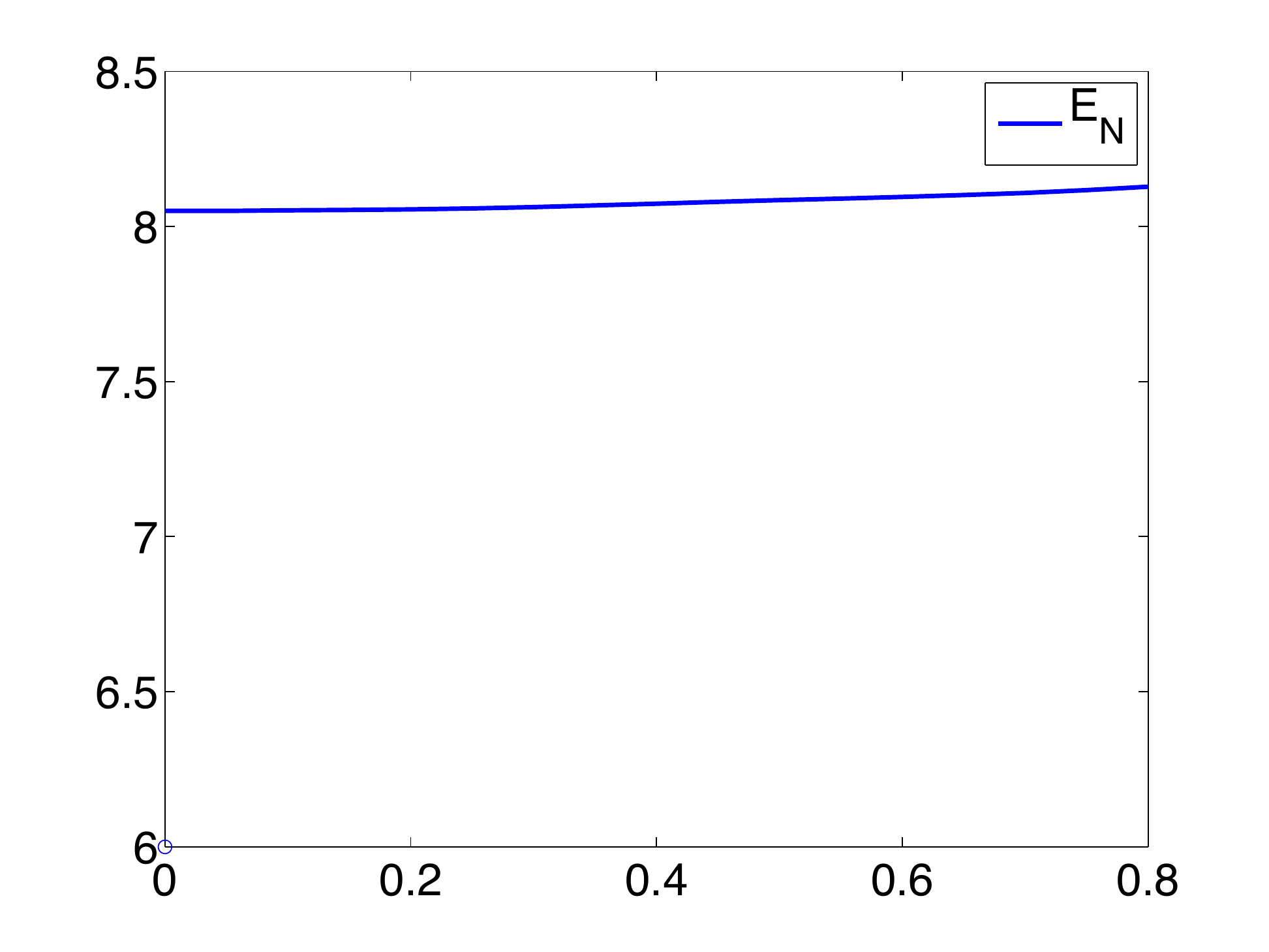}  \qquad   \includegraphics[width=1.6 in]{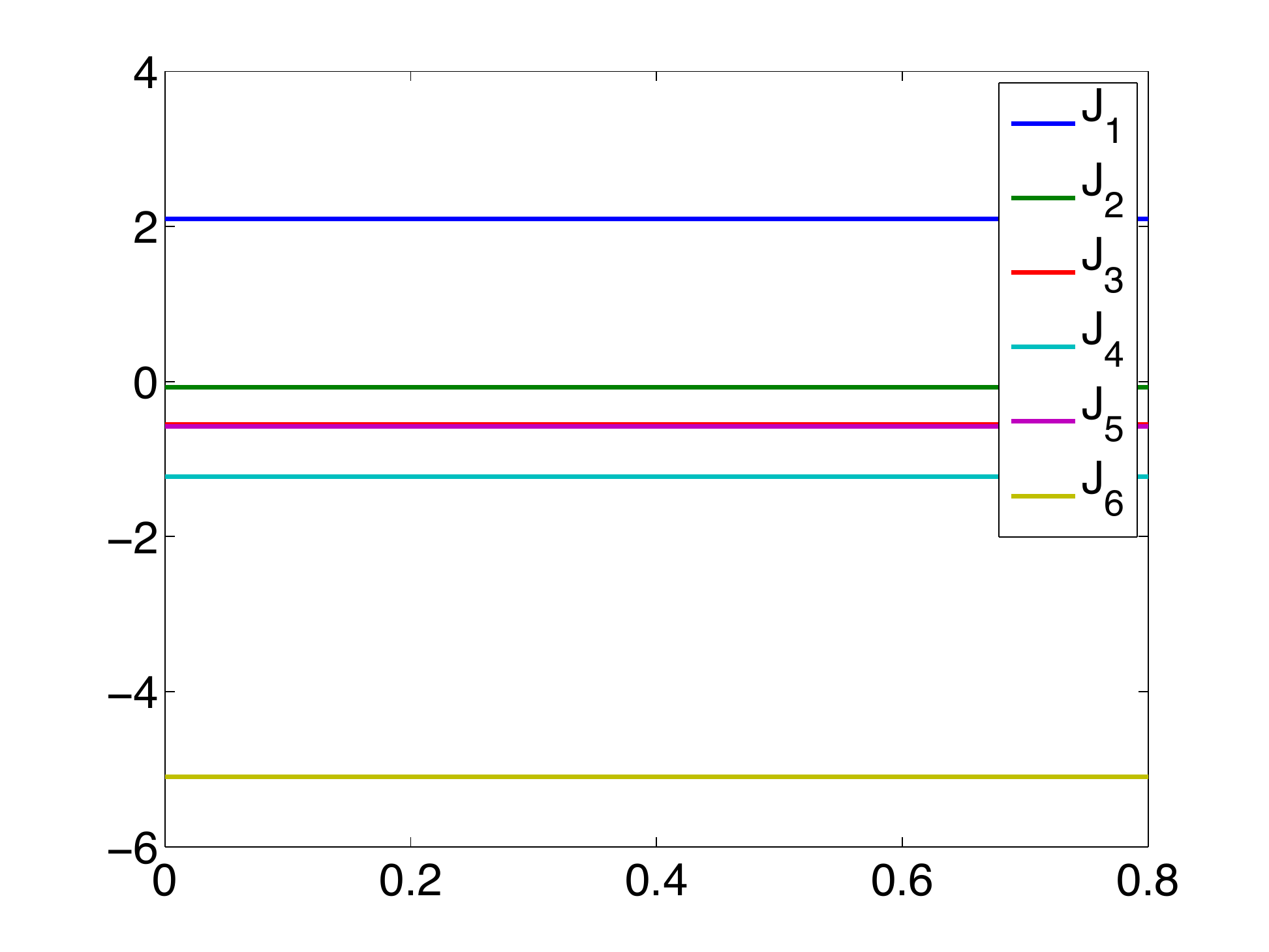}    
\end{tabular}
\vspace{-3pt}
\caption{\footnotesize Left: total ``energy'' behavior $\mathsf{E}_{\mathsf{N}_d}$. Right: conservation of the discrete momentum map $\mathbf J_{\mathsf{N}_d}=(\mathsf{J}^1,...,\mathsf{J}^6) \in \mathbb{R}  ^6 $. 
Both, during a time interval of $0.8$s.}
\label{energy_momentum_space1} 
\end{center}  
\end{figure}

\noindent\textit{Energy behavior.} The above DCEL equations, together 
with the boundary conditions, are equivalent to the DEL equations 
for $\mathsf{N}_d(\mathbf{g}_a , \boldsymbol{\eta}_a )$ in \eqref{model_N_d}; see the discussion in 
\S\ref{subsec:discrete_Lagrangian}. In particular, the solution 
of the discrete scheme defines a discrete symplectic flow in space 
$( \mathbf{g} _a  , \boldsymbol{\eta }_a ) \mapsto 
( \mathbf{g}_{a+1}, \boldsymbol{\eta} _{a+1})$. As a consequence, 
the ``energy'' $\mathsf{E}_{\mathsf{N}_d}$ (see \eqref{formulas_reconstruction_time_space}) of the Lagrangian 
$\mathsf{N}_d$ associated to the spatial evolution description 
is approximately conserved, as illustrated in 
Fig.~\ref{energy_momentum_space1} left, below.

\medskip

\noindent\textit{Momentum map conservation.} Since the discrete 
Lagrangian density is $SE(3)$-invariant, the discrete covariant 
Noether theorem $\mathscr{J}_{B,C}^{K,L}(g _d  ) = 0$ is verified; 
see \S\ref{DCMP}. Since the discrete Lagrangian $\mathsf{N}_d$ is $SE(3)$-invariant, the discrete momentum maps coincide: 
$\mathsf{J}^+_{\mathsf{N}_d}=\mathsf{J}^-_{\mathsf{N}_d}= 
\mathsf{J}_{\mathsf{N}_d}$, and we have
\[
\mathsf{J} _{\mathsf{N}_d}(\mathbf{g}_a, \boldsymbol{\eta}_a)
=\sum_{j=0}^{N-1} \Delta t \mathrm{Ad}^*_{(g_a^j)^{-1}} \lambda_a^j;
\]
see  \eqref{time_momentum_map} and \eqref{discrete_momentum_map}. 
In view of the boundary conditions used here, it follows that the 
discrete momentum map $\mathsf{J} _{\mathsf{N}_d}$ is exactly 
preserved as illustrated in Fig.~\ref{energy_momentum_space1} right. 
This can be seen as a consequence of the covariant discrete Noether 
theorem $\mathscr{J}_{B,C}^{0,N-1}( g _d )=0$. We also checked numerically that the discrete covariant Noether theorem \eqref{DCN} is verified. For example, for $B=K=0$, $C=A-1$, $L=N-1$, we found
\[
\sum_{a=1}^{A-1}\Delta s \left(- \mathrm{Ad}^*_{(g_a^{0})^{-1}} \mu_a^{0} + \mathrm{Ad}^*_{(g_{a}^{N-1})^{-1}} \mu_{a}^{N-1} \right) +\Delta t \left(\mathrm{Ad}^*_{(g_{A-1}^{0})^{-1}} \lambda_{A-1}^{0} - \mathrm{Ad}^*_{(g_0^0)^{-1}} \lambda_0^0\right) = \mathbf{0},
\]
up to round-off error.
 
\begin{remark} \label{boundary_pb} {\rm
We note, in Figures \ref{space_evolution1} and \ref{space_evolution10}, the inappropriate behavior of the trajectory of the section $s=1m$ at time $t=0s$ and $t=10s$. This
problem, which could result in numerical instability for the spatial
algorithm at long distances,
is the subject of future work. This problem does not appear in the
time evolution.}
\end{remark}

\begin{figure}[H]
\centering
\begin{center}
\begin{tabular}{cc}
\includegraphics[width=1.5 in]{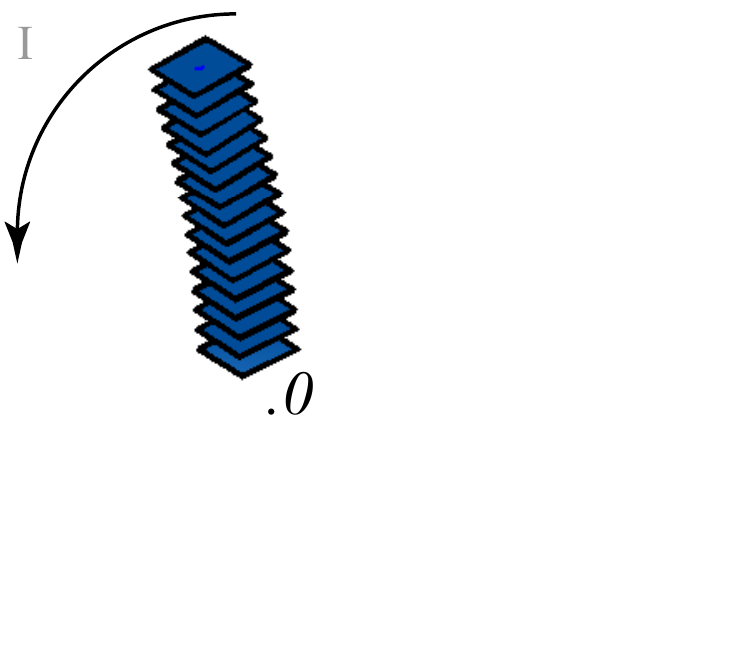}  \includegraphics[width=1.5 in]{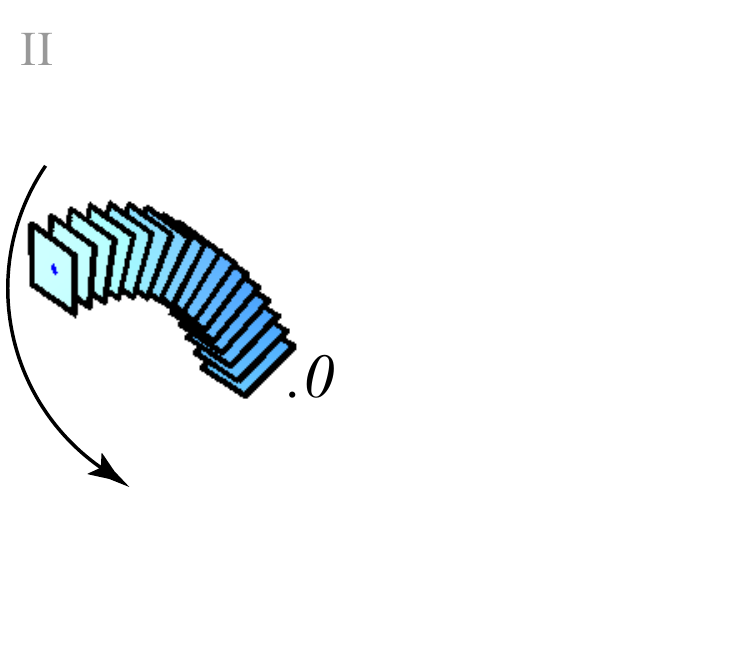}  
\includegraphics[width=1.5 in]{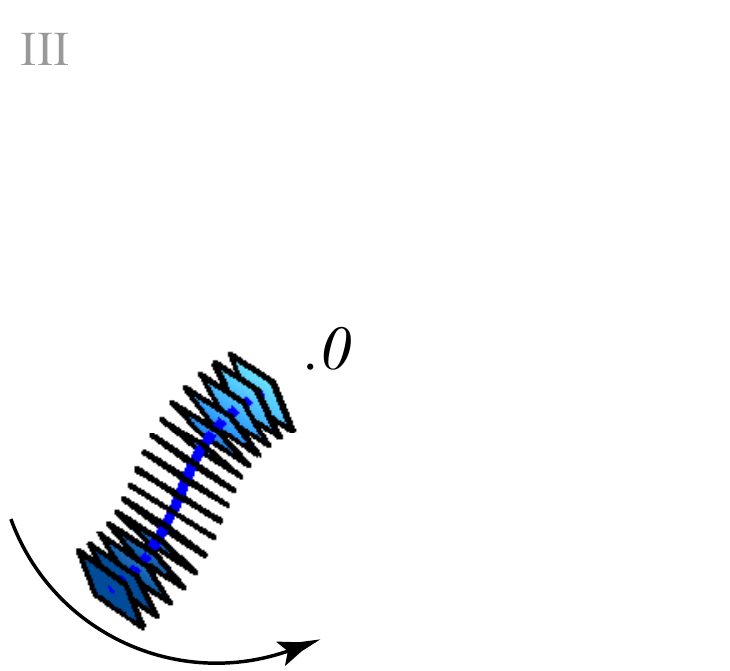}    \includegraphics[width=1.5 in]{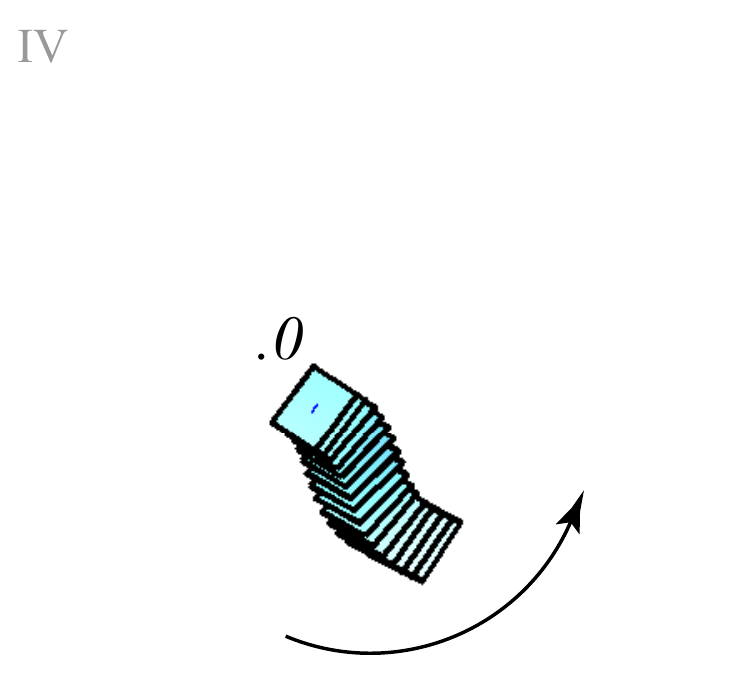}  \\
\includegraphics[width=1.5 in]{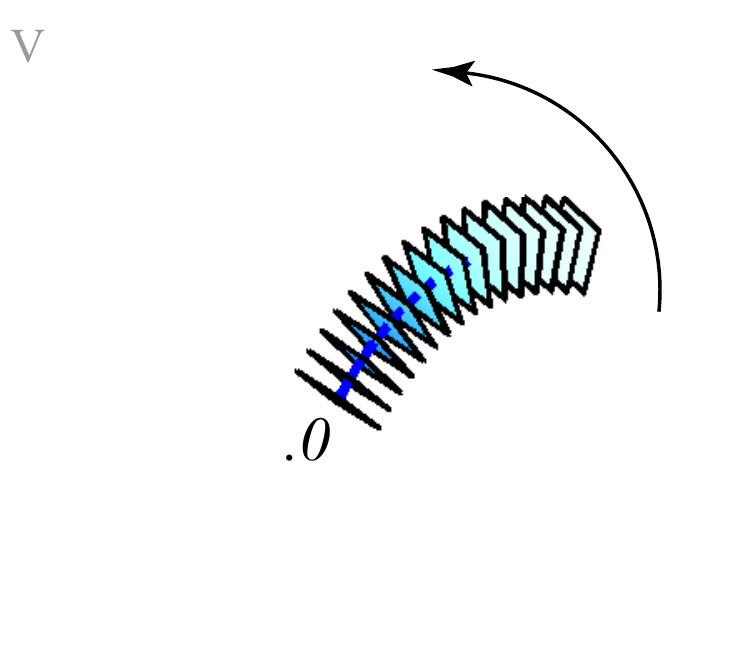} \includegraphics[width=1.5 in]{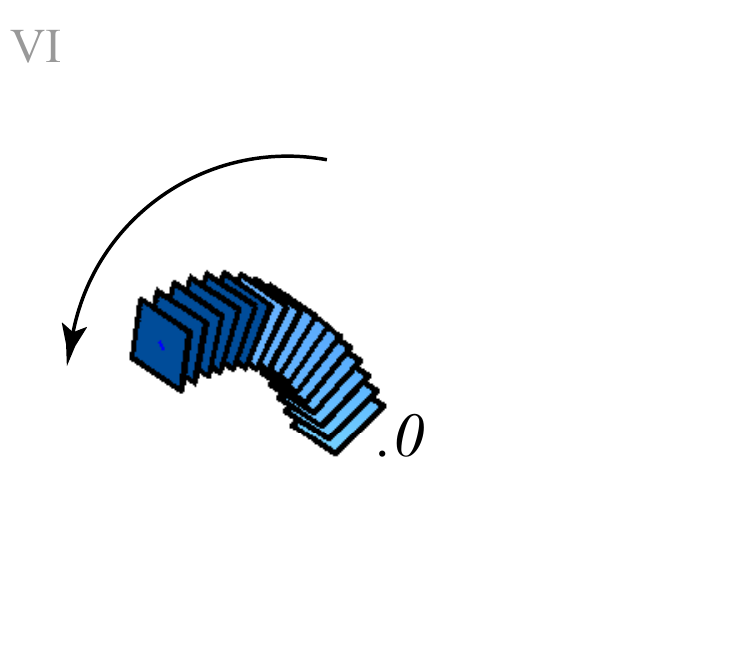} \includegraphics[width=1.5 in]{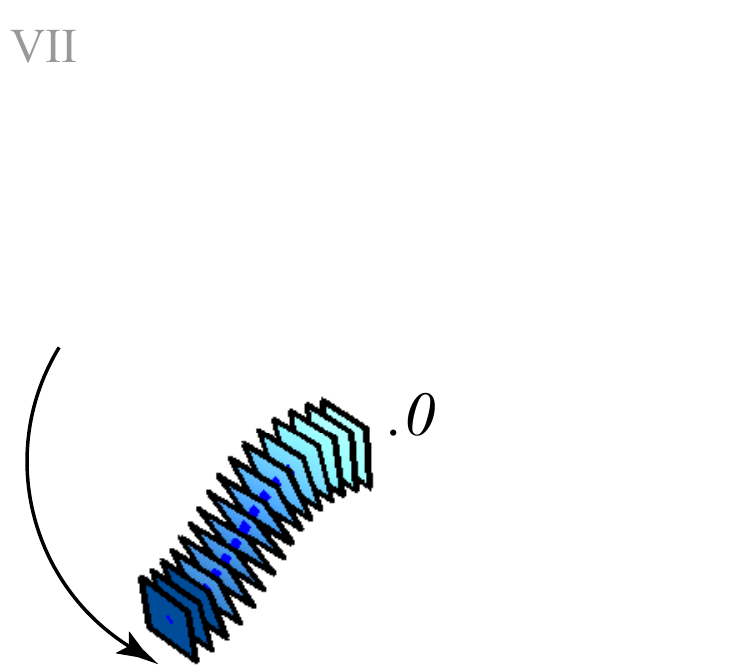} \includegraphics[width=1.5 in]{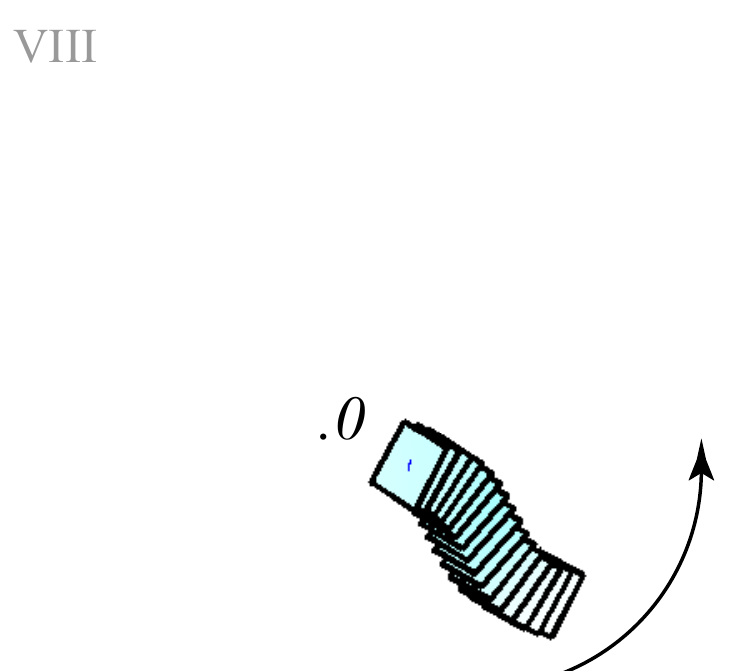}
\end{tabular}
\caption{\footnotesize  
From left to right, and top to bottom: space-integration of the beam sections at times $t=0.1\,s, 0.3\,s, 1\,s, 2\,s,  2.8\,s, 3.45\,s, 4\,s, 5.25\,s$. }\label{reconstruction_space_time1}  
\end{center}
\end{figure}

\paragraph{Time-reconstruction.} Of course, the set $\mathbf{g}_1,..., 
\mathbf{g}_{A-1}$ of time evolutions for each node $a=1,...,A-1$, 
obtained above in Fig.~\ref{space_evolution1}, can be used to 
reconstruct the set $\mathbf{g}^1,..., \mathbf{g}^{N-1}$ of beam 
configurations at each time $j=1,...,N-1$. The resulting motion 
is depicted in Fig.~\ref{reconstruction_space_time1}.
 
We note that the discrete energy $\mathsf{E}_{\mathsf{L}_d}$ and momentum maps $\mathbf{J}_{\mathsf{L}_d}^\pm$ associated to the temporal evolution need not be conserved for this problem, as is already the case in the continuous setting. Their behavior is illustrated in Fig.~\ref{energy_momentum_reconstruction1} below, where we observe periodicity due to rotations.
\begin{figure}[H]
\centering
\begin{center}
\begin{tabular}{cc}
  \includegraphics[width=3 in]{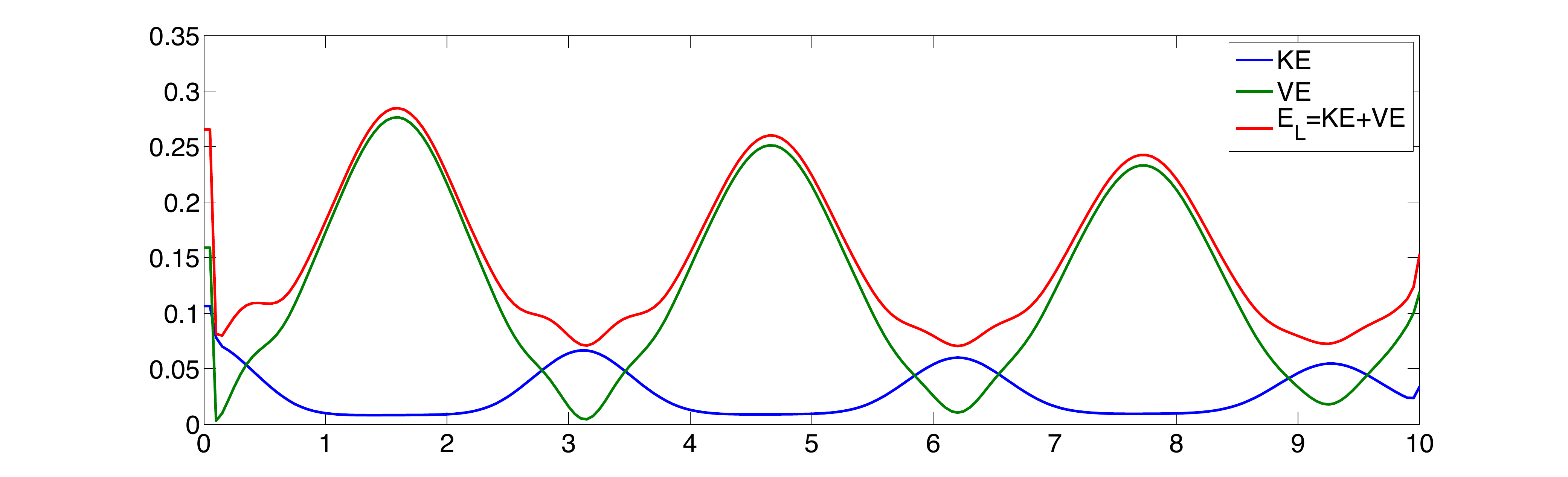}\includegraphics[width=2.95 in]{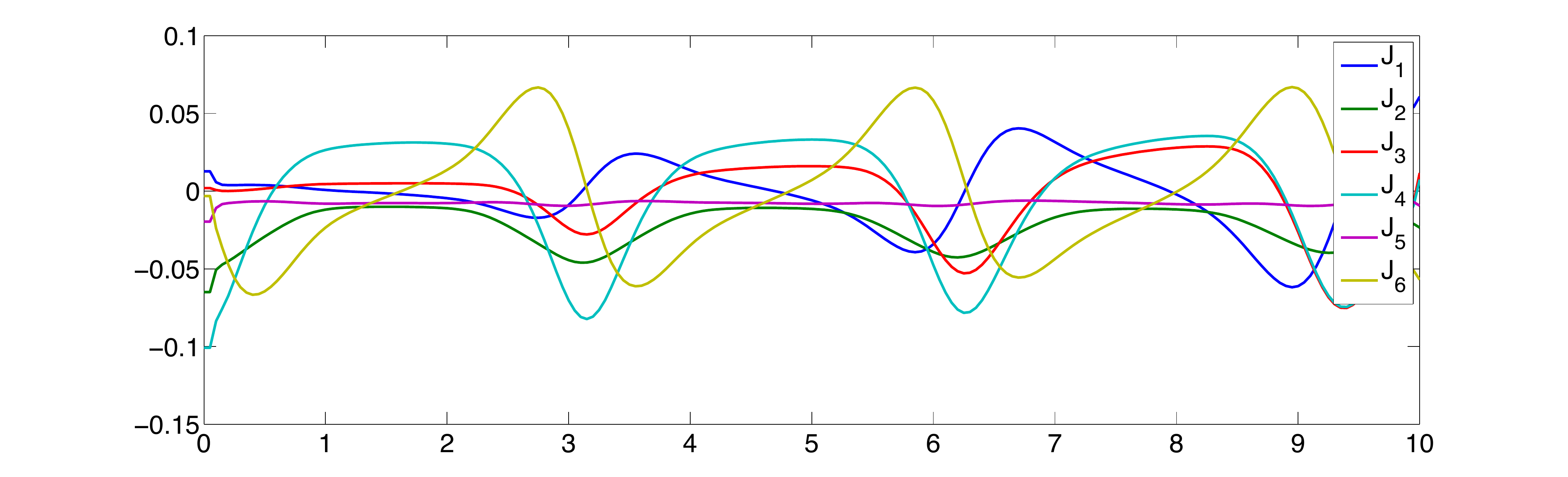}
\end{tabular}
\vspace{-3pt}
\caption{\footnotesize Left: total energy behavior
$\mathsf{E}_{\mathsf{L}_d}$. Right:  momentum map behavior
$\mathbf{J}_{\mathsf{L}_d}=( \mathsf{J} ^1 , ..., \mathsf{J}^6)$. Both, during a time interval of $10$s.}\label{energy_momentum_reconstruction1} 
\end{center}  
\end{figure}

\subsubsection{Scenario B} 

In this example, we employ a finer mesh with $\Delta s=0.02$ and $\Delta t =0.04$. The length of the beam is $L=0.8\,m$ and the total simulation time is $T = 1\,s$. The characteristics of the material are $\rho = 10^3 \,kg/m^3$, $M=10^{-1}\,kg/m$, $E=5.10^4 \,N/m^2$, $\nu=0.35$. 
  
\begin{figure}[H]
  \centering
  \begin{center}
\begin{tabular}{cc}
\includegraphics[width=2 in]{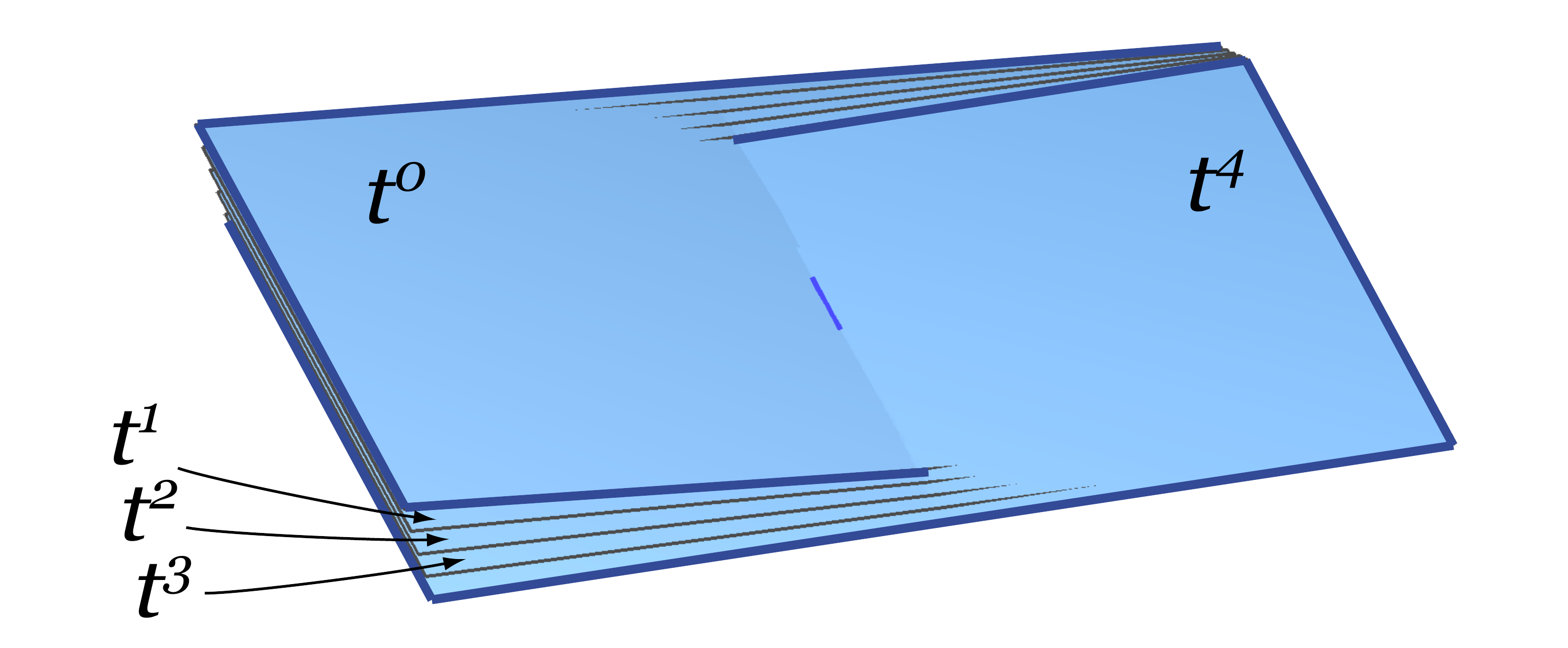} 
\end{tabular}
\caption{\footnotesize  
From left to right: initial conditions $\mathbf{g}_0^j$ (enlarged), when $j\in\{0,1,2,3,4\}$.}\label{initial_condition2}  
\end{center}
\end{figure}

\begin{figure}[H]
\centering
\begin{center}
\begin{tabular}{cc}
\includegraphics[width=1.3 in]{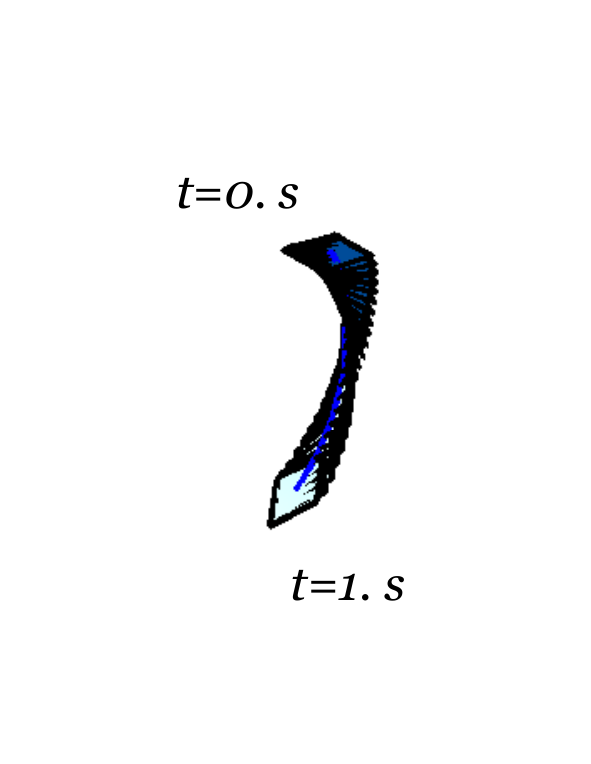}  \includegraphics[width=1.3 in]{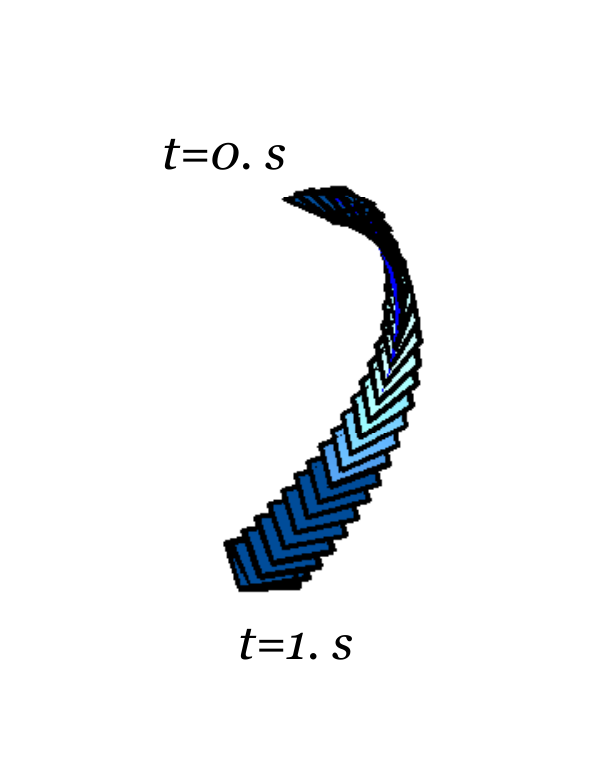}  
\includegraphics[width=1.3 in]{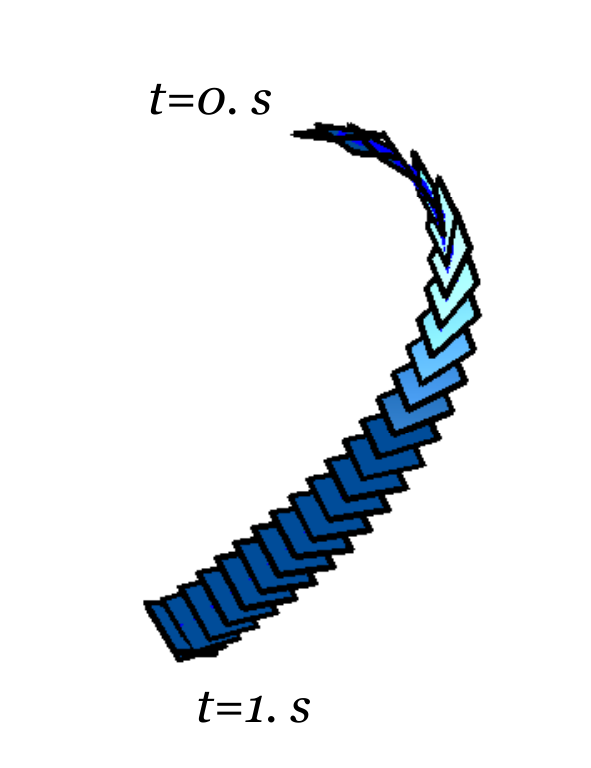}     \includegraphics[width=1.3 in]{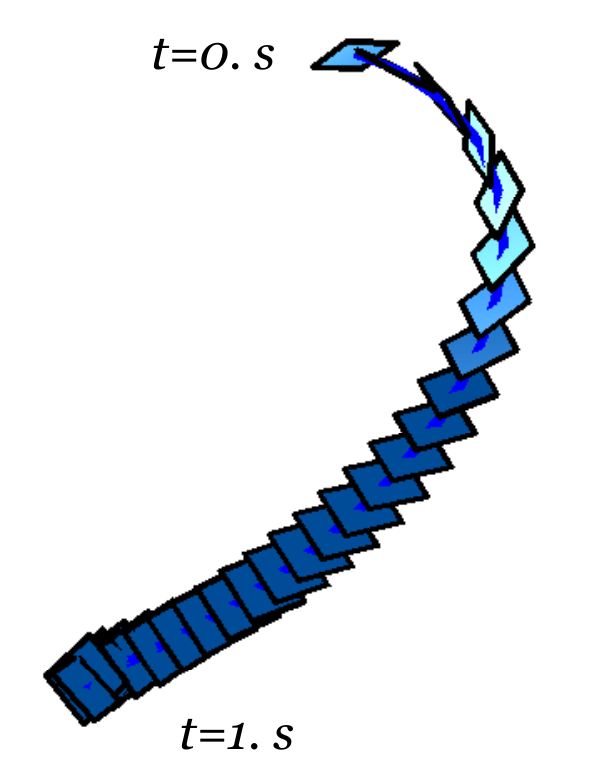}
\end{tabular}
\caption{\footnotesize  
From left to right  : displacement in space within $s=0.2\,m, \,  0.4\,m, \, 0.6\,m, \, 0.8\,m$.}\label{space_evolution2}  
\end{center}
\end{figure}

\paragraph{Space-integration.} The initial conditions $(\mathbf{g}_0,\boldsymbol{\eta}_0)$ are shown on Fig.~\ref{initial_condition2}. In this example we choose the following configuration and strain:
\[
g_0^0=(\mathrm{Id}, (0,0,0)),\quad g_0^{j+1} = g_0^j\, \tau
(\Delta t \xi _0 ^j ),\quad \text{for all $j=0,...,N-1$},
\]
where $ \xi _0 ^j = (0,-0.5,0,0,-0.1,0)$, for all $j=0,...,N-1$, and
\[
\eta_0^j= \frac{1}{\Delta s} \, \tau^{-1}\left((g_0^j)^{-1} g^j_1 \right),\quad \text{for all $j=0,...,N-1$},
\]
where $g^0_1  =(\mathrm{Id}, (0,0,\Delta s))$ and $g^{j+1}_1= g^j_1 \, 
\tau(\Delta t \xi ^j _1 )$, for all $j=0,...,N-1$,  with 
$\xi ^j _1  = (0.06,-0.499, -0.04,-0.03,-0.1,0)$.
As in Scenario A, we do not impose the zero-traction boundary
conditions.
The algorithm  produces the configurations 
$\mathbf{g}_1,..., \mathbf{g}_{A-1}$ depicted in Fig.~\ref{space_evolution2}.

As explained in Scenario A above, the discrete ``energy'' $\mathsf{E}_{\mathsf{N}_d}$ is approximately conserved due to symplecticity in space, and the discrete momentum map $\mathbf{J} _{\mathsf{N}_d}^+(\mathbf{g}_a, \boldsymbol{\eta}_a)$ is exactly preserved. This is illustrated in Fig.~\ref{energy_momentum_space2} below.
\begin{figure}[H]
\centering
\begin{center}
\begin{tabular}{cc}
 \includegraphics[width=1.75 in]{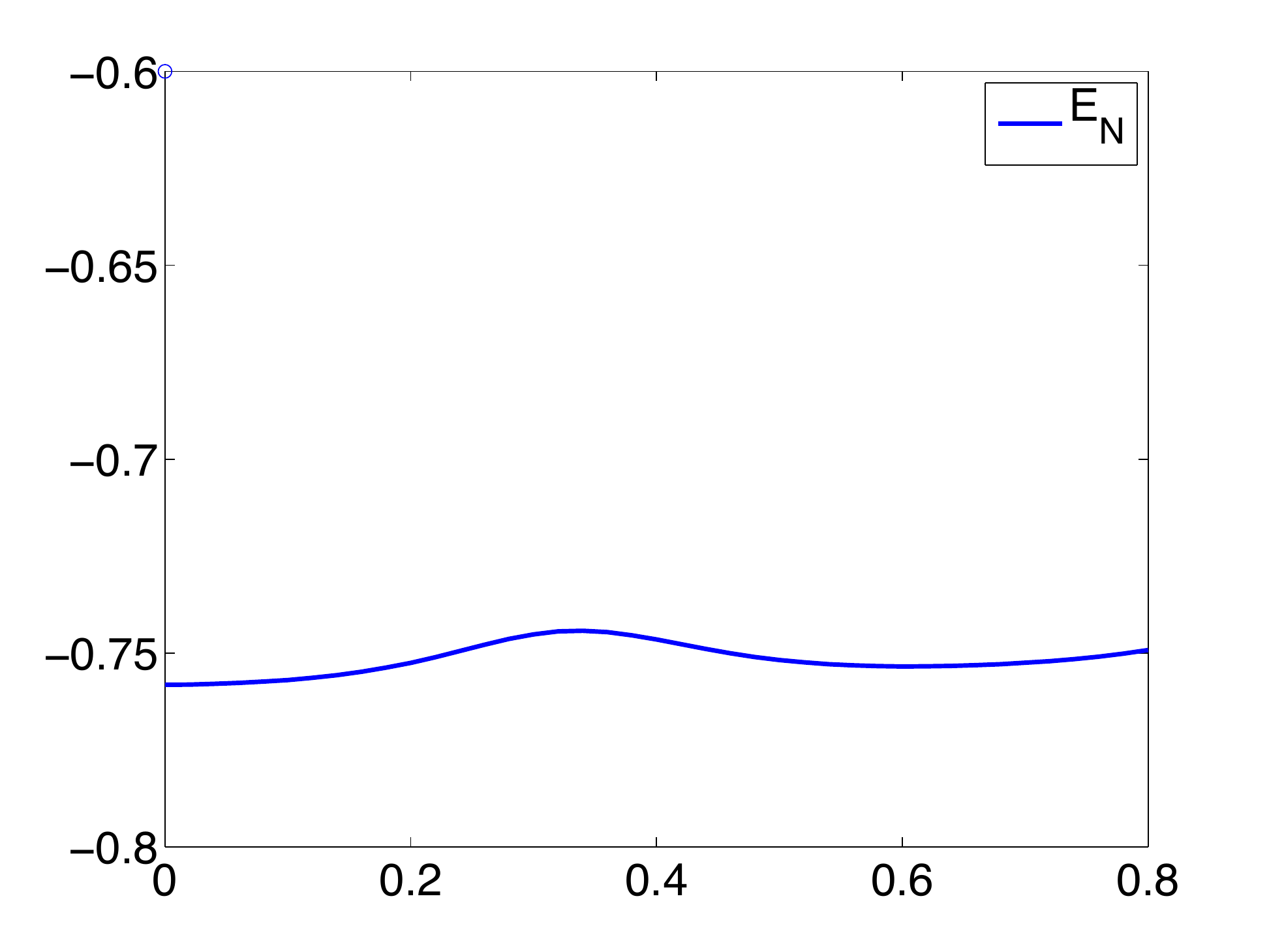}  \qquad   \includegraphics[width=1.75 in]{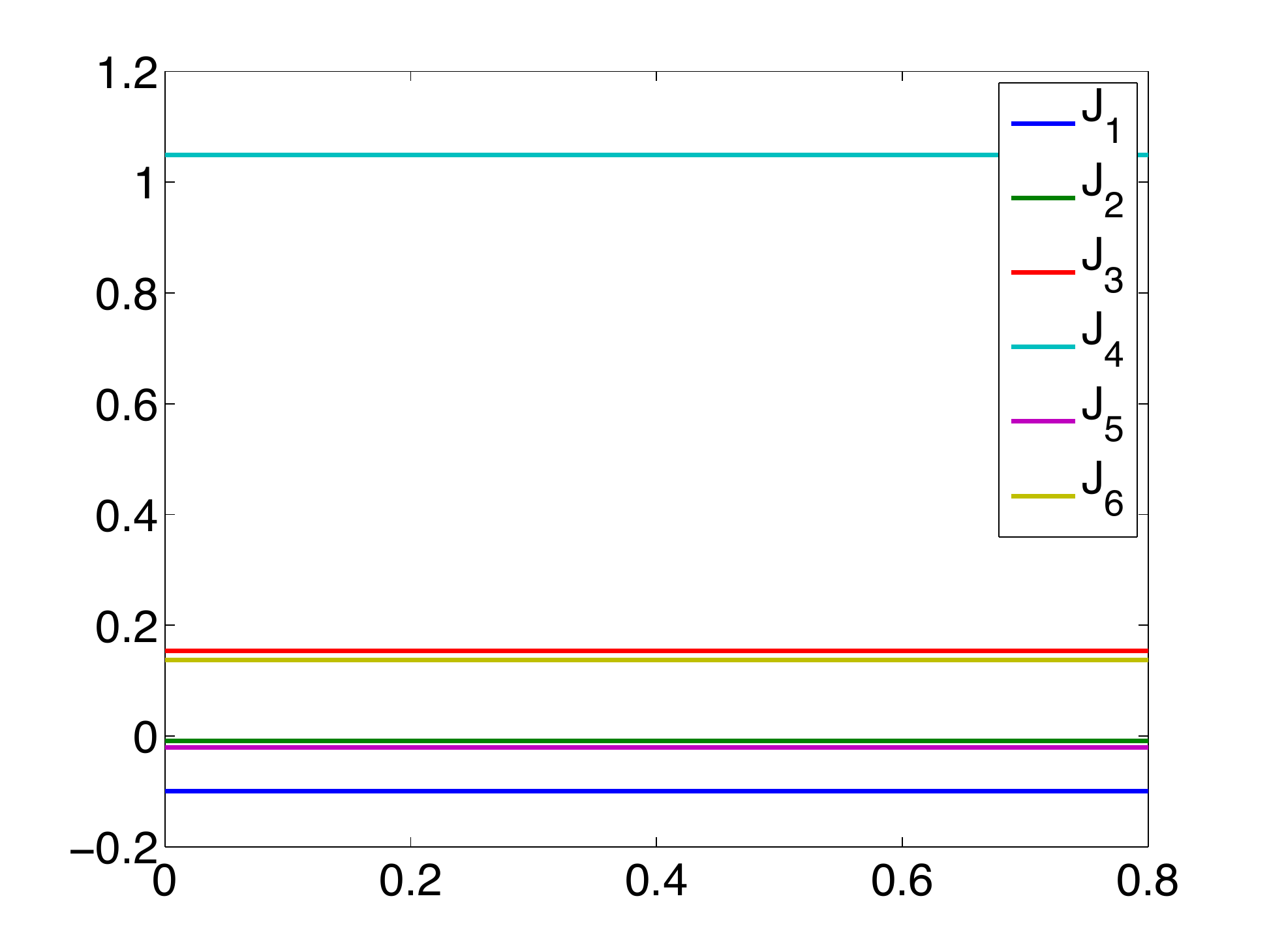}    
\end{tabular}
\vspace{-3pt}
\caption{\footnotesize Left: total ``energy'' behavior 
$\mathsf{E}_{\mathsf{N}_d}$. Right: conservation of the discrete 
momentum map $\mathbf{J}_{\mathsf{N}_d}=(\mathsf{J}^1,...,\mathsf{J}^6)$.}\label{energy_momentum_space2} 
\end{center}  
\end{figure}
As in Scenario A, we checked numerically that the covariant Noether theorem is verified. For example, we have
\[
\sum_{a=1}^{A-1}\Delta s \left(- \mathrm{Ad}^*_{(g_a^{0})^{-1}} \mu_a^{0} + \mathrm{Ad}^*_{(g_{a}^{N-1})^{-1}} \mu_{a}^{N-1} \right) +\Delta t \left(\mathrm{Ad}^*_{(g_{A-1}^{0})^{-1}} \lambda_{A-1}^{0} - \mathrm{Ad}^*_{(g_0^0)^{-1}} \lambda_0^0\right) = \mathbf{0},
\]
up to round-off error.

\begin{figure}[H]
\centering
\begin{center}
\begin{tabular}{cc}
\includegraphics[width=1.1 in]{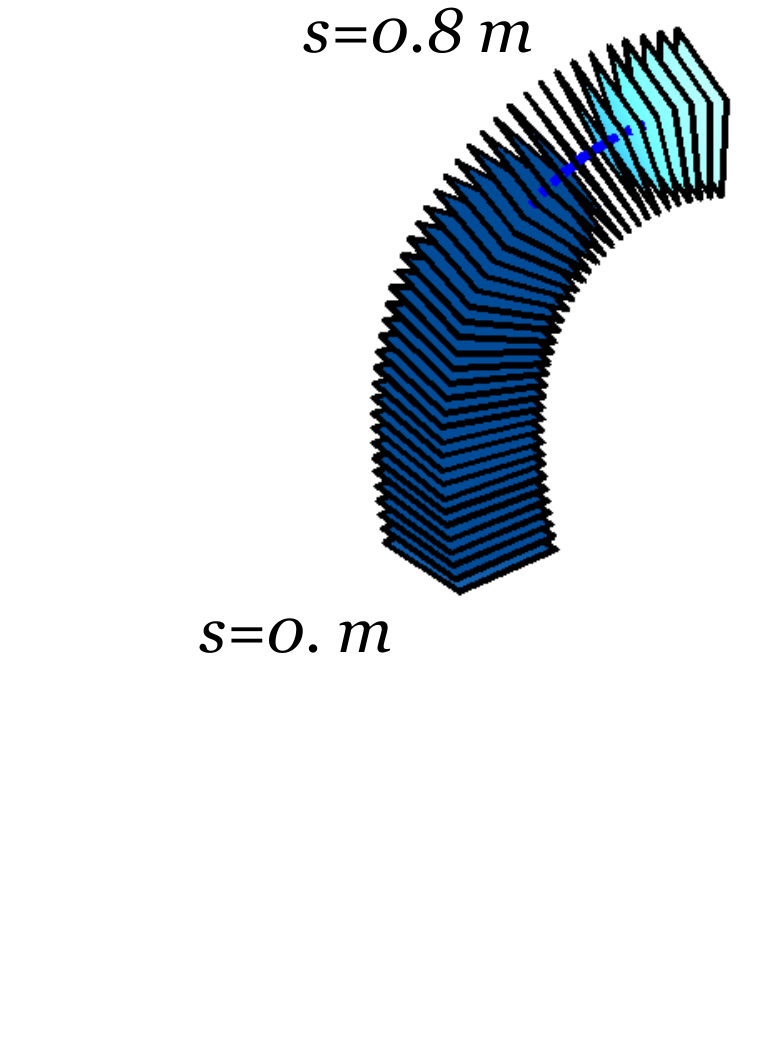}  \includegraphics[width=1.1 in]{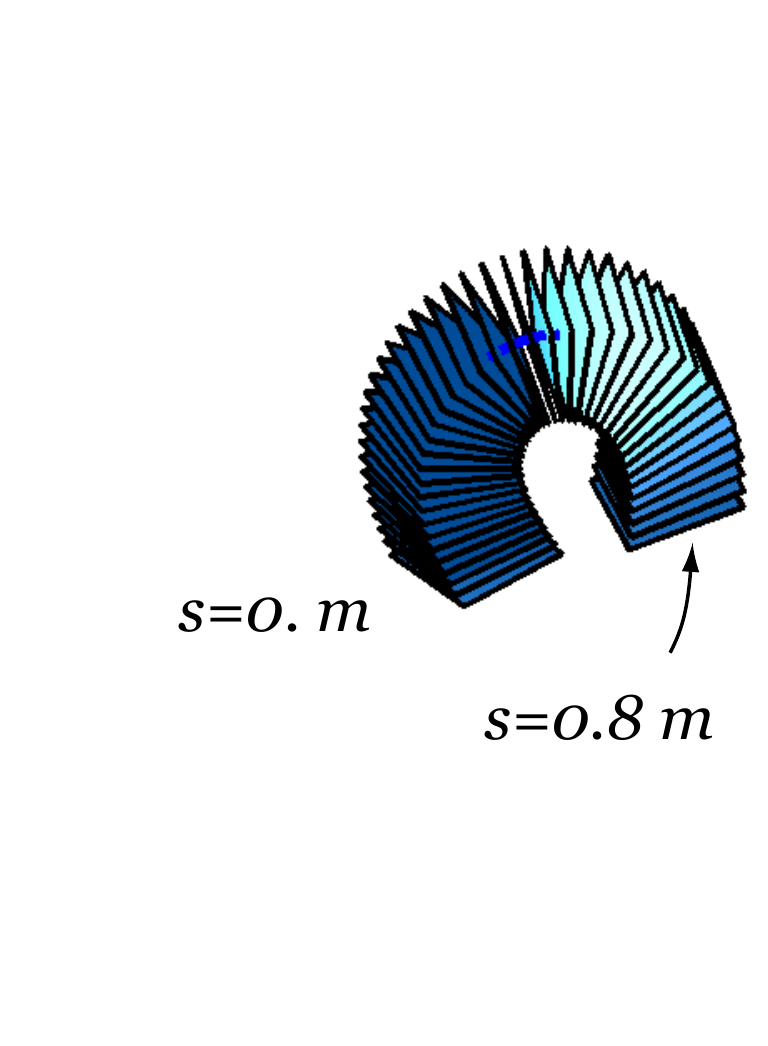}  
\includegraphics[width=1.1 in]{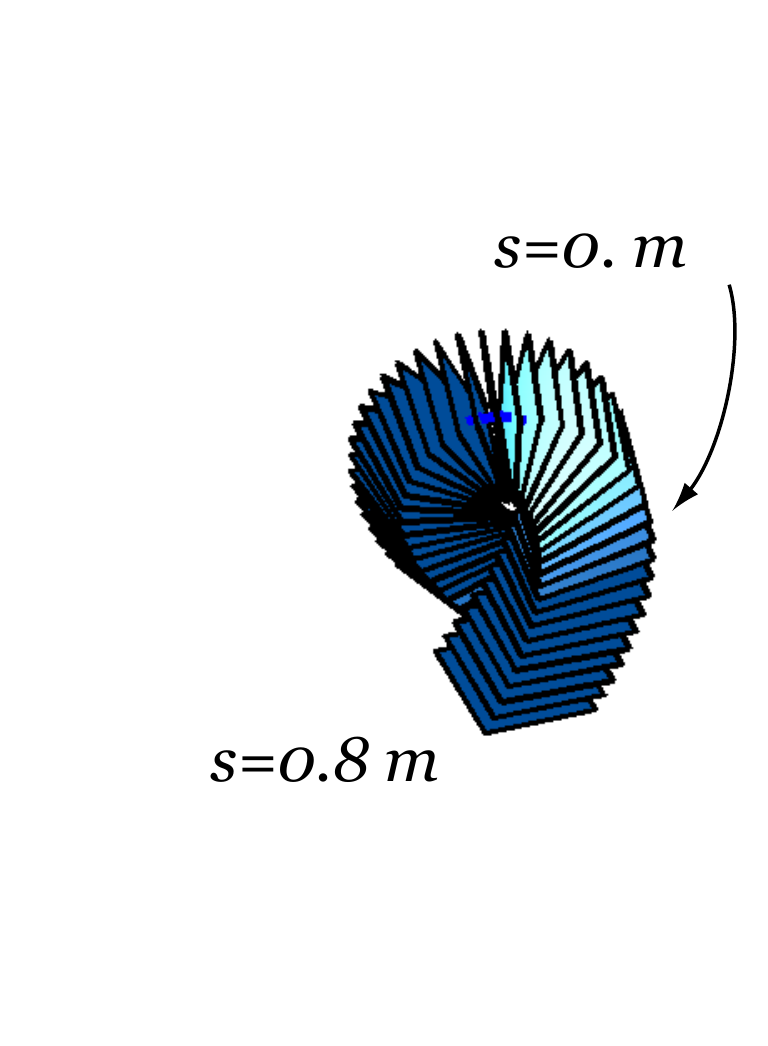}    \includegraphics[width=1.1 in]{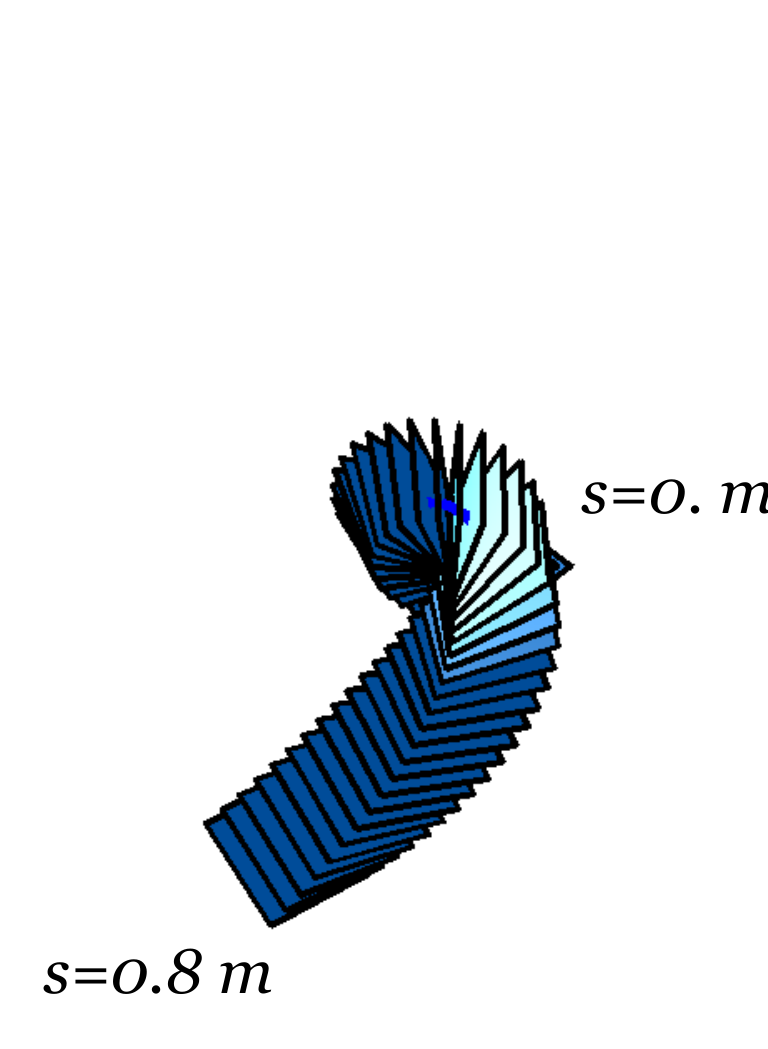} \quad  \includegraphics[width=1.1 in]{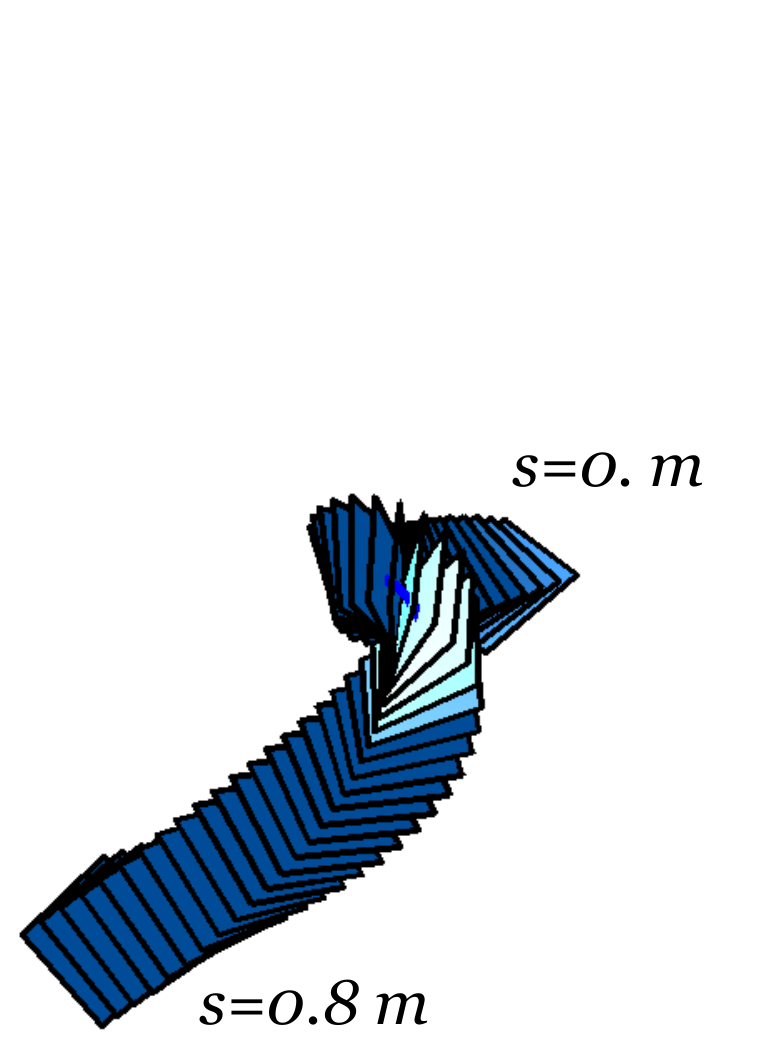}
\end{tabular}
\caption{\footnotesize  
From left to right: reconstruction of the trajectories in space of the sections, at times $t=0.16\,s, 0.36\,s, 0.52\,s, 0.76\,s, 1\,s$. }\label{reconstruction_space_time2}  
\end{center}
\end{figure}

\paragraph{Time-reconstruction.} The above computed space evolution provides a set of configurations $\mathbf{g}_0,..., \mathbf{g}_{A-1}$
for the length $s = 0.8\, m$. Repackaging this data yields 
a set of time configurations $\mathbf{g}^1,..., \mathbf{g}^{N-1}$  for the duration of $1\,s$, where 
$\mathbf{g}^j= (g_0^j,...,g_{A-1}^j)$ (see 
Figure~\ref{mesh_recons}). The obtained spatial trajectories of the sections are depicted in Figure \ref{reconstruction_space_time2}.

\begin{remark}
\label{last_remark}
{\rm
The time step is set to a fraction of the Courant limit CFL \cite{CoFrLe1928}, and computed as
$\Delta t = \frac{d}{10\,c}$,
where $d$ is the radius of the largest ball contained in the mesh element and $c$ is the nominal dilational wave speed of the material (function of the Young modulus and Poisson ratio). In our example, $d = \Delta s$, and $c=\sqrt{ \frac{\lambda +2 \mu}{\rho}}$, where $\lambda, \mu$ are the Lam\'e parameters, and $\rho$ is the density of the material. So, if the space step $\Delta s$ is reduced, then the time step $\Delta t$ is reduced. For homogeneous and isotropic materials, the Lam\'e parameters are proportional to the Young modulus. So, if the Young modulus increases, then the time step $\Delta t$ decreases.
}
\end{remark}

\section{Conclusion}

In this paper, we have introduced a discrete spacetime multisymplectic variational integrator counterpart of the continuous covariant beam formulation. We verified, through numerical examples, that the symplectic integrators in time and in space preserve the momentum maps, and that the energy oscillates around its nominal value.

We showed that the tests presented in this work validate
the general theory developed in \cite{DeGBRa2013}. In particular, 
the global discrete Noether theorem is always verified and, depending on the boundary conditions used, it implies a discrete classical Noether theorem in space or in time.

\medskip

We point out some unresolved issues that will be addressed in future work
concerning multisymplectic integrators in order to get an even more
accurate numerical tool.
\begin{itemize} 
\item[{\rm (i)}]
Solve the inappropriate behavior of the boundaries at time $t=0$ and 
$t=N$ when the integrator is updated in space, as noted in Remark
\ref{boundary_pb}.
\item[{\rm (ii)}]
We noted that the integration algorithm performs better in space 
than in time. Find similar necessary conditions for the integration 
in space, as explained in Remark \ref{last_remark}.   
\end{itemize}

\section*{Appendix}\label{Appendix} 

In this appendix we quickly explain how to obtain the explicit expressions of the equations \eqref{mu_lambda_formula}--\eqref{boun_cond3} for the case $G=SE(3)$. Recall that in this case we write $g_a^j=(\Lambda_a^j,\mathbf{r}_a^j) \in SE(3)$ and $\xi_a^j=(\boldsymbol{\omega}_a^j, \gamma_a^j)^T, \; \eta_a^j=(\boldsymbol{\Omega}_a^j, \Gamma_a^j)^T\in \mathfrak{se}(3)$, and we choose the approximation of the exponential map given by the map $ \tau : \mathfrak{se}(3) \rightarrow SE(3)$ in \eqref{cay}.
Using the formula for $\left(\left({\rm d}\tau_{(\omega , \boldsymbol{\gamma})}\right)^{-1}\right)^*$ derived in  \S\ref{Cayey_section}, we obtain that the discrete momenta $\mu_a^j, \lambda_a^j$ in \eqref{mu_lambda_formula} read explicitly
\begin{align*}
\mu_a^j& =  \begin{bmatrix}
\mathbf{I}_3 - \frac{1}{2} \Delta t \,\omega_a^j +
\frac{1}{4} \Delta t^2 \, \boldsymbol{\omega}_a^j (\boldsymbol{\omega}_a^j) ^T & \mathbf{0}_3 \\[6pt]
-\frac{1}{2}\left(\mathbf{I}_3 -
\frac{1}{2} \Delta t\,\omega_a^j \right) \Delta t\, \gamma_a^j & \mathbf{I}_3 -\frac{1}{2} \Delta t\, \omega_a^j \end{bmatrix} ^T  \mathbb{J} \left[ \begin{array}{c} \boldsymbol{\omega}_a^j \\ \gamma_a^j \end{array}\right], 
\\
\lambda_a^j &=  \begin{bmatrix}
\mathbf{I}_3 - \frac{1}{2}\Delta s\,\Omega_a^j +
\frac{1}{4}\Delta s^2\, \boldsymbol{\Omega}_a^j (\boldsymbol{\Omega}_a^j) ^T & \mathbf{0}_3 \\[6pt]
-\frac{1}{2}\left(\mathbf{I}_3 -
\frac{1}{2}\Delta s\,\Omega_a^j \right)\Delta s\, \Gamma_a^j & \mathbf{I}_3 -
\frac{1}{2}\Delta s\, \Omega_a^j \end{bmatrix} ^T  \mathbb{C} \left[ \begin{array}{c} \boldsymbol{\Omega}_a^j \\ \Gamma_a^j \end{array}\right],
\end{align*}
where 
\begin{align*}
\left[ \begin{array}{c} \boldsymbol{\omega}_a^j \\ \gamma_a^j \end{array}\right] &= \frac{1}{\Delta t} \tau^{-1} \left( (\Lambda_a^j)^{-1}\Lambda_a^{j+1}, (\Lambda_a^j)^{-1}(\mathbf{r}_a^{j+1}- \mathbf{r}_a^j) \right), 
\\
\left[ \begin{array}{c} \boldsymbol{\Omega}_a^j \\ \Gamma_a^j \end{array}\right] &= \frac{1}{\Delta s} \tau^{-1} \left( (\Lambda_a^j)^{-1}\Lambda_{a+1}^{j}, (\Lambda_a^j)^{-1}(\mathbf{r}_{a+1}^{j}- \mathbf{r}_a^j) \right).
\end{align*}
Equations \eqref{CDEL_beam}--\eqref{boun_cond3} can be explicitly written for $ \mathfrak{g}  = \mathfrak{se}(3)$, by making use of the formulas
\begin{align*}
\mathrm{Ad}^*_{\tau(\Delta t \xi_a^{j})}  \mu_a^j &=  \begin{bmatrix}
\mathbf{I}_3 + \frac{1}{2} \Delta t \,\omega_a^j +
\frac{1}{4} \Delta t^2 \, \boldsymbol{\omega}_a^j (\boldsymbol{\omega}_a^j) ^T & \mathbf{0}_3 \\[6pt]
 \frac{1}{2}\left(\mathbf{I}_3 +
 \frac{1}{2} \Delta t\,\omega_a^j \right) \Delta t\, \gamma_a^j & \mathbf{I}_3 +\frac{1}{2} \Delta t\, \omega_a^j \end{bmatrix} ^T  \mathbb{J} \left[ \begin{array}{c} \boldsymbol{\omega}_a^j \\ \gamma_a^j \end{array}\right],
\\
\mathrm{Ad}^*_{\tau(\Delta s \eta_{a}^j)} \lambda_a^j & =  \begin{bmatrix}
\mathbf{I}_3 + \frac{1}{2}\Delta s\,\Omega_a^j +\frac{1}{4}\Delta s^2\, \boldsymbol{\Omega}_a^j (\boldsymbol{\Omega}_a^j) ^T & \mathbf{0}_3 \\[6pt]
\frac{1}{2}\left(\mathbf{I}_3 +\frac{1}{2}\Delta s\,\Omega_a^j \right)\Delta s\, \Gamma_a^j & \mathbf{I}_3 +\frac{1}{2}\Delta s\, \Omega_a^j \end{bmatrix} ^T  \mathbb{C} \left[ \begin{array}{c} \boldsymbol{\Omega}_a^j \\ \Gamma_a^j \end{array}\right],
\end{align*}
and
\[
(g_a^j)^{-1} \Pi_{g_a^j}(g_a^j)= \left[ \begin{array}{c} 0 \\ (\Lambda_a^j)^{-1}\mathbf{q}_a^j \end{array}\right],
\] 
for $\Pi(\Lambda_a^j, \mathbf{r}_a^j)=\langle \mathbf{q}_a, \mathbf{r}_a^j\rangle$.

{\footnotesize

\bibliographystyle{new}

}

\end{document}